\documentclass[11pt]{article}
\usepackage[T1]{fontenc}
\usepackage{lmodern,amsmath,amssymb,amsthm,mathtools}
\usepackage[margin=1in]{geometry}
\usepackage{microtype}
\usepackage[hidelinks]{hyperref}
\hypersetup{
 pdftitle={Critical exponents and spectral functions of totally nonnegative matrices},
 pdfauthor={Wei Xie},
 pdfkeywords={total positivity, oscillatory matrix, matrix function, critical exponent}
}
\newtheorem{theorem}{Theorem}[section]
\newtheorem{lemma}[theorem]{Lemma}
\newtheorem{corollary}[theorem]{Corollary}
\newtheorem{proposition}[theorem]{Proposition}
\newtheorem{conjecture}{Conjecture}
\theoremstyle{definition}
\newtheorem{definition}[theorem]{Definition}
\theoremstyle{remark}
\newtheorem{remark}{Remark}[section]
\newtheorem*{example}{Example}
\DeclareMathOperator{\adj}{adj}
\DeclareMathOperator{\diag}{diag}
\newcommand{\TN}{\mathrm{TN}}
\newcommand{\TP}{\mathrm{TP}}
\newcommand{\Osc}{\mathrm{Osc}}
\newcommand{\dd}{\,\mathrm d}
\numberwithin{equation}{section}
\title{Critical exponents and spectral functions of totally nonnegative matrices}
\author{Wei Xie}
\date{26 November 2026}

\begin{document}
\maketitle

\begin{abstract}
The ordinary critical exponent of the totally positive matrices of
order $n\ge2$ is $n-2$, attaining the lower bound of Johnson and
Walch. For every oscillatory matrix $A$, the power $A^t$ is totally
positive for all $t>n-2$. On nonsingular totally nonnegative inputs,
separate cuts at zeros immediately above and below the diagonal give
a sharp uniform threshold $d-2$ for total nonnegativity when the
resulting blocks have order at most $d$, $2\le d\le n$. Above this
threshold the cuts determine all zero minors. For $n\ge4$, the exponents
$t$ for which $A^t$ is a $P$-matrix for every nonsingular totally
nonnegative $A$ of order $n$ form $\mathbb Z\cup\{t:|t|\ge n-3\}$.
In these dimensions, only four corner minors can fail to be positive
on totally positive inputs when $n-3<t<n-2$. All $25$ admissible joint sign
patterns occur on every prescribed simple positive spectrum, and exactly
five occur on symmetric inputs of that spectrum. The positivity estimates
extend to $A^\alpha e^{\psi(A)}$ for real operator-monotone $\psi$.
In fixed order $n$, the positive analytic preservers of total
nonnegativity that fail strict positivity on some oscillatory input
are exactly $cx^k$ with $c>0$ and $k=0,\ldots,n-2$. In all orders,
the positive functions preserving total nonnegativity on symmetric
totally positive inputs are exactly the monomials $cx^k$, $c>0$,
$k\in\mathbb Z_{\ge0}$.
\end{abstract}


\section{Introduction}

A real matrix is totally nonnegative if all its square minors are nonnegative, and totally positive if all those minors are positive. These properties will be abbreviated as TN and TP. We write $\TN_n$ and $\TP_n$ for the corresponding classes of matrices of order $n$. The Cauchy--Binet formula shows that ordinary multiplication preserves total nonnegativity and total positivity for square matrices. Positive integer powers therefore preserve both classes. For real exponents, we ask which powers preserve TN or TP matrices of a fixed order.

An $n\times n$ TN matrix $A$ is oscillatory if some positive integer power is TP. The classical oscillation theorem states that $A$ has distinct positive eigenvalues and that $A^{n-1}$ is TP \cite[Theorems 2.6.6 and 5.2.1]{FJ}. Thus an oscillatory matrix has an ordinary real power
\[
 A^t=X\diag(\lambda_1^t,\ldots,\lambda_n^t)X^{-1},
 \qquad \lambda_i>0,
\]
for every real $t$.

For a matrix class $\mathcal C$ on which real powers are defined, write $g(\mathcal C)$ for the infimum of the numbers $\gamma$ such that $A^t\in\mathcal C$ for every $A\in\mathcal C$ and every $t>\gamma$. This is the ordinary critical exponent of $\mathcal C$. Johnson and Walch \cite[Theorems 3.8--3.9]{JW} proved that the ordinary critical exponent of $\TP_n$ is finite and satisfies
\[
 n-2\le g(\TP_n)<\infty.
\]
Their lower bound uses irreducible positive definite tridiagonal matrices and approximation by TP matrices. Fallat suggested $n-2$ as a possible general value when an ordinary critical exponent exists \cite[p.~3088]{F17}. The question is whether $n-2$ also controls every minor order simultaneously.

Guillot, Khare, and Rajaratnam \cite[Corollary 2.3]{GKR} proved
that the ordinary critical exponent of doubly nonnegative matrices
of order $n$ is $n-2$, using the divided-difference criterion of
Micchelli and Willoughby \cite[Corollary~3.1]{MW79}. A doubly
nonnegative matrix is symmetric positive semidefinite and has
nonnegative entries. For a symmetric nonsingular TN matrix $A$,
the $p$th compound $C_p(A)$ is doubly nonnegative of order
$\binom np$. Its entries are the $p$th-order minors of $A$, and
$C_p(A^t)=C_p(A)^t$. Thus the doubly nonnegative result guarantees
nonnegativity of these minors for $t\ge\binom np-2$, $1\le p<n$.
The resulting bounds depend on the compound orders $\binom np$,
which can be much larger than $n$.

Han, Johnson, and Paparella \cite[Corollary~4.5]{HJP} studied
ordinary powers in the larger class of diagonalizable entrywise
nonnegative matrices with nonnegative spectrum. In order three the
critical exponent of this class is $2$. The compounds of a TP
matrix have positive entries. Fallat, Gupta, and Johnson
\cite[proof of Theorem~A]{FGJ25} proved that for every TP matrix $A$, each
$C_p(A)$ with $2\le p\le n-2$ has a negative $3\times3$ minor.

The question for oscillatory matrices includes inputs with zero minors.
A zero immediately above or below the main diagonal of a nonsingular TN
matrix imposes a triangular block structure on its powers. If the two
patterns of zeros differ, which blocks govern the critical exponent
and which minors can become positive?

There are two natural ways to weaken the demand on the output. One
can prescribe the order of a minor and require one of its index sets
to be the first or last indices. One can also require only the principal
minors to be positive, the defining condition for a $P$-matrix.
For integer powers, Fallat and Liu \cite{FL07} characterized the
oscillatory matrices attaining the classical bound $n-1$. Zarai
and Margaliot \cite{ZM21} found exact exponents for classes described
by bidiagonal factorizations.

We study the exponents for which ordinary real powers preserve
total positivity or total nonnegativity, the effect of adjacent
zero entries, and positivity when the order or the index sets of
a minor are prescribed. This includes the case of principal
minors. For spectral functions we consider $x^\alpha e^{\psi(x)}$
with $\psi$ operator monotone, and ask when a positive analytic TN
preserver in a fixed order maps oscillatory matrices to TP
matrices. We also ask which positive functions preserve TN on
symmetric TP matrices of every order.

\section{Matrix functions and the main results}\label{mainsec}

Index sets are always listed in increasing order. For a matrix $A$,
we write $A[I,J]$ for the submatrix on rows $I$ and columns $J$,
and use $a:b=\{a,a+1,\ldots,b\}$ for a consecutive index set.
Complements are taken relative to the full index set of the matrix
in question; for a single index $i$, we write $i^c=\{i\}^c$.
The symbol $\Osc_n$ denotes the class of oscillatory matrices of order $n$.

Every nonsingular TN matrix has positive spectrum, although it need not be diagonalizable. For such a matrix we use the principal power
\[
 A^t=\exp(t\log A).
\]
Equivalently, on a Jordan block $\lambda I+N$ of size $d$, with $\lambda>0$ and $N^d=0$, its power is
\[
 (\lambda I+N)^t=\sum_{j=0}^{d-1}\binom tj\lambda^{t-j}N^j.
\]
For a function $f$ analytic near the positive axis, the same Jordan
decomposition defines $f(A)$ by the finite Taylor expansion
$f(\lambda I+N)=\sum_{j=0}^{d-1}f^{(j)}(\lambda)N^j/j!$.
These matrix functions depend continuously on $A$ while its spectrum
remains in their domain of analyticity.

\subsection{Total nonnegativity and ordinary powers}

\begin{theorem}\label{main}
Let $n\ge2$ and $A\in\Osc_n$. Then $A^t$ is TP for every real $t>n-2$. The threshold cannot be decreased, and $t=n-2$ cannot be included for all oscillatory matrices.
\end{theorem}

In particular, $g(\TP_n)=n-2$, attaining the lower bound of Johnson--Walch \cite{JW}.
The proof expresses the far corner minors as positive sums of
evaluation determinants built from boundary eigenvectors. Comparing
these determinants with those for the integer power $A^{n-1}$ gives
the sign for every real $t>n-2$.
If $A$ is TP and $n\ge3$, then $A^{n-2}$ is TP by
Cauchy--Binet. If the oscillatory matrix $A$ has upper bandwidth one,
then $(A^{n-2})_{1n}=0$, since each step in a product of entries
can increase the index by at most one.

\begin{corollary}\label{sets}
For $n\ge2$, the real exponents $t$ for which $A^t$ is TN for every nonsingular $A\in\TN_n$ form exactly
\[
 \mathbb Z_{\ge0}\cup[n-2,\infty).
\]
For $n\ge3$, the exponents preserving $\TP_n$ form exactly $\mathbb Z_{>0}\cup[n-2,\infty)$; for $n=2$ they form $(0,\infty)$. The exponents making $A^t$ TP for every $A\in\Osc_n$ form exactly $(n-2,\infty)$.
\end{corollary}

For singular TN matrices with $q\ge2$ positive eigenvalues counted
with algebraic multiplicity and largest Jordan block of size $\nu$ at
zero, the sharp threshold for total nonnegativity is $q+\nu-2$
(Proposition~\ref{singularpowers}).

The following theorem extends the positivity estimate to spectral multipliers.
A real function $\psi$ on $(0,\infty)$ is operator monotone if, for positive definite Hermitian matrices $B,C$ of every order, positivity of $C-B$ implies positivity of $\psi(C)-\psi(B)$ in the positive semidefinite sense. Its analytic continuation and integral representation are recalled in Lemma~\ref{Pickrep}.

\begin{samepage}
\begin{theorem}\label{multipliers}
Let $n\ge2$, let $\psi$ be real and operator monotone on $(0,\infty)$,
and let $\Phi(x)=x^\alpha e^{\psi(x)}$ with $\alpha\ge n-2$.
Then $\Phi(A)=A^\alpha\exp(\psi(A))$ is TN for every nonsingular
$A\in\TN_n$. If $A\in\Osc_n$, it is TP whenever $\alpha>n-2$
or $\psi$ is nonconstant. If $\Phi$ extends analytically across
zero, then $\Phi(A)$ is TN also for singular $A\in\TN_n$.
\end{theorem}
\end{samepage}

\subsection{Zero minors and smaller blocks}

\begin{definition}\label{directionblocks}
For a nonsingular TN matrix $A$ of order $n$, let $\mathcal U(A)$
be the partition of $1:n$ obtained by cutting between $k$ and
$k+1$ whenever $a_{k,k+1}=0$. Define $\mathcal L(A)$ in the same
way using $a_{k+1,k}=0$, and put
\[
 d(A)=\max\{|C|:C\in\mathcal U(A)\cup\mathcal L(A)\}.
\]
\end{definition}

The two partitions can differ. An oscillatory matrix has $d(A)=n$;
the case $d(A)=1$ consists of positive diagonal matrices.

\begin{theorem}\label{zeropattern}\label{localmother}
Let $A=(a_{ij})\in\TN_n$ be nonsingular, $n\ge2$, and put
$d=d(A)$. Let $\Phi(x)=x^\alpha e^{\psi(x)}$, where $\psi$ is
real and operator monotone on $(0,\infty)$.
For $d\ge2$, suppose that $\alpha\ge d-2$. Then $\Phi(A)$ is TN.
If $\alpha>d-2$ or $\psi$ is nonconstant, then, for increasing
index sets $I=\{i_1<\cdots<i_p\}$ and $J=\{j_1<\cdots<j_p\}$,
$1\le p\le n$,
\[
 \det\Phi(A)[I,J]>0
 \quad\Longleftrightarrow\quad
 \det A^{d-1}[I,J]>0
 \quad\Longleftrightarrow\quad
 \det A^{n-1}[I,J]>0.
\]
These determinants are positive exactly when, for each $r$,
\begin{equation}\label{zeropath}
 \begin{cases}
 a_{k,k+1}>0&\text{for every }i_r\le k<j_r,\quad\text{if }i_r<j_r,\\
 a_{k+1,k}>0&\text{for every }j_r\le k<i_r,\quad\text{if }i_r>j_r.
 \end{cases}
\end{equation}
There is no condition when $i_r=j_r$. If the criterion fails,
all three determinants are zero.
For $d=1$, every positive analytic function of $A$ is positive
diagonal, and the same zero-minor criterion holds.
\end{theorem}

Since $d(A)\le n$, this includes the strict range in
Theorem~\ref{multipliers}. The zero positions depend only on the
zeros immediately above and below the diagonal of $A$, also when
eigenvalues coincide. In particular, $\Phi(A)$ is TP in the
strict range if and only if $A$ is oscillatory.

\begin{corollary}\label{localpowers}
Fix $n\ge2$ and $2\le d\le n$. The real exponents $t$ for which
$A^t$ is TN for every nonsingular $A\in\TN_n$ with $d(A)\le d$
form exactly
\[
 \mathbb Z_{\ge0}\cup[d-2,\infty).
\]
The exponents for which, for every such $A$, the matrix $A^t$ is
TN and has exactly the same zero minors as $A^{n-1}$ form exactly
$(d-2,\infty)$. If $d=1$, both sets of exponents are
$\mathbb R$.
\end{corollary}

\subsection{Four boundary minors and P-matrix preservation}

\begin{theorem}\label{iminternal}
Let $n\ge4$ and
\[
 \Phi(x)=x^\alpha e^{\psi(x)},\qquad \alpha\ge n-3,
\]
where $\psi$ is real and operator monotone on $(0,\infty)$.
For every nonsingular $A\in\TN_n$, the matrix $F=\Phi(A)$
is a $P$-matrix. If $A\in\TP_n$, every nonempty square minor
of $F$ is strictly positive except possibly
\begin{equation}\label{fourcorners}
 F_{1n},\quad F_{n1},\quad
 \det F[1:n-1,2:n],\quad \det F[2:n,1:n-1].
\end{equation}
For such $A$, the matrix $F$ is TN if and only if the four
quantities in \eqref{fourcorners} are nonnegative, and TP if
and only if they are positive.
As a bound valid for every operator-monotone $\psi$, the condition
$\alpha\ge n-3$ cannot be decreased.
\end{theorem}

For a nonsingular TN matrix $A$, approximation by TP matrices
shows that every minor of $\Phi(A)$ outside \eqref{fourcorners}
is nonnegative. Thus $\Phi(A)$ is TN if and only if the four
quantities in \eqref{fourcorners} are nonnegative.
Theorem~\ref{fourpatterns} classifies the joint signs of the four
exceptional minors on TP inputs of any fixed positive simple
spectrum. For $n-3<t<n-2$, all $25$ admissible patterns occur on
each such spectrum (Corollary~\ref{fournecessary}).

\begin{theorem}\label{imPexact}\label{imPtail}
Let $T_n$ be the set of real exponents $t$ for which $A^t$ is a
$P$-matrix for every nonsingular $A\in\TN_n$. Then
\[
 T_n=\begin{cases}
 \mathbb R,&1\le n\le3,\\[1mm]
 \mathbb Z\cup\{t\in\mathbb R:|t|\ge n-3\},&n\ge4.
 \end{cases}
\]
For $n\ge4$, let $H_n=\{1,n\}$. Then
\[
 T_n=\{t\in\mathbb R:\det A^t[H_n,H_n]>0
                         \text{ for every }A\in\TP_n\}.
\]
For every nonempty compact subset $K$ of one component of
$\mathbb R\setminus T_n$, there is a single $A\in\TP_n$
such that $\det A^t[H_n,H_n]<0$ for all $t\in K$.
\end{theorem}

The same sets result if $A$ is restricted to $\TP_n$ or $\Osc_n$.
In particular, $T_4=(-\infty,-1]\cup\{0\}\cup[1,\infty)$.

For $1\le p<n$, define
\[
 \begin{aligned}
 \mathcal E_{n,p}=\bigl\{t>0:\ &\det A^t[I,J]>0
       \text{ for every }A\in\TP_n\\
       &\text{and all }I,J\subseteq1:n\text{ with }|I|=|J|=p\bigr\},
 \end{aligned}
\]
and put
\[
 g_{n,p}=\inf\{\gamma\ge0:(\gamma,\infty)\subseteq
                                  \mathcal E_{n,p}\}.
\]

\begin{corollary}\label{imexact}
Let $n\ge4$ and $A\in\TP_n$.  Every minor of $A^t$ of orders
$2,\ldots,n-2$ is positive when $t\ge n-3$.  Moreover,
\[
 \mathcal E_{n,2}=\mathcal E_{n,n-2}
 =\{1,\ldots,n-4\}\cup[n-3,\infty),\qquad
 g_{n,2}=g_{n,n-2}=n-3.
\]
For $n=4$, the finite set preceding the interval is empty.
\end{corollary}

\subsection{Strictness for arbitrary preservers}

\begin{theorem}\label{strict}
Fix $n\ge2$. Let $f:(0,\infty)\to(0,\infty)$ be real analytic, and suppose that $f(A)$ is TN for every nonsingular $A\in\TN_n$. Then $f(A)$ is TP for every $A\in\Osc_n$ if and only if
\[
 f(x)\not\equiv cx^k\quad\text{for all }c>0\text{ and }k\in\{0,\ldots,n-2\}.
\]
\end{theorem}

The matrix order is fixed in this assertion; the function need not
preserve TN matrices of every order or belong to the operator-monotone
family above. For $n\ge3$, under the same preservation hypothesis,
the strictness condition is equivalent to
$f'(1)>(n-2)f(1)$, by Proposition~\ref{rigidity}.
For functions satisfying this condition,
Proposition~\ref{stuniform} gives a quantitative form: on compact
sets of simple positive spectra, the far corners of $f(A)$ are
uniformly comparable to those of $A^{n-1}$.
Corollary~\ref{stzeropattern} determines all zero minors of $f(A)$
under a spectral condition on the directional blocks. For $n\ge4$,
only blocks of order at least four need have simple spectrum.

\subsection{Corner minors and standard facts}

For Theorems~\ref{main} and \ref{multipliers}, we use the lower-left and upper-right corner minors
\begin{equation}\label{corners}
 F_p(A;t)=\det A^t[n-p+1:n,1:p],\qquad
 G_p(A;t)=\det A^t[1:p,n-p+1:n],
 \quad 1\le p<n.
\end{equation}
We also write $F_p(X)$ and $G_p(X)$ for the corresponding corner minors of a single matrix $X$, so that $F_p(A;t)=F_p(A^t)$ and $G_p(A;t)=G_p(A^t)$.
To prove that a TN matrix is TP, we use the following classical
recognition theorem.

\begin{lemma}[{\cite[Theorem 3.1.10]{FJ}}]\label{recognition}
A square TN matrix is TP if and only if all its lower-left and upper-right corner minors, including its determinant, are positive.
\end{lemma}

Every TN matrix, including a singular one, is a limit of TP matrices
\cite{FJ}. A nonsingular TN matrix $H$ is a $P$-matrix. Indeed, the TN LU
factorization \cite[Theorem~2.10]{GL} writes $H=LU$, with $L$
unit lower triangular TN and $U$ upper triangular TN with positive
diagonal. Cauchy--Binet gives, for every nonempty index set $I$,
\[
 \det H[I,I]\ge\det L[I,I]\det U[I,I]
             =\prod_{i\in I}U_{ii}>0.
\]

Equivalently, the factorization can be written $H=LDU$, where
$L,U$ are unit lower and upper triangular TN and $D$ is positive
diagonal. This form is unique.

We also use standard oscillation facts \cite[Theorems 5.2.1, 5.3.3, and 5.5.2]{FJ}: an oscillatory matrix has distinct positive eigenvalues; its left and right eigenvectors have nonzero first and last components; and deletion of the first or last row and column gives strict eigenvalue interlacing. In particular, if $x_a,y_a$ are dual right and left eigenvectors, then
\begin{equation}\label{projection}
 (x_a)_1(y_a)_1
 =\frac{\det(\lambda_a I-A[2:n,2:n])}
        {\left.\dfrac{\mathrm d}{\mathrm dz}\det(zI-A)\right|_{z=\lambda_a}}>0.
\end{equation}
The sign follows directly from strict interlacing. The corresponding statement holds at the last component. For a rectangular TN matrix $T$, variation diminution gives $s^-(Tx)\le s^-(x)$, where $s^-$ counts sign changes after zeros are deleted \cite[Theorem 4.2.2]{FJ}.

Sections~\ref{spectralsec}--\ref{powersec} give a direct proof of
Theorem~\ref{main}. Section~\ref{coeffsec} proves the
operator-monotone multiplier theorem, the principal-block estimates,
and the zero-minor and local-reduction results.
Section~\ref{principalSection} treats principal minors and real powers.
Sections~\ref{cssection} and~\ref{imsection} prove positivity of
the internal-order minors, and Section~\ref{endpointsec} treats
endpoint flag minors. The arbitrary preserving functions of
Theorem~\ref{strict} are studied in Section~\ref{stsection};
the final section gives fixed multipliers and finite spectral interpolation.

\section{Boundary minors of the eigenvector matrix}\label{spectralsec}

The corner minors of a spectral function will be expanded using rows of
adjugate matrices. We divide these rows by their first components and
express the resulting entries as rational functions with negative real
poles. The signs of their initial evaluation minors give the positive
coefficients in the corner expansion.

For increasing nodes $z_1<\cdots<z_k$, put $\Delta(z)=\prod_{i<j}(z_j-z_i)$. Empty products and empty determinants are $1$.

\begin{lemma}\label{blocks}
Let $B$ be TP of order $n\ge2$, and suppose that the rectangular matrix $[B;v]$ obtained by appending the row $v$ is TP. Write
\begin{equation}\label{blockdef}
 B=\begin{pmatrix}b&r\\c&D\end{pmatrix},\quad
 v=v_1(1,u),\quad C=cu-D,\quad d=bu-r,\quad h=d+uC.
\end{equation}
Then $C$, $[u;C]$, and $[d;C]$ are TP. If
\[
 R(z)=v\adj(zI-B),\qquad q(z)=R_1(z),
\]
then
\begin{equation}\label{ratio}
 q(z)=v_1\det(zI+C),\qquad
 \frac{R_{j+1}(z)}{q(z)}=u_j-h(zI+C)^{-1}e_j.
\end{equation}
\end{lemma}

\begin{proof}
Row operations give the following identities for tail index sets of the appropriate sizes:
\begin{align*}
 \det C[I,J]&=\frac1{v_1}
 \det\begin{pmatrix}c[I]&D[I,J]\\v_1&v_1u[J]\end{pmatrix},\\
 \det\begin{pmatrix}u[J]\\C[I,J]\end{pmatrix}
 &=\frac1{v_1}\det\begin{pmatrix}D[I,J]\\v_1u[J]\end{pmatrix},\\
 \det\begin{pmatrix}d[J]\\C[I,J]\end{pmatrix}
 &=\frac1{v_1}\det\begin{pmatrix}b&r[J]\\c[I]&D[I,J]\\v_1&v_1u[J]\end{pmatrix}.
\end{align*}
Subtract multiples of the last row in each determinant on the right. The signs from factoring out $-1$ and moving the last row cancel. Every determinant on the right is a minor of $[B;v]$.

The first identity in \eqref{ratio} follows by expanding the bordered determinant for $v\adj(zI-B)e_1$. To obtain the second, use $R(z)(zI-B)=\det(zI-B)v$. Subtracting its first component times $u$ from its last $n-1$ components gives
\[
 R_{2:n}(z)(zI+C)=q(z)(zu-d).
\]
Since $zu-d=u(zI+C)-h$, the formula follows.
\end{proof}

\begin{lemma}\label{spectral}
Let $B\in\TP_n$. For every row $v$ such that $[B;v]$ is TP, put
$R(z)=v\adj(zI-B)$ and $q(z)=R_1(z)$. Then $q(z)>0$ for $z>0$ and
\[
 \det\left[\frac{R_j(z_i)}{q(z_i)}\right]_{i,j=1}^p>0
 \qquad(1\le p\le n,\quad 0<z_1<\cdots<z_p).
\]
For $n\ge2$, the coefficient matrix $W$ in the partial fraction
expansion \eqref{residues} has positive initial-column flag minors.

Let $0<\lambda_1<\cdots<\lambda_n$ be the eigenvalues of $B$ and let
$\ell_i$ be the corresponding left eigenvectors, normalized by
$(\ell_i)_1=1$. Taking $z_a=\lambda_{i_a}$ gives
\begin{equation}\label{prefix}
 \det[(\ell_{i_a})_j]_{a,j=1}^p>0
 \qquad(1\le p\le n,\quad 1\le i_1<\cdots<i_p\le n).
\end{equation}
This spectral inequality also follows from the eigenvector factorization of Dopico and Koev \cite[Theorems~2 and~3]{DopicoKoev}.
\end{lemma}

\begin{proof}
We use induction on $n$. For $n=1$ one has $q(z)=v_1>0$ and $R/q=1$.
Suppose $n\ge2$ and fix any $v$ such that $[B;v]$ is TP.
Use Lemma~\ref{blocks}, and put $m=n-1$. Let $0<\mu_1<\cdots<\mu_m$ be the eigenvalues of $C$. Choose dual eigenvectors $x_a,y_a$ of $C$ with $(x_a)_1=1$. Equation~\eqref{projection} gives $(y_a)_1>0$. If $[w;C]$ is TN, variation diminution applied to
\[
 [w;C]x_a=\begin{pmatrix}wx_a\\\mu_ax_a\end{pmatrix}
\]
gives $wx_a\ge0$, since a negative first component would add a sign change to $x_a$. For a TP extension, $w-\varepsilon e_1^T$ is still a TP extension for sufficiently small $\varepsilon>0$, so $wx_a>0$. Apply this to $u,d$ to obtain
\[
 hx_a=dx_a+\mu_aux_a>0.
\]
The partial fractions in \eqref{ratio} are therefore
\begin{equation}\label{residues}
 \frac{R_{j+1}(z)}{q(z)}
 =u_j-\sum_{a=1}^m\frac{W_{aj}}{z+\mu_a},
 \quad W_{aj}=(hx_a)(y_a)_j,
 \quad q(z)=v_1\prod_{a=1}^m(z+\mu_a).
\end{equation}
The spectral assertion \eqref{prefix} of the induction hypothesis, applied to $C$, gives
\begin{equation}\label{Wflag}
 \det W[T,1:k]>0\qquad(|T|=k,\ 1\le k\le m).
\end{equation}
Indeed, each row is a positive multiple of $y_a/(y_a)_1$.

Fix $p$ positive increasing nodes $z_1<\cdots<z_p$. Divide the evaluation rows by $q(z_i)$ and subtract multiples of the first column. Cauchy--Binet gives
\begin{equation}\label{cv}
 \frac{\det[R_j(z_i)]_{i,j=1}^p}{\Delta(z)}
 =\sum_{|T|=p-1}a_T
   \prod_{i=1}^p\frac{q(z_i)}{\prod_{a\in T}(z_i+\mu_a)},
 \quad a_T=\Delta(\mu_T)\det W[T,1:p-1]>0.
\end{equation}
Here the Cauchy--Vandermonde determinant has the exact sign
\[
 \det\left[1,\frac{-1}{z_i+\mu_a}\ (a\in T)\right]_{i=1}^p
 =\frac{\Delta(z)\Delta(\mu_T)}
        {\prod_{i=1}^p\prod_{a\in T}(z_i+\mu_a)}.
\]
All terms in \eqref{cv} are positive. This includes $p=n$, when the only set $T$ is $1:m$.

Every TP matrix $B$ admits such a row $v$. Take $v_j=L^j$.
In every minor containing $v$, the coefficient of the highest power of $L$ is a positive
minor of $B$, so a sufficiently large $L$ makes $[B;v]$ TP.
At each eigenvalue $\lambda_i$ one has $R(\lambda_i)(\lambda_i I-B)=0$
and $q(\lambda_i)>0$, hence $R(\lambda_i)/q(\lambda_i)=\ell_i$.
Taking any increasing spectral subset as the nodes gives \eqref{prefix}.
\end{proof}

\section{Mixed evaluation determinants and scalar kernels}\label{cornersec}

We now use the adjugate rows from the preceding section to express
the corner minors of $f(A)$ as sums of scalar kernels with positive
coefficients.

Let $A\in\TP_n$, $n\ge3$.  For the moment put
\[
 B=A[1:n-1,1:n-1],\qquad v=A[n,1:n-1],\qquad
 R_A(z)=(e_n^T\adj(zI-A))[1:n-1].
\]
A block determinant identity gives
\begin{equation}\label{rowidentity}
 R_A(z)=v\adj(zI-B).
\end{equation}
Indeed, both sides are the first $n-1$ entries of the last adjugate row,
as follows by expanding the corresponding cofactors along the last row.
Since $[B;v]$ is TP, Lemma~\ref{spectral} gives the expansion
\eqref{residues} and the minor inequalities \eqref{Wflag} for this row,
with $m=n-2$.

Apply the same construction to $A^T$ and define
\begin{equation}\label{shared}
 E_A(z)=\adj(zI-A)_{1n}
       =b\prod_{a=1}^{m}(z+\mu_a),\qquad
 m=n-2,\quad b>0,\quad 0<\mu_1<\cdots<\mu_m.
\end{equation}
For $1\le j\le n-1$, put
$\phi_j(z)=\adj(zI-A)_{jn}/E_A(z)$.  Thus
\begin{equation}\label{actualflag}
 \phi_1=1,\qquad
 \phi_{j+1}(z)=u_j-H_j(z),\qquad
 H_j(z)=\sum_{a=1}^{m}\frac{W_{aj}}{z+\mu_a}
 \quad(1\le j\le m),
\end{equation}
where every initial-column flag minor of $W$ is positive:
$\det W[T,1:k]>0$ for $|T|=k$.
Here $W$ is the residue matrix obtained from \eqref{residues}
for $A^T$.  If $C$ denotes the TP matrix in that construction, then
$W=\diag(hx_a)Y_C$, where the rows of $Y_C$ are dual left
eigenvectors of $C$.  Consequently
\[
 W^{-1}\diag(\mu_1,\ldots,\mu_m)W=C\in\TP_m.
\]
The expansions below use only the positivity of the initial flag
minors.  The identity with $C$ is also needed for the general
preserver argument in Section~\ref{stsection}.

Order the eigenvalues of $A$ as
$0<\lambda_1<\cdots<\lambda_n$ and normalize its right eigenvector
matrix $X$ by $X_{1a}=1$.  At a spectral point,
\[
 \adj(\lambda_a I-A)
   =\chi_A'(\lambda_a)X[1:n,a](X^{-1})[a,1:n],
 \qquad \chi_A(z)=\det(zI-A),
\]
so $\phi_j(\lambda_a)=X_{ja}$ for $j<n$.
Lemma~\ref{spectral}, applied to $A^T$, shows that every initial-column
minor of $X^T$ is positive, including $\det X>0$.
This positivity statement concerns evaluations at the eigenvalues
of $A$.  The rational functions in \eqref{actualflag} represent
only the first $n-1$ rows of $X$ at these nodes.  The last row is
used only through its evaluation minors.

For functions $g_1,\ldots,g_k$ and nodes $z=(z_1,\ldots,z_k)$, write
\[
 \det[g_1,\ldots,g_k](z_i)
 =\det[g_1,\ldots,g_k]_z
 =\det[g_j(z_i)]_{i,j=1}^k.
\]
For functions with derivatives through order $k-1$, the Wronskian is
\[
 W(g_1,\ldots,g_k)(x)
 =\det[g_j^{(i-1)}(x)]_{i,j=1}^k.
\]
For an ordered system $g_1,\ldots,g_N$, the Wronskians
$W(g_1,\ldots,g_j)$, $1\le j\le N$, are called its initial
Wronskians.

Let $p+q=n$, with $1\le p,q<n$, and define, at increasing positive nodes,
\begin{equation}\label{sourcemixed}
 D_A^p[f](z_1,\ldots,z_n)
 =\det[\phi_1,\ldots,\phi_q\mid
                  f\phi_1,\ldots,f\phi_p](z_i).
\end{equation}
Subtracting $u_j$ times the first column within each block replaces
its nonconstant functions by $-H_j$ without changing the determinant.
At the eigenvalues of $A$, this determinant satisfies
\begin{equation}\label{sourcecorner}
 D_A^p[f](\lambda_1,\ldots,\lambda_n)
       =\det X\,G_p(f(A)).
\end{equation}
To verify it, the evaluation matrix in \eqref{sourcemixed} is
\[
 [X^T[1:n,1:q]\mid\diag(f(\lambda_a))X^T[1:n,1:p]].
\]
Multiply on the left by $X^{-T}$.  Its first $q$ columns become
$[e_1,\ldots,e_q]$, and the last $p$ columns become
$f(A)^T[1:n,1:p]$.  Expanding along the first $q$ columns leaves
$\det f(A)^T[q+1:n,1:p]=G_p(f(A))$, with positive sign.
The same identity applied to $A^T$ gives the lower-left corner $F_p$.
The formulas hold for arbitrary values of $f$ on the spectrum.

For positive increasing nodes $\lambda_1,\ldots,\lambda_n$, put
\[
 \chi_\lambda(z)=\prod_{a=1}^n(z-\lambda_a),\qquad
 \epsilon_p=(-1)^{p(p-1)/2},
\]
and denote the divided difference of order $n-1$ by
\[
 [g]_\lambda
 =\sum_{a=1}^n\frac{g(\lambda_a)}{\chi_\lambda'(\lambda_a)}.
\]
For $g\in C^{n-1}(I)$ and nodes in an interval $I$, repeated nodes
are understood by continuity. If $(T_1,\ldots,T_n)$ is uniform on the
simplex $T_i\ge0$, $\sum_i T_i=1$, the Hermite--Genocchi formula
\cite[formula~(52)]{deBoorDD} gives
\begin{equation}\label{hgdivided}
 [g]_\lambda=\frac{1}{(n-1)!}\,
  \mathbb E\,g^{(n-1)}\left(\sum_{i=1}^n T_i\lambda_i\right).
\end{equation}
For polynomials $E,D$ with $D(\lambda_a)\ne0$ for every $a$,
and $1\le p<n$, define
\begin{equation}\label{Hdef}
 H_{p,E,D}[f]
 =\epsilon_p\det\bigl([x^{i+j}E(x)f(x)/D(x)]_\lambda\bigr)_{i,j=0}^{p-1}.
\end{equation}
Thus the scalar kernel is defined at arbitrary distinct positive nodes; when
these are the spectrum of $A$, the denominator in the divided difference
is $\chi_A'(\lambda_a)$.

\begin{lemma}\label{kernel}
Suppose that $D$ and $E$ are finite products of factors $x+\mu$, with $\mu\ge0$, allowing positive constant factors, and put $r=\deg D$. For any function $f$ on the nodes $\lambda_1,\ldots,\lambda_n$,
\begin{equation}\label{quotient}
 H_{p,E,D}[f]=
 \frac{\det\bigl[\lambda_i^kD(\lambda_i)\ (0\le k<n-p)
 \mid\lambda_i^jE(\lambda_i)f(\lambda_i)\ (0\le j<p)\bigr]_{i=1}^n}
 {\Delta(\lambda)\prod_{i=1}^nD(\lambda_i)}.
\end{equation}
In particular, $H_{p,E,D}[x^t]>0$ when $t>n-p-1+r$.
\end{lemma}

\begin{proof}
Expanding the numerator of \eqref{quotient} along its last $p$ columns and dividing by the denominator yields
\[
 \epsilon_p\sum_{|S|=p}\Delta(\lambda_S)^2
 \prod_{s\in S}\frac{f(\lambda_s)E(\lambda_s)}{\chi_\lambda'(\lambda_s)D(\lambda_s)}.
\]
This is \eqref{Hdef}, by Cauchy--Binet.

Now take $f(x)=x^t$. Multiplication of a polynomial by $x+\mu$ acts on its ascending coefficients by a TN lower bidiagonal matrix. Hence the coefficient matrices of each of the two column families
\[
 D,xD,\ldots,x^{n-p-1}D,\qquad E,xE,\ldots,x^{p-1}E
\]
have nonnegative maximal minors. In the numerator of \eqref{quotient}, the first block has exponents at
most $n-p-1+r$, and the factor $x^t$ makes the exponents in the second
block at least $t$. Thus the two sets of exponents are strictly ordered
when $t>n-p-1+r$. Expanding by their coefficient matrices writes the numerator in \eqref{quotient} as a sum of nonnegative coefficients times generalized Vandermonde determinants
\[
 \det[\lambda_i^{\alpha_j}]_{i,j=1}^n>0,
 \qquad 0<\lambda_1<\cdots<\lambda_n,\quad \alpha_1<\cdots<\alpha_n.
\]
This last fact follows by setting $\lambda_i=e^{s_i}$: repeated Rolle's theorem gives nonvanishing for the exponential system, and coalescing nodes determines the positive sign from its Wronskian. The columns in each coefficient matrix are linearly independent, since multiplication by a nonzero polynomial is injective. Hence each coefficient matrix has at least one positive maximal minor, and at least one term of the expansion is positive. The denominator is positive as well.

\end{proof}

For a set $V\subseteq1:m$, put
\[
 D_V(x)=\prod_{a\in V}(x+\mu_a),\qquad
 c_V=\Delta(\mu_V)\det W[V,1:|V|]>0,
\]
with $D_\varnothing=c_\varnothing=1$.
For $|T|=q-1$ and $|U|=p-1$, define
\begin{equation}\label{ED}
 D_{TU}^A(x)=D_{U\setminus T}(x),\qquad
 E_{TU}^A(x)=D_{T\setminus U}(x).
\end{equation}
If $h=|T\cap U|$ and $r=\deg D_{TU}^A$, then
\begin{equation}\label{degrees}
 \begin{aligned}
 r&=p-1-h\le p-1,\qquad
 \deg E_{TU}^A=q-1-h=n-2p+r,\\
 q-1+r&=p-1+\deg E_{TU}^A=n-2-h\le n-2.
 \end{aligned}
\end{equation}
The superscript $A$ indicates that the poles and residues are obtained from $A$.

\begin{lemma}\label{cornerexp}
Let $A\in\TP_n$, $n\ge3$, and $p+q=n$ with $1\le p,q<n$.
With the data just defined, at any increasing positive nodes
$z_1<\cdots<z_n$ one has
\begin{equation}\label{paired}
 D_A^p[f](z)
 =\sum_{\substack{|T|=q-1\\|U|=p-1}}c_Tc_U
 \frac{\det\bigl[z_i^kD_{U\setminus T}(z_i)\ (0\le k<q)
       \mid f(z_i)z_i^jD_{T\setminus U}(z_i)\ (0\le j<p)\bigr]_{i=1}^n}
      {\prod_{i=1}^nD_{T\cup U}(z_i)}.
\end{equation}
At the eigenvalues of $A$, put
\[
 w_{TU}(A)=\frac{\Delta(\lambda)c_Tc_U}
                 {\det X\prod_{a=1}^nD_T(\lambda_a)}>0.
\]
Then
\begin{equation}\label{uppercornerH}
 G_p(f(A))
 =\sum_{\substack{|T|=q-1\\|U|=p-1}}
       w_{TU}(A)H_{p,E_{TU}^{A},D_{TU}^{A}}[f],
\end{equation}
and, using the same construction for $A^T$,
\begin{equation}\label{cornerH}
 F_p(f(A))
 =\sum_{\substack{|T|=q-1\\|U|=p-1}}
       w_{TU}(A^T)H_{p,E_{TU}^{A^T},D_{TU}^{A^T}}[f].
\end{equation}
The coefficients are independent of $f$.
\end{lemma}

\begin{proof}
Within the first block of \eqref{sourcemixed}, multilinearity and
Cauchy--Binet expand the nonconstant columns $-H_1,\ldots,-H_{q-1}$
by sets $T$ of size $q-1$, with coefficient $\det W[T,1:q-1]$.
The second block similarly uses $U$, with coefficient
$\det W[U,1:p-1]$.
For an ordered set $T$, the rational functions
\[
 1,\quad\frac{-1}{x+\mu_a}\quad(a\in T)
\]
are obtained from $1,x,\ldots,x^{|T|}$, divided by $D_T(x)$,
by a constant change of basis with determinant $+\Delta(\mu_T)$.
Indeed, evaluate both bases at $|T|+1$ increasing positive nodes;
the Cauchy--Vandermonde formula in the proof of
Lemma~\ref{spectral} gives precisely this determinant and sign.
The analogous change in the second block has determinant
$+\Delta(\mu_U)$.  The resulting term before clearing denominators is
therefore
\[
 c_Tc_U\det\bigl[z_i^k/D_T(z_i)\ (0\le k<q)
          \mid f(z_i)z_i^j/D_U(z_i)\ (0\le j<p)\bigr]_{i=1}^n.
\]
Multiplying row $i$ by $D_{T\cup U}(z_i)$ gives
$D_{U\setminus T}$ in the first block and
$D_{T\setminus U}$ in the second.  This proves \eqref{paired};
there is no interchange of the two blocks.

Now take $z=\lambda$.  By \eqref{quotient}, the numerator in each
term of \eqref{paired} equals
\[
 \Delta(\lambda)\prod_{a=1}^nD_{U\setminus T}(\lambda_a)
       H_{p,E_{TU}^A,D_{TU}^A}[f].
\]
Since $D_{T\cup U}=D_TD_{U\setminus T}$, its quotient by the
row denominator is
$\Delta(\lambda)H_{p,E_{TU}^A,D_{TU}^A}[f]/\prod_aD_T(\lambda_a)$.
Divide \eqref{paired} by $\det X$ and use \eqref{sourcecorner} to
obtain \eqref{uppercornerH}.  Applying that identity to $A^T$
gives \eqref{cornerH}, because $A$ and $A^T$ have the same spectrum.
\end{proof}

The derivation of \eqref{paired} applies to any two prefixes of lengths
$q,p$ of the same rational flag, evaluated at $p+q$ increasing positive
nodes. Only the passage to \eqref{uppercornerH} uses the full spectrum
of the matrix producing the flag.

\begin{lemma}\label{kernelcompare}
Let $n\ge3$, $1\le p<n$, and let $K$ be a compact subset of
$\{0<\lambda_1\le\cdots\le\lambda_n\}$.
Suppose that $D,E$ are products of factors $x+\mu$, $\mu\ge0$,
with positive constant factors allowed, and that
\[
 0\le r=\deg D\le p-1,\qquad
 \deg E=n-2p+r\ge0.
\]
For each fixed $t>n-2$, there are constants $0<c_t\le C_t<\infty$,
depending only on $n,p,K,t$, such that
\begin{equation}\label{purekernelcompare}
 c_tH_{p,E,D}[x^{n-1}]
 \le H_{p,E,D}[x^t]
 \le C_tH_{p,E,D}[x^{n-1}]
\end{equation}
uniformly over all such $D,E$ and all $\lambda\in K$.
There are also $c,\eta>0$, depending only on $n,p,K$, such that
\begin{equation}\label{kernelcompareend}
 H_{p,E,D}[x^{n-2+\beta}]
 \ge c\beta H_{p,E,D}[x^{n-1}]
 \qquad(0<\beta\le\eta).
\end{equation}
\end{lemma}

\begin{proof}
Both polynomial blocks in \eqref{quotient}, before multiplication
of the second block by $x^\tau$, have degree at most $n-2$.
By Lemma~\ref{kernel}, it suffices to compare
\[
 V_{I,J}(\lambda;\tau)
 =\det[\lambda_a^i\ (i\in I)\mid\lambda_a^{\tau+j}\ (j\in J)]_{a=1}^n,
 \qquad I,J\subseteq\{0,\ldots,n-2\},\quad |I|=n-p,\quad |J|=p.
\]
The coefficients and the denominator in the expansion do not depend
on $\tau$. Enlarge $K$ to the weakly ordered tuples in a compact
positive interval containing its nodes in the interior. We first take
distinct nodes. Choose $h\ge0$ so that $s_a=\log\lambda_a+h\ge0$
on this set. List the elements of $I$ and of $J$ in increasing
order, and let $b_1,\ldots,b_n$ be the resulting exponents. For $\tau>n-2$ they are strictly increasing
and nonnegative. Newton interpolation gives the identity
\[
 \frac{V_{I,J}(\lambda;\tau)}{\Delta(s)}
 =e^{-h\sum_j b_j}
   \det\bigl([e^{b_j x}]_{s_1,\ldots,s_a}\bigr)_{a,j=1}^n,
\]
so this quotient extends continuously to repeated nodes.
Expanding the exponential in finite Taylor sums and applying
Cauchy--Binet gives
\[
 \frac{\det[e^{s_a b_j}]_{a,j=1}^n}{\Delta(s)}
 \ge\frac{\Delta(b)}{\prod_{j=0}^{n-1}j!}.
\]
In fact every term is nonnegative by the generalized Vandermonde
inequality, and the displayed term uses exponents $0,\ldots,n-1$.
The finite Taylor matrices converge entrywise, and the bound
passes to repeated nodes by continuity.

For $\tau=t$ and $\tau=n-1$, the quotients by $\Delta(s)$ have
positive lower bounds and finite upper bounds there for each
index pair. There are only finitely many pairs, giving
\eqref{purekernelcompare} for distinct nodes. For $\tau=n-2+\beta$,
$0<\beta\le1$, the exponent list has at most one difference below
$1$, and that difference is $\beta$. Consequently $\Delta(b)\ge
\beta$, and $e^{-h\sum_j b_j}$ has a positive lower bound.
Comparing with the upper bound for the normalized determinant
at $n-1$ gives \eqref{kernelcompareend}. The divided-difference
formula \eqref{Hdef} extends both inequalities to repeated nodes.
\end{proof}

\section{Ordinary powers and the critical endpoint}\label{powersec}

\begin{proof}[Proof of Theorem~\ref{main}]
Let $A$ be oscillatory and fix $t>n-2$. For $n\ge3$, choose TP
matrices tending to $A$ with spectra in a fixed compact positive
interval. In the positive-weight expansion
\eqref{cornerH}, compare the two powers term by term using
\eqref{purekernelcompare}. Passage to the limit gives, for both
far corners of each order, a lower bound by a positive constant
times the corresponding positive minor of $A^{n-1}$. For $n=2$, the two
off-diagonal entries of $A^t$ are those of $A$ multiplied by the
positive divided difference of $x^t$ at its two eigenvalues.
Thus the far corners stay positive throughout $(n-2,\infty)$.
The set of parameters in this interval for which $A^t$ is TP is
open and contains $n-1$. At any limit point $A^t$ is TN, with
positive determinant and positive far corners, so it is TP by
Lemma~\ref{recognition}. The set is also relatively closed, proving
the assertion.

For sharpness when $n\ge3$, take
\[
 J_n(a)=\operatorname{tridiag}(1,a,1),\qquad a>2.
\]
It is positive definite and TN by its positive-pivot bidiagonal
$LDL^T$ factorization, and its positive adjacent entries make it
oscillatory.
If $\lambda_1,\ldots,\lambda_n$ are its eigenvalues, interpolation
and bandwidth give
\begin{equation}\label{sharp}
 (J_n(a)^t)_{1n}=[x^t]_\lambda.
\end{equation}
Indeed, that entry vanishes for $J_n(a)^k$, $0\le k\le n-2$, and
equals $1$ for $k=n-1$, so only the leading interpolation
coefficient contributes.
For $n-3<t<n-2$,
\[
 \frac{\mathrm d^{n-1}}{\mathrm dx^{n-1}}x^t
 =\prod_{j=0}^{n-2}(t-j)x^{t-n+1}<0.
\]
The divided-difference mean value theorem makes \eqref{sharp}
negative, whereas at $t=n-2$ the same entry is zero.
For $n=2$, the excluded endpoint is $A^0=I$.
\end{proof}

\begin{proof}[Proof of Corollary~\ref{sets}]
A nonsingular TN matrix is a limit of TP matrices. Its spectrum stays
away from zero, so Theorem~\ref{main} and continuity give TN powers
for $t>n-2$, and then for $t=n-2$. Nonnegative integer powers are
TN directly. TP matrices have TP positive integer powers; this
includes $n-2$ when $n\ge3$.

To exclude any noninteger $0<t<n-2$, put
$d=\lfloor t\rfloor+3\le n$ and use $J_d(a)\oplus I_{n-d}$.
Equation~\eqref{sharp} gives a negative entry because
$d-3<t<d-2$. A sufficiently small TP perturbation retains this
negative entry, so the obstruction also holds for TP matrices.

For $t<0$, take a symmetric positive definite $2\times2$ matrix
with positive off-diagonal entries and eigenvalues
$0<\lambda_1<\lambda_2$. Its off-diagonal entry after taking the
$t$th power is a positive multiple of
$(\lambda_2^t-\lambda_1^t)/(\lambda_2-\lambda_1)<0$.
Embed it as a leading block and perturb to TP.
The value $t=0$ gives $I$, which is not TP.
Finally, $J_n(a)^k$ has zero $(1,n)$ entry for every integer
$0\le k\le n-2$. These observations prove all the assertions.
\end{proof}

\subsection{Powers of singular matrices}

Let $\nu$ be the size of the largest Jordan block at zero. For
noninteger $t$, the ordinary power of a singular TN matrix is defined
when $x^t$ has right derivatives through order $\nu-1$ at zero.
This holds exactly for $t>\nu-1$, and the power then vanishes on
the generalized zero eigenspace. Nonnegative integer powers have
their usual meaning, including $A^0=I$.

\begin{lemma}[{\cite[Theorem~2.9]{GL}}]\label{tnrankfactorization}
Let $M$ be an $m\times m$ TN matrix of rank $r\ge1$. Then $M=BC$
for TN matrices $B\in\mathbb R^{m\times r}$ and
$C\in\mathbb R^{r\times m}$.
\end{lemma}

\begin{proposition}\label{singularpowers}
Fix integers $q\ge2$, $\nu\ge1$, and $n\ge q+\nu$. Let $A$ range
over the matrices in $\TN_n$ with $q$ positive eigenvalues counted
with algebraic multiplicity and largest Jordan block of size $\nu$
at zero. The real exponents for which $A^t$ is defined and TN for
every such $A$ form exactly
\[
 \mathbb Z_{\ge0}\cup[q+\nu-2,\infty).
\]
\end{proposition}

\begin{proof}
By Lemma~\ref{tnrankfactorization}, write $A=BC$. The matrix $CB$
is TN by Cauchy--Binet. Every column of $A$ is a linear combination
of the columns of $B$. Since $B$ has $\operatorname{rank}A$ columns,
they form a basis of the column space of $A$. The identity
$AB=B(CB)$ shows that $CB$ represents the restriction of $A$ to this
space in that basis. Each Jordan block at zero loses one dimension,
while the Jordan blocks at positive eigenvalues are preserved.
Iterating this factorization $\nu$ times leaves a nonsingular TN
matrix $D$ of order $q$.

For noninteger $t>\nu-1$, let $h$ be a Hermite interpolant of
$x^{t-\nu}$ at the positive eigenvalues, matching derivatives up
to one less than the largest block size at each eigenvalue.
Then $x^\nu h(x)$ has the same required derivatives as $x^t$ at
every eigenvalue of $A$, including zero. Hence $A^t=A^\nu h(A)$.
Applying $(BC)^m h(BC)=B(CB)^{m-1}h(CB)C$ at each step gives
\[
 A^t=P D^{\,t-\nu} Q,
\]
where $P$ and $Q$ are products of the rectangular TN factors.
For $t\ge q+\nu-2$, Corollary~\ref{sets} gives $D^{t-\nu}\in\TN_q$.
Cauchy--Binet gives $A^t\in\TN_n$, as it does for every nonnegative
integer power.

For sharpness, take the upper bidiagonal TN matrix $H$ whose diagonal
consists of $\nu$ zeros followed by $q$ ones, with all superdiagonal
entries one. Its diagonal blocks are a single nilpotent Jordan block
and $U=I_q+S_q$, where $S_q$ is the upper shift. Since
$e_1^T H^\nu=e_{\nu+1}^T$, the last $q$ entries of the first row of
$H^t=H^\nu h(H)$ form the first row of $h(U)=U^{t-\nu}$. Thus
\[
 (H^t)_{1,\nu+k}=\binom{t-\nu}{k-1},\qquad 1\le k\le q.
\]
For noninteger $\nu-1<t<q+\nu-2$, set $k=\lfloor t-\nu\rfloor+3$.
Then $2\le k\le q$ and the binomial coefficient has exactly one
negative factor in its numerator. Adjoining a zero block gives the
same entry in any larger order.
\end{proof}

\begin{corollary}\label{singularrank}
Fix $2\le r<n$. For every $A\in\TN_n$ of rank $r$ and every
$t\ge r-1$ for which $A^t$ is defined, $A^t$ is TN. The threshold
$r-1$ is sharp even for symmetric matrices. The sharp threshold
for simultaneous existence and total nonnegativity over the full
rank-$r$ class in $\TN_n$ is $r$.
\end{corollary}

\begin{proof}
If the Jordan block sizes at zero are $s_1,\ldots,s_\ell$, then
$r=q+\sum_{i=1}^\ell(s_i-1)$. For $q\ge2$ this gives
$q+\nu-2\le r-1$. For $q=1$ the same factorization gives
$A^t=\lambda^{t-\nu}A^\nu$ for noninteger $t>\nu-1$, where
$\lambda$ is the sole positive eigenvalue. For $q=0$ every defined
noninteger power is zero. In all cases $\nu\le r+1$, with equality
for the nilpotent upper shift of order $r+1$, padded by zeros if
necessary. Its noninteger powers require $t>r$.

For sharpness of $r-1$, let $B$ be the $r\times(r+1)$ matrix with
$B_{i,i}=B_{i,i+1}=1$ and all other entries zero. Then
$A=B^TB$ is symmetric TN of rank $r$, and $BB^T=J_r(2)$ is
positive definite. For $t>0$, the one-step factorization above gives
$A^t=B^T J_r(2)^{t-1}B$. Hence its $(1,r+1)$ entry equals
$(J_r(2)^{t-1})_{1r}$. By the bandwidth argument of \eqref{sharp},
this is the divided difference of $x^{t-1}$ at the $r$ positive
eigenvalues of $J_r(2)$. It is negative for $r-2<t<r-1$ by the
divided-difference mean value theorem. Adjoining zeros covers every
$n>r+1$.
\end{proof}

\subsection{Zeros at the critical endpoint}

The comparison below determines how a vanishing corner minor crosses zero at the endpoint of the universal range.

\begin{proposition}\label{linear}
Let $n\ge2$, $1\le p<n$, and let $A\in\TN_n$ be nonsingular.
There are $c,\eta>0$ such that
\begin{equation}\label{linearbound}
 F_p(A;n-2+\beta)\ge c\beta F_p(A;n-1)
 \qquad(0<\beta\le\eta).
\end{equation}
The same assertion holds for $G_p$. If either corner minor vanishes
at $t=n-2$ and is positive at $t=n-1$, then it has a simple zero
with positive derivative at $n-2$ and is negative immediately to
the left.
\end{proposition}

\begin{proof}
For $n\ge3$, take a TP approximation $A_\nu\to A$ whose spectra
lie in a fixed compact positive interval. The endpoint estimate
\eqref{kernelcompareend} applies to every kernel in \eqref{cornerH}
with the same constants, regardless of the distances between
the eigenvalues.
The positive coefficients in that expansion give
\[
 F_p(A_\nu;n-2+\beta)
 \ge c\beta F_p(A_\nu;n-1)
 \qquad(0<\beta\le\eta).
\]
Passing to the limit proves \eqref{linearbound}.

For $n=2$, $F_1(A;t)=A_{21}[x^t]_\lambda$,
including the confluent case. As $t\to0$, the quotient
$[x^t]_\lambda/t$ tends to $[\log x]_\lambda>0$.
Since $F_1(A;1)=A_{21}$, this gives \eqref{linearbound}.
Transposition proves the estimates for $G_p$.

The corner minors are real analytic in $t$.
If $F_p(A;n-2)=0<F_p(A;n-1)$, division of \eqref{linearbound} by $\beta$
and passage to the limit give
\[
 \left.\frac{\mathrm d}{\mathrm dt}F_p(A;t)\right|_{t=n-2}
 \ge cF_p(A;n-1)>0.
\]
It therefore has a simple zero and is negative immediately to the
left.  The argument for $G_p$ is identical.
\end{proof}

\begin{lemma}\label{integerprofilebound}
Let $A$ be a nonsingular TN matrix of order $n\ge2$ with
$a_{d,d+1}>0$ for $1\le d<n$. Then $G_p(A^k)>0$ for $1\le p<n$ and
integers $k\ge\max(p,n-p)$. If $A\in\Osc_n$ and $n\ge4$, every
minor of $A^{n-2}$ of order $2\le p\le n-2$ is positive.
\end{lemma}

\begin{proof}
In the elementary bidiagonal factorization of $A$
\cite[Theorem~2.2.2]{FJ}, retain one factor $I+uE_{d,d+1}$ with
$u>0$ at each cut $d$, and denote its position in the retained order
by $\pi(d)$. Here $E_{ij}$ is the matrix unit. Without such a factor, $A$ is block lower
triangular across $d$, whereas $A_{d,d+1}>0$. In Cauchy--Binet each
retained factor permits a move from a selected $d$ to a vacant $d+1$,
and all factors permit the index set to stay fixed.

Put $q=n-p$. To transform $1:p$ into $q+1:n$, number selected indices
from right to left and vacancies from left to right. Selected index
$a$ exchanges with vacancy $b$ at $d=p-a+b$ after exchanges
$(a-1,b)$ and $(a,b-1)$, whenever present. The preceding moves put
$a$ at $d$ and $b$ at $d+1$. Define $h(p)=0$ and
$h(d+1)-h(d)=\mathbf{1}_{\{\pi(d+1)<\pi(d)\}}$ for $1\le d<n-1$.
Use factor $\pi(d)$ in copy $a+h(d)$ of $A$ for exchange $(a,b)$.
From either predecessor the copy number increases by one or stays
fixed. In the latter case the factor position increases. Exchanges at the same cut
have distinct $a$ and hence distinct copies. The first and last
exchanges use copies $1$ and $p+h(q)\le\max(p,q)$, since the
increments of $h$ are $0$ or $1$. Cauchy--Binet gives $G_p(A^k)>0$
by keeping the indices fixed at unused factors.

If $A$ is oscillatory and $2\le p\le n-2$, then $\max(p,q)\le n-2$.
The preceding argument for $A$ and $A^T$ gives positivity of the two
far corners of order $p$ in $A^{n-2}$. By Pinkus's shadow theorem
\cite[Theorem~1.1]{Pinkus08}, any zero minor of that order would have
a principal subminor lying in a right or left shadow. In a right
shadow, moving each pair of indices $(i_s,j_s)$ to $(s,n-p+s)$ keeps the
corresponding subminor in that shadow, since $s\le i_s$ and $j_s\le n-p+s$.
Thus the upper-right corner would vanish. The left shadow gives the
same implication for the lower-left corner.
\end{proof}

\begin{proposition}\label{endpointfour}
Let $A\in\Osc_n$, $n\ge4$, and $m=n-2$. Every nonempty minor of $A^m$
is positive except possibly the four in \eqref{fourcorners} with $F=A^m$. In each
direction, at most one of the far corners of orders $1$ and $n-1$ can
vanish. All nine resulting zero patterns occur. Each vanishing minor has
a simple zero with positive derivative in $t$ at $m$. For $t<m$ sufficiently
close to $m$, the negative minors of $A^t$ are precisely the zero minors of $A^m$.
\end{proposition}

\begin{proof}
Lemma~\ref{integerprofilebound} gives positivity of all minors
of orders $2\le p\le n-2$ at $m$.
The positive adjacent and diagonal entries give a path of length $m$
from $i$ to $j$ whenever $|i-j|\le m$, so $(A^m)_{ij}>0$ there.
If $a_{ij}>0$ for some $j-i\ge2$, that step shortens the adjacent
path from $1$ to $n$ to at most $m$ steps, and diagonal steps pad
it to length $m$. If no such entry exists, each step increases the
index by at most one, so there is no path of length $m$ from $1$
to $n$. Hence $(A^m)_{1n}=0$ exactly when $A$ has upper bandwidth
one. Transposition gives the lower bandwidth criterion.

Put $E=\diag(1,-1,1,\ldots)$ and $B=EA^{-1}E$.
Jacobi's identity makes $B$ TN and $B^{n-1}$ TP, hence $B$ is oscillatory.
Every order-$(n-1)$ minor of $A^m$ is $\det(A^m)$ times its complementary
entry of $B^m$ with transposed indices. In particular,
$G_{n-1}(A^m)=\det(A^m)(B^m)_{1n}$. If $(A^m)_{1n}=0$, then
$A$ has upper bandwidth one and $A[1:n-1,2:n]$ is lower triangular
with positive diagonal. Jacobi's identity gives $B_{1n}>0$, hence
$(B^m)_{1n}\ge B_{11}^{m-1}B_{1n}>0$. The two upper corners
therefore cannot both vanish, and transposition handles the lower corners.

To realize all patterns, let $N$ be the upper shift and take $A=LU$,
choosing $U$ and $L^T$ independently from
\[
 I+N,\qquad (I-N)^{-1},\qquad (I+N)^2 .
\]
These factors are TN by bidiagonality, checkerboard inversion and
multiplication, respectively. Their positive adjacent entries make
$A$ oscillatory. Write $\widehat T=ET^{-1}E$, so $B=\widehat U\widehat L$.
By triangular support and positive diagonals, the upper bandwidths of $A$
and $B$ equal those of $U$ and $\widehat U$, respectively. The first two
choices are exchanged by checkerboard inversion and only $I+N$ has
bandwidth one. The third choice and its checkerboard inverse
$(I-N)^{-2}$ both have bandwidth at least two. The same argument for
$A^T$ makes the lower choice independent. This realizes all nine patterns.
Proposition~\ref{linear} gives the simple zeros, and continuity keeps
all other minors positive in a neighborhood of $m$.
\end{proof}

For an oscillatory matrix $A$ of order $n\ge2$, define its individual ordinary critical exponent and its integer oscillation exponent by
\[
 g(A)=\inf\{\gamma:A^t\in\TP_n\text{ for every real }t>\gamma\},
 \qquad r(A)=\min\{k\in\mathbb Z_{>0}:A^k\in\TP_n\}.
\]
Since $A^0=I$ is not TP, Theorem~\ref{main} gives
$0\le g(A)\le n-2$. The integer bound gives $r(A)\le n-1$.

\begin{corollary}\label{attainment}
For $A\in\Osc_n$ with $n\ge2$, the following are equivalent:
$g(A)=n-2$, $r(A)=n-1$, and $A^t\notin\TN_n$ for every $t$ in some
left neighborhood of $n-2$.
\end{corollary}

\begin{proof}
For $n=2$, $A^0=I$ and Proposition~\ref{linear} give the assertion.
Suppose $n\ge3$. Once an integer power is TP, it stays TP at every
later integer by Cauchy--Binet, since each selected column set of the
nonsingular TN right factor has full rank. If $r(A)\le n-2$, then $A^{n-2}$
is TP and remains so in a neighborhood of $n-2$ by openness. Together
with Theorem~\ref{main}, this gives $g(A)<n-2$.
If $r(A)=n-1$, then $A^{n-2}$ is TN and nonsingular but has a
vanishing proper far corner by Lemma~\ref{recognition}. That corner
is negative immediately to the left by Proposition~\ref{linear}.
Such an interval forces $g(A)\ge n-2$, with equality by Theorem~\ref{main}.
\end{proof}

\section{Operator-monotone multipliers}\label{coeffsec}\label{multsec}

To determine the signs of the scalar kernels for $x^\alpha e^{\psi(x)}$,
we express the relevant Wronskians in terms of determinants of Taylor coefficients.
For real power sums $p_1,p_2,\ldots$, define $h_j$ by
\[
 H(t)=\sum_{j\ge0}h_jt^j=\exp\left(\sum_{r\ge1}\frac{p_r}{r}t^r\right),
 \qquad h_j=0\quad(j<0).
\]

\begin{lemma}\label{translation}
Let $q\ge1$, $U\ge0$, and define
\[
 p_r=r\int_{[0,U]}v^{r-1}\,\mathrm d\nu(v),\qquad
 \nu=(q-1)\,\mathrm dv+\eta_0+\sum_{j=1}^{m} z_j\eta_j
 \quad\text{on }[0,U],
\]
where the $\eta_j$ are finite nonnegative measures and the $z_j$
are indeterminates. For increasing nonnegative integers $K=(k_0,\ldots,k_{q-1})$,
put $d=\sum_i(k_i-i)$ and $L=p_1-(q-1)U$. Let $b_K$ be the number
of ways to lower $K$ to $(0,\ldots,q-1)$, decreasing one entry by one
at each step and retaining strictly increasing nonnegative entries.
Then
\begin{equation}\label{coefficientbound}
 \det[h_{k_i-j}]_{i,j=0}^{q-1}-b_K L^d/d!
 \quad\text{has nonnegative coefficients in }z_1,\ldots,z_m .
\end{equation}
The empty lowering path gives $b_K=1$ at $d=0$. In particular
$b_K\ge1$. For scalar power sums $p_r=\sum_j a_j u_j^r$ with
$U=u_0>u_1>\cdots>0$, the scalar bound applies when every partial
sum $a_0+\cdots+a_j$ is at least $q-1$: take the density
$\sum_j a_j\mathbf1_{[0,u_j]}$.
\end{lemma}

\begin{proof}
Suppose a series $H_c(t)$ satisfies
$\partial_c H_c=\delta(c)tH_c+t^2\partial_t H_c$. Its coefficients
satisfy $h_j'=(\delta(c)+j-1)h_{j-1}$. Writing
$T_k=\det[h_{k_i-j}]_{i,j=0}^{q-1}$ gives
\begin{equation}\label{lowering}
 T_k'=\sum_i(\delta(c)+k_i-q)T_{k-e_i}.
\end{equation}
In the row derivative, split $\delta(c)+k_i-j-1$ as
$(\delta(c)+k_i-q)+(q-j-1)$. The last terms collect by columns
and replace column $j$ by a multiple of column $j+1$, so vanish.
A lowered tuple with a negative or repeated index contributes zero.
For constant $\delta$, the same differential equation describes
the translation $H_c(t)=(1-ct)^{-\delta}H_0(t/(1-ct))$.

First suppose that all $\eta_j$ have continuous densities, and let
$\varphi$ be the resulting polynomial density of $\nu$. Put
$\delta(c)=\varphi(U-c)$ and use
\[
 H_c(t)=\exp\left(\int_0^c
             \frac{\delta(s)t}{1-(c-s)t}\,\mathrm ds\right),
 \qquad
 L(c)=\int_0^c(\delta(s)-q+1)\,\mathrm ds .
\]
Differentiation verifies the differential equation, and $H_U=H$.
Compare polynomials coefficient by coefficient. The degree-zero
determinant is $1$. Every admissible coefficient in
\eqref{lowering} is $L'(c)+(k_i-1)$ and therefore dominates
$L'(c)\ge0$. When $d>0$ there is at least one legal step, and splitting all paths
by their first step gives $b_K=\sum_{i\ {\rm legal}}b_{K-e_i}$.
Induction on $d$ gives
$T_k'\ge b_K L'(c)L(c)^{d-1}/(d-1)!$ coefficient by coefficient.
Since $T_k=0$ at $c=0$ for $d>0$, integration proves the bound.

For $U>0$, approximate each $\eta_j$ weakly by continuous
nonnegative densities of the same mass. Each coefficient of the
determinant depends polynomially on finitely many moments, so
the inequalities pass to the limit. For $U=0$, add the background
measure $(q-1)\,\mathrm dv$ on $[0,\varepsilon]$ and let
$\varepsilon\downarrow0$. The factor $b_K$ is sharp uniformly in $U$.
At $U=0$, $H(t)=e^{Lt}$ and differentiation with respect to $L$ gives
$\partial_L T_k=\sum_{i\ {\rm legal}}T_{k-e_i}$, hence equality.
At $L=1$, multiplying row $i$ by $k_i!$ gives the falling factorials
of $k_i$ in degrees $0,\ldots,q-1$. The Vandermonde formula yields
$b_K/d!=\prod_{i<j}(k_j-k_i)/\prod_i k_i!$.
\end{proof}

\begin{lemma}[{\cite[Theorems 5.2 and 5.4]{Hansen}}]\label{Pickrep}
A real function $\psi$ is operator monotone on $(0,\infty)$ if and
only if it has the representation
\begin{equation}\label{Pick}
 \psi(z)=\psi(1)+\int_{[0,\infty]}
          \frac{(1+s)(z-1)}{z+s}\,\mathrm d\mu(s),\qquad
 \mu\ge0,\quad \mu([0,\infty])<\infty.
\end{equation}
At $s=\infty$ the kernel is $z-1$. Equivalently, $\psi$ has a
holomorphic continuation to the upper half-plane with nonnegative
imaginary part, extending across the positive axis. Its derivative
at $1$ equals $\mu([0,\infty])$, and if $\psi$ is nonconstant
then $\psi'(x)>0$ for all $x>0$.
\end{lemma}

We also write this representation as
\[
 \psi(z)=\psi(1)+b(z-1)
  +\int_{[0,\infty)}
      \left(\frac1{1+s}-\frac1{z+s}\right)\mathrm d\rho(s),
 \qquad b=\mu(\{\infty\}),\quad \mathrm d\rho=(1+s)^2\,\mathrm d\mu
\]
for finite $s$. In particular,
$\psi'(x)=b+\int(x+s)^{-2}\,\mathrm d\rho(s)$.

\begin{lemma}\label{ECT}
Let $\phi_1,\ldots,\phi_N$ be real analytic on an interval, with all
initial Wronskians positive. For every initial prefix of length $j$,
the normalized evaluation determinant
\[
 \frac{\det[\phi_\ell(x_i)]_{i,\ell=1}^{j}}{\Delta(x)}
\]
extends to a strictly positive continuous function of weakly
increasing nodes in the interval.
\end{lemma}

\begin{proof}
Newton interpolation expresses this quotient as the determinant
with $(i,\ell)$ entry $[x_1,\ldots,x_i]\phi_\ell$. This formula
extends continuously to repeated nodes. To see that it never
vanishes, the identity
\begin{equation}\label{csECT}
 W\big((\phi_2/\phi_1)',\ldots,(\phi_{k+1}/\phi_1)'\big)
       =\frac{W(\phi_1,\ldots,\phi_{k+1})}{\phi_1^{k+1}}>0
\end{equation}
and induction with Rolle's theorem bound the zeros of a nonzero
linear combination of the first $j$ functions by $j-1$,
counting multiplicity. Indeed, divide by $\phi_1>0$ and
differentiate unless only its coefficient is nonzero.
A singular divided-difference matrix would give a nonzero
combination with all Newton coefficients zero, and thus a
zero at each node with its assigned multiplicity. Finally,
the node domain is connected, and at a single repeated node
the determinant is the positive Wronskian divided by
$\prod_{i=0}^{j-1}i!$.
\end{proof}

\subsection{Positive scalar kernels}

\begin{lemma}\label{toeplitz-wronskian}
Let $F$ be positive and smooth on an interval, and define $h_j(x)$
as the Taylor coefficients of $F(x)/F(x-t)$ at $t=0$.
For $q,v\ge1$, one has
\begin{equation}\label{Wronskian}
 W(1,x,\ldots,x^{q-1},F,xF,\ldots,x^{v-1}F)
 =\left(\prod_{j=0}^{q+v-1}j!\right)F(x)^v
       \det[h_{v+i-j}(x)]_{i,j=0}^{q-1}.
\end{equation}
\end{lemma}

\begin{proof}
Divide every function by $F$, which divides the Wronskian by
$F^{q+v}$. Move the $v$ monomials to the front, with sign
$(-1)^{qv}$. Center the monomials in each block at $x$ and
divide derivative row $j$ by $j!$, counting rows from zero.
The monomials now give a unit block. The remaining entries
are the Taylor coefficients of $(z-x)^j/F(z)$ in degrees
$v+i$. The Taylor coefficient of degree $j$ in $F(x)/F(x+t)$
is $(-1)^jh_j(x)$. Factoring $F(x)^{-1}$ from the $q$ columns
and removing the alternating signs gives $F(x)^{-q}(-1)^{qv}$ times the indicated
Toeplitz determinant. The two signs cancel.
\end{proof}

\begin{lemma}\label{Pickkernel}
Let $n\ge2$, $1\le p<n$, and $q=n-p$. Let $D$ be a product of at
most $p-1$ factors $x+d$ with $d\ge0$, and $E$ any finite product
of factors $x+e$ with $e\ge0$. Positive constant factors are
allowed in both. If $\Phi(x)=x^{n-2}e^{\psi(x)}$ for a nonconstant
operator-monotone $\psi$, then
\begin{equation}\label{Pickdet}
 H_{p,E,D}[\Phi]
 =\frac{\det[x_a^i\ (0\le i<q)
             \mid x_a^j\Phi(x_a)E(x_a)/D(x_a)\ (0\le j<p)]_{a=1}^n}
        {\Delta(x)}
\end{equation}
has a strictly positive continuous extension to weakly increasing
positive nodes.
\end{lemma}

\begin{proof}
The identity follows from \eqref{quotient} by dividing each row by
$D(x_a)$. Put $F=\Phi E/D$. At a fixed $x>0$, set $U=1/x$,
$v_s=(x+s)^{-1}$, and
$p_k(x)=(-1)^{k-1}(\log F)^{(k)}(x)/(k-1)!$.
If $\rho_x$ is the image of $(x+s)^{-2}\,\mathrm d\rho(s)$
under $s\mapsto v_s$, then $p_k(x)=k\int v^{k-1}\,\mathrm d\nu(v)$
on $[0,U]$, with
\[
 \mathrm d\nu(v)=
 \left(n-2+\sum_e\mathbf 1_{[0,v_e]}(v)
             -\sum_d\mathbf 1_{[0,v_d]}(v)\right)\mathrm dv
       +b\,\mathrm d\delta_0(v)+\mathrm d\rho_x(v).
\]
The sums are over the factors of $E,D$ with multiplicity. Each
$\log(x+s)$ contributes $v_s^k=k\int_0^{v_s}v^{k-1}\,\mathrm dv$.
The remaining terms follow by differentiating $\psi$ in
Lemma~\ref{Pickrep}. The density displayed is at least $q-1$,
and the excess mass is at least $\psi'(x)>0$. Thus for every increasing
nonnegative $K$, with $d=\sum_i(k_i-i)$,
\begin{equation}\label{localcoefficients}
 \det[h_{k_i-j}(x)]_{i,j=0}^{q-1}\ge b_K\,\psi'(x)^d/d!.
\end{equation}

The series for these power sums is $F(x)/F(x-t)$. Taking $k_i=v+i$
in \eqref{localcoefficients} makes the determinant in \eqref{Wronskian}
positive for every $1\le v\le p$. The prefixes consisting only
of monomials have positive Wronskians as well. Lemmas
\ref{toeplitz-wronskian} and \ref{ECT} now give the conclusion.
\end{proof}

\subsection{Leading and trailing principal blocks}

Let $A=X\Lambda Y\in\TP_n$, with $Y=X^{-1}$ and with its eigenvalues
$\lambda_a$ in increasing order. Normalize $X_{1a}=1$ and put
\[
 \omega_a=Y_{a1},\qquad
 u_i(\lambda_a)=X_{i+1,a},\qquad
 v_j(\lambda_a)=Y_{a,j+1}/\omega_a,\qquad 0\le i,j<n.
\]
Equation~\eqref{projection} gives $\omega_a>0$, and $u_0=v_0=1$.
By Lemma~\ref{spectral} for $A^T$ and $A$, the evaluation determinants
of any initial prefix of $u$ and $v$ at increasing spectral subsets
are positive. For a scalar function $f$ on the spectrum,
\begin{equation}\label{properblockbiorthogonal}
 f(A)_{i+1,j+1}=
 \sum_{a=1}^n\omega_a f(\lambda_a)u_i(\lambda_a)v_j(\lambda_a).
\end{equation}
For $p,r\ge1$ and $d=p+r\le n$, Cauchy--Binet gives
\begin{align}\label{properblockGram}
 \det f(A)[1:p,r+1:d]
 =\sum_{|L|=d}\left(\prod_{a\in L}\omega_a\right)
 &\det[u_0,\ldots,u_{r-1},fu_0,\ldots,fu_{p-1}]_{\lambda_L}
 \notag\\[-1mm]
 &\times\det[v_0,\ldots,v_{d-1}]_{\lambda_L}.
\end{align}
Indeed, use \eqref{properblockbiorthogonal} with the constant function
$1$ in the first $r$ rows and with $f$ in the last $p$ rows. The first
rows are $[I_r\ 0]$, and expansion along them leaves the target block.
In the mixed determinant, the $u_i$ come from \eqref{actualflag}.
Subtracting their constant parts within the two blocks gives the
prefixes of $(1,-H_1,\ldots)$ and $f(1,-H_1,\ldots)$.
The largest needed index is $\max(p,r)-1\le n-2$.
When $n=2$ both prefixes consist of their first function alone.

Inserting \eqref{paired} into \eqref{properblockGram} and using
\eqref{quotient} expresses each upper boundary minor as a positive
linear combination of $H_{p,E,D}[f]$ at $d=p+r$ nodes, with
$\deg D\le p-1$, $\deg E\le r-1$. The coefficients are independent
of $f$. For $n=2$ this is the single kernel with $D=E=1$.

\begin{lemma}\label{mothercornercompare}
Let $p,r\ge1$, $d=p+r$, $n\ge d$, and $0<a<b$. Suppose that
$\Phi(x)=x^\alpha e^{\psi(x)}$, where $\psi$ is real operator monotone,
$\alpha\ge d-2$, and $\alpha>d-2$ or $\psi$ is nonconstant.
For nonsingular $A\in\TN_n$ with $\sigma(A)\subseteq[a,b]$, put
$M_f=\det f(A)[1:p,r+1:d]$. There are $0<c\le C<\infty$,
depending only on $d,a,b,\Phi$, such that
\[
 cM_{x^{d-1}}\le M_\Phi\le CM_{x^{d-1}}.
\]
The same holds after transposition or simultaneous reversal of the
indices. These boundary minors of $\Phi(A)$ are positive for $A\in\TP_n$.
In particular the comparison covers both proper far-corner families
when $d=n$.
\end{lemma}

\begin{proof}
Absorb $(\alpha-d+2)\log x$ into $\psi$. Both $\Phi$ and the reference
$x^{d-1}$ then satisfy Lemma~\ref{Pickkernel} in order $d$.
Normalize every factor $x+s$ in $E,D$ to $(x+s)/(1+s)$, which
multiplies both kernels by the same positive number. The ensuing
products lie in compact sets of at most $r-1$ and $p-1$ factors
of the form $tx+1-t$, $0\le t\le1$. The ratio of the two kernels
is continuous and strictly positive there and on the compact set
of weakly increasing nodes in $[a/2,2b]$, by Lemma~\ref{Pickkernel}
and its divided-difference formula. Thus
\begin{equation}\label{motherkernelbounds}
 c H_{p,E,D}[x^{d-1}]\le H_{p,E,D}[\Phi]\le C H_{p,E,D}[x^{d-1}].
\end{equation}
Summing with the positive weights above gives the comparison and
strict positivity for TP inputs. Approximation by TP matrices
and continuity of the analytic functional calculus give the
comparison on TN inputs, including at Jordan blocks. There are
only finitely many pairs $(p,r)$ for the fixed sum $d$.
\end{proof}

\begin{lemma}\label{principalbounds}
Let $A$ be a nonsingular TN matrix with ordered eigenvalues
$\lambda_1\le\cdots\le\lambda_n$, and let $B_f$ denote either its
leading or trailing principal block of order $k$ after applying a real
analytic function $f$ on $(0,\infty)$. There are $w_L\ge0$ with
$\sum_{|L|=k}w_L=1$, independent of $f$, such that
\begin{equation}\label{blockcharacteristic}
 \det(zI_k-B_f)=\sum_{|L|=k}w_L\prod_{a\in L}(z-f(\lambda_a)).
\end{equation}
All $w_L$ can be taken positive for $A\in\TP_n$.
If $f>0$ and $b_1\le\cdots\le b_n$ are the values $f(\lambda_a)$
in increasing order, then
\[
 b_1\cdots b_k\le\det B_f\le b_{n-k+1}\cdots b_n.
\]
For $A\in\TP_n$ and $k<n$, equality at either end requires all
the $f(\lambda_a)$ to coincide.
\end{lemma}

\begin{proof}
For the leading block of a TP matrix, apply Cauchy--Binet to
\eqref{properblockbiorthogonal} with $z-f$. The weights are the
products of the $\omega_a$ and the two prefix evaluation determinants
on $\lambda_L$. They sum to $1$ by comparing leading coefficients.
For general $A$, take a TP approximation and pass to a subsequence
on which all weights converge in the finite simplex. This choice
depends only on the approximating matrices. Continuity of analytic
functional calculus gives \eqref{blockcharacteristic}, including at
coincident eigenvalues. Reversing the indices gives the trailing case.
At $z=0$ the determinant is a weighted average of the $k$-fold
products. If all weights are positive, equality forces these products
to coincide. For $k<n$, compare subsets differing in one index.
\end{proof}

Every real eigenvalue of $B_f$ lies between the extreme values of
$f(\lambda_a)$, by \eqref{blockcharacteristic}. For $k=n-1$, all
eigenvalues of $B_f$ are real and interlace the ordered values
$f(\lambda_a)$, denoted by $b_1\le\cdots\le b_n$. Dividing
\eqref{blockcharacteristic} by $\prod_i(z-b_i)$ gives
$\sum_i w_{i^c}/(z-b_i)$ with the weights relabelled accordingly.
For positive weights and distinct $b_i$, this sum decreases from
$+\infty$ to $-\infty$ in each gap, giving strict interlacing.
A limit in the nodes and weights gives weak interlacing in general.

\begin{theorem}\label{properblock}
Let $n\ge N\ge2$, $A\in\TP_n$, and let $\Phi(x)=x^\alpha e^{\psi(x)}$,
where $\psi$ is real and operator monotone on $(0,\infty)$ and
$\alpha\ge N-2$. The leading and trailing principal blocks of order
$N$ of $\Phi(A)$ are TP, except when $N=2$, $\alpha=0$ and $\psi$
is constant, in which case $\Phi(A)$ is a positive multiple of the identity.
\end{theorem}

The blocks are taken after applying $\Phi$ to $A$.
We use the initial-minor criterion for total positivity. It tests all
minors on consecutive rows and columns with row set or column set
starting at the first index \cite[Theorem~3.1.4]{FJ}.

\begin{proof}
At $\alpha=N-2$ with constant $\psi$, use the integer power
$A^{N-2}$ for $N\ge3$ and the identity for $N=2$. Otherwise,
Lemma~\ref{principalbounds} makes the leading principal minors
positive, and Lemma~\ref{mothercornercompare} gives every other
initial minor of the leading block, using $p+r\le N$ and
transposition. The initial-minor criterion gives total positivity.
Reversal treats the trailing block.
\end{proof}

\begin{proof}[Total nonnegativity in Theorem~\ref{multipliers}]
Theorem~\ref{properblock} with $N=n$ gives the TN assertion on
TP inputs, including the constant-function case in order two.
Approximate by TP matrices and use continuity of analytic functional
calculus. Analyticity across zero allows the same argument for
singular TN inputs.
\end{proof}

The block size, rather than the ambient order, determines the
corresponding pure-power exponent set.

\begin{corollary}\label{blockset}
Let $n\ge N\ge1$.  The real numbers $t$ for which the leading
and trailing $N\times N$ blocks of $A^t$ are TP for every
$A\in\TP_n$ form exactly
\[
 \begin{cases}
 \mathbb R,&N=1,\\
 (0,\infty),&N=2,\\
 \mathbb Z_{>0}\cup[N-2,\infty),&N\ge3.
 \end{cases}
\]
\end{corollary}
\begin{proof}
For $N=1$, the spectral Gram formula makes the two endpoint
diagonal entries positive for every real $t$.  For $N\ge2$,
Theorem~\ref{properblock} and positive integer powers give sufficiency.
If $0<t<N-2$ is not an integer, put $d=\lfloor t\rfloor+3\le N$
and take $A_0=(I_d+S_d)\oplus I_{n-d}$, with $S_d$ the upper shift.
The entry $(1,d)$ of $A_0^t$ is $\binom t{d-1}<0$.
For $t<0$ use $d=2$, whose $(1,2)$ entry is $t<0$.
These are invertible TN matrices with all principal minors equal to
$1$.  If $P_\varepsilon\in\TP_n$ tends to $I$, then
$P_\varepsilon A_0$ is TP by Cauchy--Binet, and continuity of the
principal power near the spectrum $\{1\}$ preserves the negative
entry.  At $t=0$, $A^0=I$ is not TP when $N\ge2$.
Reversal gives the same obstruction for the trailing block.
\end{proof}

\subsection{Zero blocks and local reduction}

For a nonsingular TN matrix, Pinkus's shadow criterion describes zero
minors in terms of zero minors with consecutive rows and columns
\cite[Theorem~1.1]{Pinkus08}. Cuts at zeros immediately above or below
the diagonal give obstructions to nonzero minors for every analytic
function of a fixed matrix.

\begin{lemma}\label{zerocuts}
Let $H=(h_{ij})$ be a nonsingular TN matrix of order $n$. Then
\[
 h_{k,k+1}=0\ \Longrightarrow\ H[1:k,k+1:n]=0,
 \qquad
 h_{k+1,k}=0\ \Longrightarrow\ H[k+1:n,1:k]=0.
\]
Every scalar function analytic near $\sigma(H)$ preserves these zero
blocks. If increasing equal-size index sets $I,J$ fail
\eqref{zeropath}, with $H$ in place of $A$, then
$\det g(H)[I,J]=0$ for every such function $g$.
\end{lemma}

\begin{proof}
The diagonal entries of $H$ are positive. If $h_{k,k+1}=0$, then
the minor on rows $\{i,k\}$ and columns $\{k,k+1\}$ forces
$h_{i,k+1}=0$ for every $i<k$. The minor on rows $\{i,k+1\}$
and columns $\{k+1,j\}$ then forces $h_{ij}=0$ for $i\le k$ and
$j>k+1$. Transposition gives the lower-left zero block. Analytic
functional calculus preserves each block triangular form and
applies the function separately to its diagonal blocks.

A failure $i_r\le k<j_r$ across an upper cut leaves the first
$r$ rows of $g(H)[I,J]$ supported on at most $r-1$ columns.
A failure across a lower cut is the transposed case. Hence the
minor vanishes.
\end{proof}

\begin{lemma}[{\cite[Theorem~1]{SStriangular}}]
\label{triangularcornercriterion}
Let $U$ be a nonsingular upper triangular TN matrix of order $d$.
If $G_p(U)>0$ for $1\le p<d$, then, for increasing index sets
$R=\{r_1<\cdots<r_q\}$ and $S=\{s_1<\cdots<s_q\}$,
\[
 \det U[R,S]>0
 \quad\Longleftrightarrow\quad
 r_j\le s_j\quad(1\le j\le q).
\]
The corresponding assertion for a lower triangular matrix is
obtained by transposition.
\end{lemma}

\begin{proposition}\label{localreduction}
Let $A$ be nonsingular and TN, let $f$ be positive and real analytic
on $(0,\infty)$, and put $F_C=f(A[C,C])$. Then $f(A)$ is TN
if and only if $F_C$ is TN for each
$C\in\mathcal U(A)\cup\mathcal L(A)$.
In that case the positive minors of $f(A)$ are exactly those
specified by \eqref{zeropath} if and only if every proper upper-right
corner of $F_C$ is positive for $C\in\mathcal U(A)$ and every
proper lower-left corner of $F_C$ is positive for $C\in\mathcal L(A)$.

For $n\ge d\ge2$, the function $f$ preserves every nonsingular
$A\in\TN_n$ with $d(A)\le d$ if and only if it preserves every
nonsingular matrix in $\TN_d$.
\end{proposition}

\begin{proof}
Put $F=f(A)$. By Lemma~\ref{zerocuts}, $F$ is block lower triangular
for $\mathcal U(A)$, block upper triangular for $\mathcal L(A)$,
and $F[C,C]=F_C$ for blocks in either partition. Necessity of the
local TN conditions follows by taking these principal blocks.

Suppose that all $F_C$ are TN. They are nonsingular because $f$ is
positive on their input spectra. For each $C\in\mathcal U(A)$,
write the TN factorization $F_C=L_CD_CU_C$ with unit triangular
factors and $D_C$ positive diagonal. Set
$U=\bigoplus_{C\in\mathcal U(A)}U_C$ and
$D=\bigoplus_{C\in\mathcal U(A)}D_C$. Then $L=FU^{-1}D^{-1}$
is block lower triangular with diagonal blocks $L_C$, so it is unit
lower triangular. This gives the unique factorization $F=LDU$.
Transposing the argument for $\mathcal L(A)$ shows that $L$ is
block diagonal with TN blocks, just as $U$ is for $\mathcal U(A)$.
Cauchy--Binet gives $F\in\TN_n$.

For $C\in\mathcal U(A)$, the upper-right corner of $F_C=L_CD_CU_C$
of order $p$ is the corresponding corner of $U_C$ multiplied by
the product of the first $p$ entries of $D_C$. Only the first
$p$ intermediate indices survive in Cauchy--Binet. Thus the local
corner hypothesis makes each permitted triangular minor of $U_C$
positive by Lemma~\ref{triangularcornercriterion}. The lower
blocks of $L$ have the transposed property.

For $I=\{i_1<\cdots<i_p\}$ and $J=\{j_1<\cdots<j_p\}$ satisfying
\eqref{zeropath}, put $K=\{\min(i_r,j_r):1\le r\le p\}$. Its
entries are strictly increasing. Each $i_r$ and $\min(i_r,j_r)$
belong to the same lower block, and $\min(i_r,j_r)$ and $j_r$
belong to the same upper block. Hence $\det L[I,K]$ and
$\det U[K,J]$ are products of positive permitted block minors.
The term indexed by $K$ in Cauchy--Binet for $LDU$ is positive,
so $\det F[I,J]>0$. Lemma~\ref{zerocuts} gives all the remaining
zeros. Conversely, each local corner in the hypothesis is
itself a minor allowed by \eqref{zeropath}.

For the final equivalence, embed smaller matrices in order $d$
by taking direct sums with identities, then use the local
criterion. Conversely, embed any nonsingular $B\in\TN_d$ as
$A=B\oplus I_{n-d}$. Then $d(A)\le d$ and $f(B)$ is a
principal block of $f(A)$.
\end{proof}

\subsection{Strict positivity for spectral multipliers}

\begin{lemma}\label{directioncorners}
Let $B=(b_{ij})$ be a nonsingular TN matrix of order $d\ge2$, and
let $h(x)=x^\beta e^{\chi(x)}$, where $\chi$ is real and operator
monotone, $\beta\ge d-2$, and $\beta>d-2$ or $\chi$ is
nonconstant. If $b_{k,k+1}>0$ for every $1\le k<d$, then
$G_p(h(B))>0$ for $1\le p<d$. If instead all $b_{k+1,k}$ are
positive, then all $F_p(h(B))$ are positive.
\end{lemma}

\begin{proof}
If $b_{k,k+1}>0$ for every $k$, then $G_p(B^{d-1})>0$ by
Lemma~\ref{integerprofilebound}, since $d-1\ge\max(p,d-p)$.
Lemma~\ref{mothercornercompare} gives $G_p(h(B))>0$.
Apply the same argument to $B^T$ for the lower-left corners.
\end{proof}

\begin{proof}[Proof of Theorem~\ref{zeropattern}]
If $d=d(A)=1$, Lemma~\ref{zerocuts} makes $A$ positive diagonal.
Every positive analytic function of $A$ is then positive diagonal,
and the criterion reduces to $I=J$.
Suppose $d\ge2$. For each $C\in\mathcal U(A)\cup\mathcal L(A)$
of order $q\ge2$, the matrix $B=A[C,C]$ is nonsingular TN
and $\alpha\ge d-2\ge q-2$. Thus $\Phi(B)$ is TN by
Theorem~\ref{multipliers}. Scalar blocks give positive values.
Proposition~\ref{localreduction} now gives $\Phi(A)\in\TN_n$.

If $\alpha>d-2$ or $\psi$ is nonconstant, then the strict
hypotheses of Lemma~\ref{directioncorners} hold in each nonscalar
block. In each upper block of $A$ the upper adjacent entries are
positive, and in each lower block the lower adjacent entries
are positive. The lemma gives the corresponding local corner
inequalities, so Proposition~\ref{localreduction} gives the
claimed positive minors of $\Phi(A)$. The same argument applies
to $x^{d-1}$ and $x^{n-1}$.
\end{proof}

\begin{proof}[Strictness in Theorem~\ref{multipliers}]
All adjacent entries of an oscillatory matrix are positive, so
Theorem~\ref{zeropattern} gives every minor positive under the
stated strictness condition.
\end{proof}

\begin{proof}[Proof of Corollary~\ref{localpowers}]
By Proposition~\ref{localreduction}, the TN-preserving exponents
for the stated source class are exactly those for nonsingular
$\TN_d$. Corollary~\ref{sets}, in order $d$, gives
$\mathbb Z_{\ge0}\cup[d-2,\infty)$.
Theorem~\ref{zeropattern} gives the full positivity criterion
for every $t>d-2$.

Among the remaining TN-preserving exponents, only integers
$0\le k\le d-2$ need to be excluded from the stable zero-pattern
assertion. Let $S_d$ be the upper shift and take
\[
 A=\diag(I_d+S_d,I_{n-d}),
\]
omitting the second block if $d=n$. Then $d(A)=d$, the $(1,d)$
entry of $A^k$ is zero, and that of $A^{n-1}$ is
$\binom{n-1}{d-1}>0$. This excludes all these integers.
The case $d=1$ follows from the positive diagonal assertion of
Theorem~\ref{zeropattern}.
\end{proof}

\begin{corollary}\label{analyticmaxrank}
Under the strictness hypotheses of Theorem~\ref{zeropattern}, a
minor of $\Phi(A)$ vanishes exactly when the corresponding minor
of $g(A)$ vanishes for every scalar function $g$ analytic near
$\sigma(A)$. For any row and column sets $R,C$, join $i\in R$
to $j\in C$ when $\ell_i\le j\le u_i$, where $\ell_i$ is the
first index of the block of $\mathcal L(A)$ containing $i$ and
$u_i$ is the last index of its block in $\mathcal U(A)$. If
$\nu_A(R,C)$ denotes the maximum size of a matching in this
bipartite graph, then
\[
 \max_g\operatorname{rank}g(A)[R,C]
 =\operatorname{rank}\Phi(A)[R,C]
 =\operatorname{rank}A^{n-1}[R,C]
 =\nu_A(R,C).
\]
For positive diagonal $A$, the same assertions hold for every
positive analytic $\Phi$.
\end{corollary}

\begin{proof}
Theorem~\ref{zeropattern} and Lemma~\ref{zerocuts} give the
assertion about vanishing minors. Every minor of $g(A)$ that
is nonzero satisfies \eqref{zeropath} and hence gives a matching
of the selected rows and columns. Conversely, both endpoints
$\ell_i,u_i$ are nondecreasing in $i$. If a matching pairs
$i<i'$ with $j>j'$, then
$\ell_i\le\ell_{i'}\le j'<j\le u_i\le u_{i'}$,
so the columns can be exchanged. Successive exchanges of
adjacent inverted pairs give an increasing matching of the same
size. Its row and column sets satisfy \eqref{zeropath}, giving a
positive minor of both $\Phi(A)$ and $A^{n-1}$. This also applies
when $A$ is positive diagonal, for which the graph has only
edges between equal indices.
\end{proof}

\begin{remark}
For fixed $A,I,J$ with $d(A)\ge2$, the power minor
$\det A^t[I,J]$ is either positive for every $t>d(A)-2$,
or zero for every real $t$. The alternatives are distinguished
by \eqref{zeropath} and follow from Theorem~\ref{zeropattern}
and Lemma~\ref{zerocuts}, respectively. The first allows zeros
or negative values below the stated threshold.
\end{remark}

\begin{corollary}\label{localattainment}
Let $A$ be nonsingular and TN, with $d=d(A)\ge2$, and put $m=d-2$.
Its individual TN threshold $\inf\{\gamma:A^t\in\TN_n\text{ for all }t>\gamma\}$
equals $m$ if and only if a minor permitted by \eqref{zeropath}
vanishes in $A^m$. Equivalently, $d-1$ is the first nonnegative
integer $k$ at which every permitted minor of $A^k$ is positive.
In this case some permitted minor has a simple zero at $m$ with positive
derivative, and $A^t$ fails to be TN throughout a left neighborhood of $m$.
\end{corollary}

\begin{proof}
At integer times a positive minor persists: the Cauchy--Binet term
indexed by $J$ for $\det A^{k+1}[I,J]$ contains the positive factor
$\det A[J,J]$. At $d-1$ every permitted minor is positive by
Theorem~\ref{zeropattern}. If all are positive at $m$, they stay
positive in a neighborhood of $m$, and all forbidden minors vanish
identically by Lemma~\ref{zerocuts}. The threshold is then less than $m$.

Otherwise Proposition~\ref{localreduction} applied to $x^m$ gives
a vanishing far corner of $B^m$ in the direction of some block
$B=A[C,C]$. Its order must be $d$, since $m>q-2$ on every smaller
nonscalar block of order $q$ and Lemma~\ref{directioncorners} gives
strict positivity there. The same lemma makes this corner positive
in $B^{d-1}$. Proposition~\ref{linear} gives a simple zero with
positive derivative at $m$, also for $d=2$. Since $A^t[C,C]=B^t$,
the corner is negative in $A^t$ immediately to the left of $m$.
Theorem~\ref{zeropattern} gives $A^t\in\TN_n$ for $t>m$.
\end{proof}

Other permitted minors can have multiple zeros. If $A=(I_d+S_d)\oplus
(I_d+S_d)$, with $S_d$ the upper shift, the minor of $A^t$ on rows
$\{1,d+1\}$ and columns $\{d,2d\}$ equals $\binom{t}{d-1}^2$,
which has a double zero at $d-2$.

\section{Principal minors and real powers}\label{principalSection}

For $\Phi(x)=x^\alpha e^{\psi(x)}$, let $E=\diag(1,-1,1,\ldots)$
and put
\[
 \Phi^\#(x)=1/\Phi(1/x)=x^\alpha e^{-\psi(1/x)}.
\]
Here $x\mapsto-\psi(1/x)$ is operator monotone, since inversion and
negation each reverse the positive definite order. If $A$ is nonsingular
TN, so is $B=EA^{-1}E$. Functional calculus and Jacobi's identity give
\begin{equation}\label{imchecker}
 \Phi^\#(B)=E\Phi(A)^{-1}E,\qquad
 \det\Phi^\#(B)[I,J]
   =\frac{\det\Phi(A)[J^c,I^c]}{\det\Phi(A)}.
\end{equation}
The factors from $E$ cancel the signs in Jacobi's identity.

\begin{proof}[Proof of the $P$-matrix assertion in Theorem~\ref{iminternal}]
The argument works for $n\ge3$ and $\alpha\ge n-3$.
If $\Phi$ is constant, the result is immediate. Otherwise apply
Theorem~\ref{properblock} with $N=n-1$ to TP approximations of $A$
and $B$. The leading order-$(n-1)$ blocks of $F=\Phi(A)$ and
$G=\Phi^\#(B)$ are TN by continuity and nonsingular by
Lemma~\ref{principalbounds}, so they are $P$-matrices. A proper
principal minor of $F$ either lies in this block or is complementary
to a principal minor in the leading block of $G$. In the latter case
\eqref{imchecker} gives positivity because $\det F=\prod_i\Phi(\lambda_i)>0$.
\end{proof}

\begin{lemma}\label{imentries}\label{imouterprincipal}
For every nonsingular $A\in\TN_n$, $n\ge2$, and $t\in\mathbb R$,
all principal minors of $A^t$ of orders $1,n-1,n$ are positive.
If $A$ is TP and $t>0$, its adjacent entries $(A^t)_{ij}$,
$|i-j|=1$, are also positive.
\end{lemma}

\begin{proof}
For $t>0$, take an integer $k>t$, put $B=A^k$ and $\beta=t/k$,
and use the integral representation
\begin{equation}\label{imstieltjes}
 A^t=\frac{\sin\pi\beta}{\pi}\int_0^\infty s^{\beta-1}K_s\,\dd s,
 \qquad K_s=B(B+sI)^{-1}=I-s(B+sI)^{-1}.
\end{equation}
Expanding the determinant and the diagonal cofactor gives
\begin{equation}\label{imdiagonalkernel}
 (K_s)_{ii}=\frac{\displaystyle\sum_{S\ni i}s^{n-|S|}\det B[S,S]}
                  {\det(B+sI)}>0,
\end{equation}
since all coefficients are nonnegative and $\det B>0$.
For adjacent $i,j$, the surviving terms from $sI$ in
$(B+sI)[j^c,i^c]$ lie on its diagonal. Thus its determinant has
nonnegative coefficients and, for TP $B$, a positive constant term.
The cofactor sign is $-1$, so $(K_s)_{ij}>0$. Integration gives both
entry assertions. For $t<0$, the diagonal positivity follows by
applying the positive-time result to $EA^{-1}E$, and $t=0$ is immediate.
Now $\det A^t=\prod_i\lambda_i^t>0$, and Jacobi gives
$\det A^t[i^c,i^c]=\det(A^t)(A^{-t})_{ii}>0$.
\end{proof}

\subsection{The excluded exponents}

The excluded positive exponents can be detected at the same two
indices in every order. Write $H_j=\{1,j\}$ for $j\ge2$.

\begin{lemma}\label{im:P-necessary}
Let $n\ge4$, let $k$ be an integer with $0\le k\le n-4$, and let
$K\subset(k,k+1)$ be nonempty and compact. There exists
$A\in\TP_n$ such that $\det A^t[H_n,H_n]<0$ for every $t\in K$.
\end{lemma}

\begin{proof}
Put $d=k+4$, $m=d-3$, and $t_0=\min K$. We first construct an
oscillatory matrix $B_r$ of order $d$ for which the principal
minor of $B_r^{-t}$ on $H_d$ is negative throughout $K$. Its
smallest eigenvalue contributes a rank-one term. We estimate its
interaction with the remaining spectral terms. Choose $a>0$
so that $(a+1)t_0>1$. Let $S$ be the lower shift of order $d$,
and set
\[
 L=(I-S)^{-1},\qquad
 D_r=\diag(1,\ldots,1,r^{-a},1),\qquad
 B_r=L D_r U_r,\quad r>0,
\]
where $U_r$ is unit upper bidiagonal, with superdiagonal entries
\[
 u_i=1\quad(1\le i\le d-3),\qquad
 u_{d-2}=r,\qquad u_{d-1}=r^a.
\]
These factors are TN, $\det B_r=r^{-a}$, and all entries
immediately above and below the main diagonal of $B_r$ are
strictly positive.  Thus $B_r$ is oscillatory
and has distinct positive eigenvalues.

Define the monic polynomials
\[
 P_0=1,\qquad P_1=z-1,\qquad
 P_j=(z-1)P_{j-1}-zP_{j-2}.
\]
They are the characteristic polynomials of
$(I-S_j)^{-1}(I+S_j^T)$, where $S_j$ is the lower shift of order
$j$. These matrices are nonsingular TN with
positive adjacent entries and are therefore oscillatory; consequently
all their roots are positive and simple.  Also $P_j(0)=(-1)^j$.  The tridiagonal pencil
\[
 zI-B_r=L\bigl[z(I-S)-D_rU_r\bigr]
\]
gives, on writing $q_r(z)=\det(zI-B_r)$,
\begin{align*}
 \theta_r(z)&=(z-r^{-a})P_{d-2}(z)-rzP_m(z),\\
 q_r(z)&=(z-1)\theta_r(z)-zP_{d-2}(z).
\end{align*}
The resolvent $R_r(z)=(zI-B_r)^{-1}$ satisfies
\begin{equation}\label{im:P-resolvent}
 (R_r)_{1d}=\frac r{q_r(z)},\qquad
 (R_r)_{d1}=\frac{z^{d-2}}{q_r(z)},\qquad
 (R_r)_{dd}=\frac{\theta_r(z)}{q_r(z)}.
\end{equation}
For the middle identity one must include the right factor $I-S$:
the numerator is $z^{d-1}-z^{d-2}(z-1)=z^{d-2}$.

Put $H(z)=(z-1)P_m(z)$, and denote its $m+1$ positive roots,
counting multiplicity, by $\rho_1,\ldots,\rho_{m+1}$. As
$r\to\infty$,
\[
 \frac{q_r(z)}r\longrightarrow-zH(z)
\]
uniformly on fixed compact sets, together with derivatives.
Thus one eigenvalue $\lambda_s$ tends to zero, $d-2$ eigenvalues
converge to the roots of $H$, counting multiplicity, and the
remaining eigenvalue $\lambda_\ell$ tends to infinity.  More precisely,
\begin{equation}\label{im:P-scales}
 \lambda_s\sim r^{-(a+1)},\qquad
 q'_r(\lambda_s)\sim(-1)^m r,\qquad
 \lambda_\ell\sim r,\qquad
 q'_r(\lambda_\ell)\sim r^{d-1}.
\end{equation}
The first two assertions follow from $q_r(0)=(-1)^d r^{-a}$;
the large eigenvalue follows from
$\operatorname{tr}B_r=r+2d-3+r^{-a}$.  Root counts can be taken
on fixed disjoint contours, so repeated roots of $H$ cause no
difficulty.

Let $Q_r$ be the spectral projection of $\lambda_s$.  Taking
residues in \eqref{im:P-resolvent}, and using
$\theta_r(\lambda)=\lambda P_{d-2}(\lambda)/(\lambda-1)$ at
every root $\lambda\ne1$ of $q_r$, gives
\begin{align*}
 (Q_r)_{1d}&\longrightarrow(-1)^m,&
 (Q_r)_{d1}&=O\bigl(r^{-(a+1)(d-2)-1}\bigr),\\
 (Q_r)_{dd}&\sim r^{-a-2},&
 (Q_r)_{11}&=O\bigl(r^{-(a+1)m}\bigr).
\end{align*}
For the last estimate, use rank one and the exact identity
\[
 (Q_r)_{11}
 =\frac{r\lambda_s^{d-3}(\lambda_s-1)}
 {P_{d-2}(\lambda_s)q'_r(\lambda_s)}.
\]

Write $B_r^{-t}=\lambda_s^{-t}Q_r+G_r$. We estimate $G_r$
by separating the bounded eigenvalues from $\lambda_\ell$.
The resulting bounds are
\begin{equation}\label{im:P-remainder}
 (G_r)_{11},(G_r)_{dd},(G_r)_{1d}=O(1),\qquad
 r(G_r)_{d1}\longrightarrow
 -[x^{m-t}]_{\rho_1,\ldots,\rho_{m+1}}.
\end{equation}
For the bounded eigenvalues, integrate $z^{-t}R_r(z)$ on a fixed
contour enclosing the roots of $H$ in the open right half-plane.
There $r(R_r)_{d1}\to-z^m/H(z)$, giving the confluent divided
difference in \eqref{im:P-remainder}. The remaining entries of $R_r$
are bounded because all cofactors of the tridiagonal pencil are
$O(r)$ there: only its $(d-2,d-1)$ entry is of order $r$ and all
others are bounded. Also $q_r/r$ stays away from zero on the contour.
At the large eigenvalue, the projection entries at $(1,d),(d,1),(d,d)$
and $(1,1)$ are respectively $O(r^{2-d}),O(r^{-1}),O(r^{-1})$ and
$O(r^{2-d})$. This follows from \eqref{im:P-resolvent} and
\eqref{im:P-scales}, with the rank-one formula for $(1,1)$ used above.
Since $\lambda_\ell^{-t}\le\lambda_\ell^{-t_0}$ for large $r$, this
projection contributes $o(r^{-1})$ to $(G_r)_{d1}$ and bounded terms
to the other entries. All bounds and limits are uniform on $K$.

Since $Q_r[H_d,H_d]$ has rank at most one,
\begin{align*}
 \det B_r^{-t}[H_d,H_d]
 ={}&\lambda_s^{-t}\bigl(
 (Q_r)_{11}(G_r)_{dd}+(Q_r)_{dd}(G_r)_{11}\\
 &\hspace{23mm}-(Q_r)_{1d}(G_r)_{d1}
                  -(Q_r)_{d1}(G_r)_{1d}\bigr)
       +\det G_r[H_d,H_d].
\end{align*}
After multiplication by $r$, all terms in parentheses except the
third tend uniformly to zero. Also $\det G_r[H_d,H_d]=O(1)$,
$r^{1-(a+1)t_0}\to0$, and $(r^{a+1}\lambda_s)^{-t}\to1$
uniformly on $K$. Thus
\[
 r^{1-(a+1)t}\det B_r^{-t}[H_d,H_d]
 \longrightarrow (-1)^m[x^{m-t}]_{\rho_1,\ldots,\rho_{m+1}}
\]
uniformly on $K$. Since $0<m-t<1$, the $m$th derivative of
$x^{m-t}$ has sign $(-1)^{m-1}$, as does its divided difference,
including confluent nodes. The limit is continuous and strictly
negative on $K$. Jacobi's identity and $\det B_r=r^{-a}$ also give
\begin{equation}\label{im:P-negative-limit}
 r^{1-t}\det B_r^t[2:d-1,2:d-1]
 \longrightarrow (-1)^m[x^{m-t}]_{\rho_1,\ldots,\rho_{m+1}}<0.
\end{equation}

Choose one $r$ for which $\det B_r^{-t}[H_d,H_d]<0$ throughout
$K$. With $E$ the checker diagonal of order $d$, the nonsingular
TN matrix $EB_r^{-1}E$ has the same negative minor in its
positive $t$th power. The principal power is jointly continuous
in the matrix and $t$ near this positive spectrum. Compactness
therefore gives a TP approximation $C$ with
$\det C^t[H_d,H_d]<0$ for all $t\in K$.

To increase the order by one, insert a zero row and a zero column
at the same interior position of $C$. The resulting matrix $C_0$
is TN, and its eigenvalues are those of $C$ together with one
simple zero. Take TP matrices tending to $C_0$ by density of TP matrices in the TN class. Their
spectral
projections converge since the limiting spectrum is simple. For
the eigenvalue $\mu$ tending to zero, $\mu^t\le\mu^{t_0}$ once
$\mu<1$. The powers thus converge uniformly on $K$ to $C^t$
with the same zero row and column inserted. A sufficiently close
TP approximation retains the negative minor on the two endpoints
throughout $K$. Iteration gives order $n$.
\end{proof}

\begin{proof}[Proof of Theorem~\ref{imPexact}]
For $n\le3$, every nonempty principal order is covered by
Lemma~\ref{imouterprincipal}; for $n=1$ the assertion is scalar.
Let $n\ge4$.  The $P$-matrix part of Theorem~\ref{iminternal},
already proved above using only leading blocks, gives sufficiency
for $t\ge n-3$. Nonnegative integer powers remain nonsingular
TN and hence are $P$-matrices.  Jacobi's identity shows that a nonsingular matrix is
a $P$-matrix if and only if its inverse is.  This gives all stated
negative exponents as well. Lemma~\ref{im:P-necessary} gives the
negative minor on $H_n$ throughout each compact subset of an
excluded positive interval. For a compact subset of an excluded
negative interval, apply the lemma to its reflection and replace its
witness $A$ by $EA^{-1}E\in\TP_n$. Indeed,
$(EA^{-1}E)^{-t}=EA^tE$ has the same principal minors as $A^t$.
This proves the formula for $T_n$ and its characterization by the
single principal minor. The TP witnesses also give necessity on
$\Osc_n$ and on all nonsingular TN matrices.
\end{proof}

The bound in Theorem~\ref{iminternal} is sharp uniformly over the
operator-monotone multiplier.  If $\alpha<n-3$, choose a nonintegral
$t$ with $\max(\alpha,0)<t<n-3$ and set
$\psi(x)=(t-\alpha)\log x$.  This function is operator monotone,
while $x^\alpha e^{\psi(x)}=x^t$.  By Lemma~\ref{im:P-necessary},
there is a TP matrix $A$ for which $\det A^t[H_n,H_n]<0$.

\begin{proposition}\label{uniformP}
For $n\ge1$, let $\mathcal U_n$ consist of the real numbers $\alpha$
such that $A^\alpha e^{\psi(A)}$ is a $P$-matrix for every
nonsingular $A\in\TN_n$ and every real operator-monotone function
$\psi$ on $(0,\infty)$. Then
\[
 \mathcal U_1=\mathcal U_2=\mathbb R,\qquad
 \mathcal U_n=[n-3,\infty)\quad(n\ge3).
\]
\end{proposition}

\begin{proof}
The scalar case is immediate. For a nonsingular $A\in\TN_2$
and any positive analytic $f$, Lemma~\ref{principalbounds} gives
positive diagonal entries of $f(A)$, and
$\det f(A)=\prod_i f(\lambda_i)>0$.
Thus every real $\alpha$ belongs to $\mathcal U_2$.

For $n\ge3$, the leading-block argument above gives sufficiency
when $\alpha\ge n-3$. The preceding sharpness argument supplies
necessity for $n\ge4$.

It remains to exclude $\alpha<0$ in order three. Put
\[
 K=\begin{pmatrix}8&12&8\\12&26&20\\4&14&13\end{pmatrix}.
\]
Its entries are positive and its determinant is $64$. Its minors
of order two, with row and column pairs ordered as $12,13,23$, are
\[
 \begin{pmatrix}64&64&32\\64&72&44\\64&76&58\end{pmatrix}.
\]
Hence $K\in\TP_3$. Its characteristic polynomial is
$q(z)=z^3-47z^2+194z-64$. The signs
$q(0)=-64$, $q(1)=84$, $q(4)=24$ and
$q(17/4)=-747/64$, together with $q(z)\to\infty$, locate its middle
root $\lambda$ in $(4,17/4)$.

\begin{samepage}
The corresponding spectral projection $P_\lambda$ satisfies
\[
 (P_\lambda)_{22}
 =\frac{\lambda^2-21\lambda+72}{q'(\lambda)}<0,
\]
because the numerator is at least $13/16$ on this interval, while
$q'(\lambda)<0$ at the middle simple root.
\end{samepage}

Fix $\alpha<0$ and use the operator-monotone function
$\psi(x)=-1/x$. Every $A_k=(K/\lambda)^k$, $k\ge1$, is TP.
The function $f(x)=x^\alpha e^{-1/x}$ tends to zero both at zero
and at infinity, and $f(1)=e^{-1}$. The three spectral projections
of $A_k$ are independent of $k$, so
$f(A_k)\to e^{-1}P_\lambda$. In particular, $f(A_k)_{22}<0$
for all sufficiently large $k$. Thus $\alpha\notin\mathcal U_3$.
\end{proof}

The quantifier over $\psi$ is essential here: although every real
power of a nonsingular TN matrix of order three is a $P$-matrix, this property
does not persist for every operator-monotone exponential multiplier.

For $n\ge4$ and nonsingular $A\in\TN_n$, consider the special case
\[
 h(x)=a+\frac b x+\int_{(0,\infty)}\frac{\mathrm d\mu(s)}{x+s},
 \quad a,b\ge0,\quad \mu\ge0,\quad
 \int_{(0,\infty)}\frac{\mathrm d\mu(s)}{1+s}<\infty,
\]
with $h\not\equiv0$.  The displayed representation makes $xh(x)$
positive and operator monotone, so $\log(xh(x))$ is operator
monotone by composition with the operator logarithm.
Theorem~\ref{iminternal}, applied to
$x^{n-2}h(x)=x^{n-3}e^{\log(xh(x))}$, shows that $A^{n-2}h(A)$
is a $P$-matrix and all its minors of orders $2,\ldots,n-2$
are nonnegative. If $A$ is TP, these minors are strictly positive.

\begin{corollary}\label{localP}
Let $A$ be a nonsingular TN matrix. Let $\mathcal B(A)$ be the
common refinement of $\mathcal U(A)$ and $\mathcal L(A)$, and put
$b(A)=\max_{C\in\mathcal B(A)}|C|$. For every function $f$ analytic
near $\sigma(A)$ and every principal index set $I$,
\begin{equation}\label{localPfactor}
 \det f(A)[I,I]
 =\prod_{C\in\mathcal B(A)}
   \det f(A[C,C])[I\cap C,I\cap C],
\end{equation}
where each factor uses the indices inherited from $C$, and an empty
factor is $1$. In particular, for positive real-analytic $f$,
$f(A)$ is a $P$-matrix if and only if every $f(A[C,C])$ is.

Over all nonsingular $A\in\TN_n$ with $b(A)\le s$, $1\le s\le n$,
the exact exponent set for $A^t$ to be a $P$-matrix is $T_s$ of
Theorem~\ref{imPexact}: it is $\mathbb R$ for $s\le3$ and
$\mathbb Z\cup\{t:|t|\ge s-3\}$ for $s\ge4$. For
$A^\alpha e^{\psi(A)}$ with arbitrary real operator-monotone $\psi$ on $(0,\infty)$,
the uniform range is $\mathbb R$ for $s\le2$ and $[s-3,\infty)$ for $s\ge3$.
\end{corollary}

\begin{proof}
Each cut of $\mathcal B(A)$ occurs where at least one of the two
adjacent entries is zero. By Lemma~\ref{zerocuts}, the corresponding
two-block triangular form is preserved by $f$, which acts on its
diagonal blocks separately. Restricting to the same rows and columns
in $I$ preserves the triangular form. Successive cuts give
\eqref{localPfactor}, and taking $I\subseteq C$ gives necessity of
the blockwise $P$ condition.

Theorem~\ref{imPexact} and Proposition~\ref{uniformP} give the stated
ranges on blocks of size at most $s$. These ranges shrink with the
order, so they suffice on the whole class by \eqref{localPfactor}.
Conversely, every nonsingular $B\in\TN_s$ occurs as the leading
principal block of $B\oplus I_{n-s}$, whose common-refinement blocks have size at most $s$ (omit the identity when $s=n$).
Thus the order-$s$ ranges are necessary as well.
\end{proof}

For example, let $A=I_n+N$ with $N$ the upper shift and $n\ge3$.
Here $b(A)=1$ and $d(A)=n$. The upper triangular matrix $A^t$ has
unit diagonal for every real $t$, so all its principal minors equal $1$.
Yet $(A^{1/2})_{13}=\binom{1/2}{2}=-1/8$.

\section{Common flags and boundary blocks}\label{opsection}\label{cssection}

Theorem~\ref{properblock} gives the bound $p+r-2$ for the minor
on rows $1:p$ and columns $r+1:r+p$. The following estimate
replaces this bound by $\max(p,r)-1$.

\begin{theorem}\label{edge}\label{farcornerPhi}
Let $p,r\ge1$, $n\ge p+r$, $A\in\TP_n$, and
$\Phi(x)=x^\alpha e^{\psi(x)}$ with $\psi$ real operator monotone
on $(0,\infty)$ and $\alpha\ge\max(p,r)-1$. Then
\begin{equation}\label{csedgepositive}
 \det\Phi(A)[1:p,r+1:r+p]>0
\end{equation}
unless $p=r=1$, $\alpha=0$, and $\psi$ is constant, when the
entry vanishes. The result also holds after transposition or
simultaneous reversal. In particular, for $n\ge4$ and $\alpha\ge n-3$,
both far corners of every order $2\le p\le n-2$ are positive.
\end{theorem}

In \eqref{properblockGram} the two blocks of the mixed determinant
use prefixes of lengths $r$ and $p$ of the same rational flag.
Expanding products of flag minors lets us group the signed terms in
the Wronskian into nonnegative sums.

\subsection{A mixed determinant for two prefixes}

For a residue matrix $W$ write $P_h(W;T)=\det W[T,1:h]$,
with $P_0(W;\varnothing)=1$.

\begin{theorem}\label{csmother}
Let $q,s\ge1$, $m\ge\max(q,s)-1$, and
$0<\mu_1<\cdots<\mu_m$.  Suppose that the first
$\max(q,s)-1$ columns of $W$ have positive initial flag minors, and put
\[
 H_j(x)=\sum_{a=1}^m\frac{W_{aj}}{x+\mu_a},\qquad
 \mathcal F_d=(1,-H_1,\ldots,-H_{d-1}).
\]
Let $\Phi(x)=x^\alpha e^{\psi(x)}$, where $\psi$ is real and
operator monotone on $(0,\infty)$.  If
$\alpha\ge\max(q,s)-1$, then, for all
$0<x_1<\cdots<x_{q+s}$,
\begin{equation}\label{csevaluation}
       \det[\mathcal F_q\mid\Phi\mathcal F_s](x_i)\ge0.
\end{equation}
All initial Wronskians of $(\mathcal F_q,\Phi\mathcal F_s)$
and the determinant in \eqref{csevaluation} are strictly positive,
with precisely the following exception. If $q=s$, $\alpha=q-1$,
$\psi$ is constant, and either $q=1$ or $m=q-1$, then the proper
initial Wronskians are positive and the full evaluation matrix has
rank $2q-1$ at every set of $2q$ distinct positive nodes. Its full
Wronskian vanishes identically.
\end{theorem}

\subsection{A positive decomposition of flag products}

Draw a partial matching on an ordered set by arcs on one side of
the line. The matching is exposed if the arcs do not cross and
no unmatched vertex lies between the two ends of an arc. For a
partition $Y=I\sqcup J$, color $I$ negative and $J$ positive.
In the matching expansion for flag minors \cite[Section~3]{DKK},
a compatible matching pairs vertices of opposite colors and leaves
only positive vertices unmatched.

\begin{lemma}\label{csnetwork}
Let $m\ge l\ge k\ge0$ and let $W$ be an $m\times l$ matrix with
positive initial flag minors. Fix $S\subseteq Z\subseteq1:m$ with
$|S|\le k$ and $|Z|=k+l-|S|$, and put $Y=Z\setminus S$.
There are nonnegative coefficients $d_M(W)$, independent of the
partition $Y=I\sqcup J$, such that
\begin{equation}\label{csflagTL}
 P_k(W;S\cup I)P_l(W;S\cup J)
   =\sum_{M\ \mathrm{compatible\ with}\ I,J}d_M(W),
\end{equation}
when $|I|=k-|S|$ and $|J|=l-|S|$. The matchings are exposed,
with $k-|S|$ arcs and $l-k$ unmatched vertices.
\end{lemma}

\begin{proof}
For $l=0$ the empty matching has weight $1$. For $l>0$, eliminate the entries below the diagonal of $W$ one column at a time, using row $i-1$ to clear row $i$ in decreasing order of $i$. Before column $j$ is cleared, row $i\ge j$ is the original row $i$ plus a combination of the preceding $j-1$ rows and vanishes in columns $1:j-1$. Expanding along this row gives its $j$th entry as
\[
 \frac{P_j(W;\,i-j+1:i)}{P_{j-1}(W;\,i-j+1:i-1)}>0 .
\]
Thus $W=L\binom{U}{0}$, where $U$ is upper triangular with positive diagonal and $L$ is a product of matrices $I+tE_{i,i-1}$ with $t>0$.

Draw each factor of $L$ on $m$ horizontal lines, with one extra directed edge of weight $t$ from level $i-1$ to level $i$. Both its endpoints are inside the strip for that factor. The horizontal edges have unit weight. Place the rightmost factor first from the left, and give the entrance edge at level $j\le l$ weight $U_{jj}$. In one factor a path may stay on its level or take the extra edge to a vacant level. Cauchy--Binet over the factors therefore expresses $P_h(W;T)$ as the sum of the weights of vertex-disjoint path families from $1:h$ to $T$.

Superpose a red family from $1:k$ to $S\cup I$ and a blue family from $1:l$ to $S\cup J$. Discard edges used by both families and orient each remaining used edge forward for red and backward for blue. At each interior vertex the incoming and outgoing counts agree. Since its degree is at most three, it has either no edges left or one edge in each direction. The components are thus cycles and disjoint chains. Each chain starts at a sink in $J$ and ends at a sink in $I$ or at one of the $l-k$ extra blue sources. The sink-to-sink chains pair opposite colors on $Y$. By planarity no such pair can enclose a sink joined to a source, so these pairs form an exposed matching $M$.

For a fixed $M$, changing from one compatible coloring to another reverses the colors at the ends of some of its pairs. Swap red and blue along the corresponding chains. At each interior vertex on a chosen chain this exchanges all singly used incident edges, and every doubly used edge still belongs to both families. The resulting families are again vertex-disjoint within each color. The exchange is its own inverse and preserves the product of the weights. Thus the total weight $d_M(W)$ of the pairs of families inducing $M$ is nonnegative and the same for every compatible coloring. Summing over $M$ gives \eqref{csflagTL}.
\end{proof}

\subsection{Positivity of the common-flag Wronskians}

\begin{proof}[Proof of Theorem~\ref{csmother}]
First take $q\ge s$. Set $l=q-1$, $k=s-1$, $N=q+s$ and
$c_N=\prod_{j=0}^{N-1}j!$. Fix $x>0$ and put
$u=1/x$, $v_i=(x+\mu_i)^{-1}$,
$D_T=\prod_{i\in T}(x+\mu_i)$ and
$\Delta_T=\prod_{i<j,\ i,j\in T}(\mu_j-\mu_i)$.
The matrix $V_{ij}=v_i^2(u-v_i)^{j-1}$ satisfies
\[
 P_h(V;T)=\Delta_T/D_T^{h+1}\qquad(|T|=h)
\]
by the Vandermonde formula. Write
$w_h(T)=P_h(W;T)P_h(V;T)>0$.

For a signed measure $\nu$ on $[0,u]$, let $R_{q,s}(\nu)$ be
$\det[h_{s+i-j}]_{i,j=0}^{q-1}$ with the power sums of
Lemma~\ref{translation}. Put
$\sigma_x=b\delta_0+\rho_x$ as in the proof of
Lemma~\ref{Pickkernel}, so $\sigma_x([0,u])=\psi'(x)$, and
write $\mathrm d\eta_i=\mathbf1_{[0,v_i]}\,\mathrm dv$.
We claim that
\begin{equation}\label{opCauchyExpansion}
 \frac{W(\mathcal F_q,\Phi\mathcal F_s)(x)}{c_N\Phi(x)^s}
 =\sum_{|T|=l,\ |U|=k}w_l(T)w_k(U)
 R_{q,s}\!\left(\alpha\,\mathrm dv+\sigma_x
                          +\sum_{i\in T}\eta_i-\sum_{i\in U}\eta_i\right).
\end{equation}
Indeed, divide \eqref{paired} by the node Vandermonde and coalesce
the nodes. For $I=U\setminus T$, $J=T\setminus U$, factor the
positive function $D_I$ from every column of the numerator. This
contributes $D_I^N$ to its Wronskian and leaves the multiplier
$F=\Phi D_J/D_I$. Equation~\eqref{Wronskian}
gives $c_NF^s R_{q,s}$ with the measure displayed, since each
$\log(x+\mu_i)$ contributes $\eta_i$. Finally,
$D_I^qD_J^s/D_{T\cup U}^{N}=D_T^{-q}D_U^{-s}$,
which gives exactly the two factors from $V$.

Group \eqref{opCauchyExpansion} by $S=T\cap U$ and
$Y=T\mathbin\triangle U$. Apply Lemma~\ref{csnetwork} to each
of $W,V$ and fix one matching from each admitting a common coloring.
Draw them above and below $Y$, with $T\setminus U$ positive.
Their union consists of cycles and paths, including isolated vertices.
In each component in its own increasing order every arc joins
opposite parities: all vertices of that component inside an arc
are paired internally on that side by noncrossing and exposure.
Both path ends are positive, so its even length leaves one end
unmatched on each side. On either side the vertices before the
unmatched vertex pair internally, giving it odd rank. Thus each
path has the fixed colors $+,-,\ldots,+$ in increasing order,
while each cycle admits both alternating colorings independently.

For a component $H=\{i_1<\cdots<i_a\}$, the measure
$\eta_H=\sum_{j=1}^a(-1)^{j-1}\eta_{i_j}$ is nonnegative
since its intervals are nested. Put
\[
 \eta_0=(\alpha-q+1)\,\mathrm dv+\sigma_x+\sum_{H\text{ path}}\eta_H.
\]
For this pair of matchings the sum over colorings is the sum, at
all independent $z_C=\pm1$, of the polynomial
\[
 R_{q,s}\!\left((q-1)\,\mathrm dv+\eta_0+
                                     \sum_{C\text{ cycle}}z_C\eta_C\right).
\]
Its coefficients are nonnegative by Lemma~\ref{translation}.
Averaging retains precisely the monomials of even degree in each
variable, proving nonnegativity of every group. More generally,
weighting the average by any product of the signs selects the
corresponding parities and still gives a nonnegative result, for
every determinant of Lemma~\ref{translation}. The unweighted
average is at least $b_K\eta_0([0,u])^{d}/d!$ in its notation.

Choose now $U\subseteq T$ of sizes $k,l$. Its group consists of
one summand, whose measure exceeds $(q-1)\,\mathrm dv$ by the mass
$L=(\alpha-q+1)/x+\psi'(x)+\sum_{j\in T\setminus U}v_j$.
It follows that
\begin{equation}\label{csuniformstrict}
 W(\mathcal F_q,\Phi\mathcal F_s)(x)
 \ge c_N\Phi(x)^s w_l(T)w_k(U)L^{qs}/(qs)! .
\end{equation}
This is strict positivity when $q>s$, $\alpha>q-1$, or
$\psi$ is nonconstant.

Suppose $q=s\ge2$, $\Phi(y)=cy^{q-1}$, and $m\ge q$.
Choose $S$ of size $q-2$ and $Y=\{i,j\}$ disjoint from $S$.
Put $a=\mu_i$ and $b=\mu_j$. The two summands
in this group of \eqref{paired} have equal positive prefactors.
After removing $c^q$ and coalescing, the first numerator is the
Wronskian of $(y+a)\mathcal P_q$ and $y^{q-1}(y+b)\mathcal P_q$,
where $\mathcal P_q=(1,y,\ldots,y^{q-1})$. Its polynomials have
degree at most $2q-1$, so the Wronskian equals $c_{2q}$ times their
coefficient determinant. The lowest and highest $q-1$ degrees
give triangular blocks with diagonals $a$ and $1$, leaving the
two columns $(a,1)^T,(b,1)^T$ in the middle. Adding the term
with $a,b$ exchanged gives, apart from positive factors,
$(a-b)(a^{q-1}-b^{q-1})>0$.
If $m=q-1$, clearing all denominators turns $\mathcal F_q$ into a
basis of $\mathcal P_q$ by invertibility of the residue matrix.
Together with $y^{q-1}\mathcal P_q$ the span has degrees $0$ to
$2q-2$, giving the asserted rank and zero Wronskian. For $q=s=1$
a constant $\Phi$ makes the two columns proportional for every $m$.

For $q<s$ set $y=1/x$,
$\mu'_a=1/\mu_{m+1-a}$,
$W'_{aj}=(-1)^{j-1}W_{m+1-a,j}/\mu_{m+1-a}^2$, and
$\Phi^\#(y)=1/\Phi(1/y)=y^\alpha e^{-\psi(1/y)}$.
The new flag minors are positive because the row reversal and
column signs cancel; $-\psi(1/y)$ is operator monotone.
In each column $-H_j(1/y)$, removing its constant leaves
$(-1)^j(-H'_j(y))$, by partial fractions.
Divide the evaluation rows by $\Phi(x_i)$, exchange the blocks,
and reverse the reciprocal nodes. The net sign is positive since
$qs+\binom N2+\binom q2+\binom s2$ is even. Thus the original
determinant equals $\prod_i\Phi(x_i)$ times the determinant for
$(\mathcal F'_s,\Phi^\#\mathcal F'_q)$ at the reversed reciprocal
nodes. Coalescing gives
\[
 W(\mathcal F_q,\Phi\mathcal F_s)(x)
 =\Phi(x)^N x^{-N(N-1)}W(\mathcal F'_s,\Phi^\#\mathcal F'_q)(1/x)>0 .
\]

A prefix in the first block alone has positive Wronskian by the
confluent Cauchy--Vandermonde formula \eqref{cv}. Apply the result
just proved to a mixed proper prefix $(q,v)$, $v<s$. Its only possible
exception would require $q=v$ and $\alpha=q-1$, whereas then
$\alpha\ge s-1>q-1$. Lemma~\ref{ECT} gives all evaluation signs
in the strict cases.
\end{proof}

\subsection{Boundary minors of the matrix function}
\label{csedgesection}

\begin{proof}[Proof of Theorem~\ref{edge}]
In \eqref{properblockGram} both the weights and the $v$-prefix minors
are positive. Apply Theorem~\ref{csmother} with $(q,s)=(r,p)$:
when $r=p\ge2$ the spectral flag has $m=n-2\ge2p-2\ge p$ poles,
so even the constant-multiplier endpoint is strict. For $p=r=1$
constant $\Phi$ gives a scalar matrix. Use transposition and reversal
for the other boundaries. Taking $r=n-p$ gives the stated internal far
corners since $\max(p,n-p)-1\le n-3$ there.
\end{proof}

\begin{corollary}\label{csfinite}
Let $N\ge2$, and let $\psi$ be real analytic on $(0,\infty)$.
Suppose that, on every set of $N$ distinct positive numbers, $\psi$
agrees with a real operator-monotone function on $(0,\infty)$.
Let $p,r\ge1$, $p+r\le N$, $n\ge p+r$, and $\tau\ge0$.
For $\Phi(x)=x^\alpha e^{\tau\psi(x)}$ with
$\alpha\ge\max(p,r)-1$, every $A\in\TP_n$ satisfies
\[
             \det\Phi(A)[1:p,r+1:r+p]>0,
\]
except when $\alpha=0$ and $\tau\psi$ is constant; then $p=r=1$
and the indicated entry is zero.  The same assertion holds after
transposition and simultaneous reversal.
In particular, if $n\ge N$ and $\alpha\ge N-2$, the leading and
trailing $N\times N$ blocks of $\Phi(A)$ are TP, except when
$N=2$, $\alpha=0$, and $\tau\psi$ is constant.
The corresponding nonnegative assertions hold for nonsingular TN
matrices.
\end{corollary}

\begin{proof}
In each summand of \eqref{properblockGram}, extend its $p+r$ nodes
to $N$ distinct positive nodes and interpolate $\psi$ by an
operator-monotone $\chi$. Substituting $\tau\chi$ for $\tau\psi$
leaves that summand unchanged. Interpolation on any two nodes,
extended to $N$, makes $\psi$ nondecreasing. If $\tau\psi$ is
nonconstant, analyticity makes it strictly increasing, so every
$\tau\chi$ used here is nonconstant. Theorem~\ref{csmother} and
the same pole count as in the preceding proof give the stated signs. For a leading $N\times N$ block these boundary minors
and their transposes cover all off-diagonal initial minors;
Lemma~\ref{principalbounds} covers the leading principal ones.
Apply the initial-minor criterion and reversal. TP approximation
and continuity of the fixed analytic $\Phi$ give the TN assertion.
\end{proof}

\begin{corollary}\label{csedgepowers}
Let $p,r\ge1$, $n\ge p+r$, and $\gamma=\max(p,r)-1$.
The exponents $t>0$ for which $\det A^t[1:p,r+1:r+p]>0$ for all
$A\in\TP_n$ form
\begin{equation}\label{edgeexponents}
 \mathbb Z_{>0}\cup\bigl([\gamma,\infty)\cap(0,\infty)\bigr).
\end{equation}
For every nonempty compact subset of an excluded unit interval,
one $A\in\TP_n$ makes this minor negative throughout that subset.
\end{corollary}

\begin{proof}
Theorem~\ref{edge} and integer powers give sufficiency.
Write $C_{a,b}(M)=\det M[1:a,b+1:a+b]$ and let $S$ be the
upper shift of order $a+b$. Then
\[
 C_{a,b}((I+S)^t)
 =\det\!\left[\binom{t}{b+j-i}\right]_{i,j=0}^{a-1}
 =\prod_{j=0}^{a-1}\frac{\binom{t+j}{b}}{\binom{b+j}{b}},
\]
where $\binom th=0$ for $h<0$. The first equality is the binomial
expansion. For the second use \eqref{Wronskian} at $1$ with $q=a$,
$v=b$, $F(x)=x^{-t}$, so $h_j=(-1)^j\binom tj$. The monomial
Wronskian is the Vandermonde on the exponents
$0,\ldots,a-1,-t,\ldots,-t+b-1$; its cross factors cancel
the sign $(-1)^{ab}$ and give the product. For $b-2<t<b-1$ only
the factor $j=0$ is negative. This also gives the entry $t$ for
$a=b=1$, including negative $t$ and $t=0$.

Checker inversion $H\mapsto EH^{-1}E$ exchanges $C_{a,b}(H^t)$
and $C_{b,a}((EH^{-1}E)^t)$ up to the positive divisor $\det H^t$
by Jacobi's identity. Thus assume $r\ge p$. For compact
$K\subset(k,k+1)$, $k\in\{0,\ldots,\gamma-1\}$, put $s=k+2\le r$. The
displayed product for $(a,b)=(p,s)$ is negative on $K$. By joint
continuity and compactness this persists under TP approximation.

To increase $a$ by one from an order-$m$ witness $H$ for $(a,b)$,
put $y=(L,\ldots,L^m)^T$. For large $L$, $[H\mid y]$ is TP by
last-column expansion, while
$\det[H^t[1:a+1,b+1:m]\mid y[1:a+1]]<0$ uniformly on $K$:
its leading coefficient in $L$ is $C_{a,b}(H^t)$.
For large $\rho$ the matrix
\[
 B_\rho=\begin{pmatrix}H&(\rho I-H)y\\0&\rho\end{pmatrix}
\]
is nonsingular TN: divide its last column by $\rho$ in the top
$m$ rows and use the TP limit $[H\mid y]$, then expand along the
bottom row. Conjugation of $H\oplus\rho$ by
$\left(\begin{smallmatrix}I&y\\0&1\end{smallmatrix}\right)$ shows
that the last column of $B_\rho^t$ above its diagonal, divided by
$\rho^t$, tends uniformly to $y$ on $K$. Thus
$C_{a+1,b}(B_\rho^t)<0$ on $K$; approximate by TP matrices.
Checker inversion gives the step increasing $b$. Iterate from
$(p,s)$ to $(p,r)$ and pad by a direct sum with an identity and
TP approximation if $n>p+r$.
\end{proof}

\section{Minors of the internal orders}\label{imsection}

A cyclic interval of size $p$ consists of $p$ consecutive indices in
the cyclic order, listed in increasing order for taking minors.

\begin{lemma}\label{cyclicflag}
If the maximal minors of a real $p\times n$ matrix are nonnegative
and those on cyclic intervals are positive, then all are positive.
\end{lemma}

\begin{proof}
Suppose a zero minor exists. Choose one whose set $I$ has a shortest cyclic cover
$[a,c]$ of globally minimum length $\ell$ in steps. Among ties maximize
the sum of distances from $a$ to $I$. Both endpoints are in $I$, and
$p\le\ell\le n-2$ since cyclic minors are positive and a cover can omit one
point outside $I$. Take $b\in(a,c)\setminus I$, $d=c+1$ in cyclic order,
and $S=I\setminus\{a,c\}$. The three-term Pl\"ucker identity gives
\[
 0=\Delta_{Sac}\Delta_{Sbd}
  =\Delta_{Sab}\Delta_{Scd}+\Delta_{Sad}\Delta_{Sbc},
\]
where $\Delta_{Suv}$ is the minor on $S\cup\{u,v\}$. The signs are valid
across the cut because rotating $k$ selected columns changes a minor
by $(-1)^{k(p-1)}$, and every product uses the same multiset of columns.
The two shorter covers give $\Delta_{Sab},\Delta_{Sbc}>0$, so
$\Delta_{Scd}=0$. In this new set the successor $a'$ of $a$ in $I$
must be $a+1$, since otherwise $[a',d]$ is shorter than $[a,c]$.
But relative to $a+1$ the distance sum increases by $\ell-p+1>0$.
\end{proof}

\begin{proof}[Proof of the internal-minor assertion in Theorem~\ref{iminternal}]
Normalize $A\in\TP_n$ by positive diagonal similarity and Perron
scaling so that $A\boldsymbol1=\boldsymbol1$. If $\rho$ is its original
Perron root, replace $\psi(x)$ by $\psi(\rho x)-\psi(\rho)$, which leaves
the minor signs invariant. Along the path
\[
 F_s=h_s(A),\qquad h_s(x)=x^{(1-s)(n-3)+s\alpha}e^{s\psi(x)},\qquad 0\le s\le1,
\]
we have $F_s\boldsymbol1=\boldsymbol1$. Fix $2\le p\le n-2$ and write
$F$ for the matrix at a first zero of a $p$th minor, if one exists.
All $p$th minors are then nonnegative since $F_0=A^{n-3}$ is TP.

For any equally sized row and column sets, agreement on membership
of either endpoint gives positivity by Theorem~\ref{properblock}.
Apply the theorem directly when both sets omit that endpoint and
use \eqref{imchecker} when both contain it. Now take cyclic sets
$I,J$. If each contains one endpoint and these differ, use
Theorem~\ref{edge}. Otherwise, after transposing and
renormalizing if needed, $I$ avoids both endpoints and
$J=(1:a)\cup(n-b+1:n)$ with $a,b>0$ and $a+b=p$.

In the augmented matrix $[F[1:n,J],\boldsymbol1]$, take the maximal
minors on $\{1\}\cup I$ and $I\cup\{n\}$. Both have sign $(-1)^b$.
Indeed, replacing the last column by $F\boldsymbol1$ gives a sum of
ordered positive minors, each sharing an endpoint and each requiring
$b$ column exchanges. If $F[I,J]$ were singular, take a relation among
its rows and normalize one coefficient to one. Replace the corresponding
row of both augmented minors by that combination. Expansion along
this row now gives opposite signs from the last-column cofactors,
which have positive unsigned minors by the same endpoint criterion.
This contradicts the common strict sign.

Apply Lemma~\ref{cyclicflag} first with each cyclic row set fixed
and then, after transposing, with each column set fixed. All $p$th
minors are positive, contradicting a first zero. For orders $1,n-1$,
Theorem~\ref{properblock} makes the leading and trailing $(n-1)$-blocks
of $F=\Phi(A)$ and $E F^{-1}E$ TP. By \eqref{imchecker} this leaves
only \eqref{fourcorners}, and $\det F=\prod_i\Phi(\lambda_i)>0$.
The principal-minor assertion for nonsingular TN inputs and the
sharpness of $\alpha\ge n-3$ were proved in Section~\ref{principalSection}.
\end{proof}

\subsection{The exact power exponents at orders two and \texorpdfstring{$n-2$}{n-2}}

\begin{proof}[Proof of Corollary~\ref{imexact}]
Theorem~\ref{iminternal} covers $t\ge n-3$, and positive integer powers
preserve TP. For $p=2,n-2$, Corollary~\ref{csedgepowers} with $r=n-p$
gives a negative upper-right $p$th corner at every noninteger
$0<t<\max(p,n-p)-1=n-3$. Also $\mathcal E_{n,p}=\mathcal E_{n,n-p}$
by the TP bijection $A\mapsto EA^{-1}E$ and \eqref{imchecker}.
\end{proof}

For the upper shift $S_n$ and $m=n-3$, the single TN input $B=I+S_n$
witnesses all the excluded intervals at order two.
\begin{equation}\label{imbinomial}
 \det B^t[1:2,n-1:n]
 =\binom{t}{m+1}^{\!2}-\binom{t}{m}\binom{t}{m+2}
 =\frac{\bigl(t(t-1)\cdots(t-m+1)\bigr)^2(t-m)(t+1)}
 {(m+2)((m+1)!)^2}.
\end{equation}
For $A\in\TP_n$ and $t\le-(n-3)$, the identity
$A^t=E(EA^{-1}E)^{-t}E$ gives
$(-1)^{\sum I+\sum J}\det A^t[I,J]>0$ outside \eqref{fourcorners}.

\subsection{Joint signs of the four corner minors}

Write $\iota(A)=EA^{-1}E$ and $\Phi^\#(x)=1/\Phi(1/x)$. By
\eqref{imchecker}, the upper pair of signs is that of the $(1,n)$
entries of $F=\Phi(A)$ and $\iota(F)=\Phi^\#(\iota(A))$.

\begin{lemma}\label{orientedcorners}
Let $m=n-2\ge2$, $r=m+1$, and let $\Phi(x)=x^\alpha e^{\psi(x)}$
with $m-1\le\alpha<m$ and $\psi$ real operator monotone on $(0,\infty)$.
Let $\lambda$ be $n$ positive nodes, allowing repetitions. Put
\[
 D_j=[x^j\Phi]_{\lambda},\quad J_j=(-1)^r[x^j/\Phi]_{\lambda},\quad
 q=\Big(\prod_i\lambda_i\Big)^{m/n},\quad c=(m-\alpha)/\alpha .
\]
For $0\le j\le m$, $D_j\ge0$ if $j>0$, $J_j\ge0$ if $j<m$, and
$D_m,J_0>0$. Moreover,
\begin{equation}\label{extremecoef}
 -J_m<cqJ_0,\qquad -D_0<cq^{-1}D_m.
\end{equation}
If $A$ is nonsingular TN with spectrum $\lambda$ and
$E_A(x)=\adj(xI-A)_{1n}=\sum_{j=0}^m e_j x^j$ has $e_0,e_m>0$, then
\begin{equation}\label{cornersep}
 \frac{F_{1n}}{D_m}+\frac{\iota(F)_{1n}}{qJ_0}>0.
\end{equation}
Thus each orientation has at least one positive corner for $A\in\TP_n$.
\end{lemma}

\begin{proof}
By the Pick representation, $\psi'$ is completely monotone, so
$L=e^{-\psi}$ is completely monotone by induction from $L'=-\psi'L$.
Thus $x^{j-\alpha}L$ is completely monotone for $j<m$, and \eqref{hgdivided} gives $J_j\ge0$ with $J_0>0$.
Inverting the nodes in the Lagrange weights gives
$J_j(\Phi^\#,\lambda^{-1})=(\prod_i\lambda_i)D_{m-j}(\Phi,\lambda)$,
including repetitions by continuity. This yields the claimed $D_j$ signs.

Put $W=(-1)^r(x^{-\alpha}L)^{(r)}$, $p=\alpha+r$, and $\eta=x^pW$.
Leibniz's rule gives $-(-1)^r(x^{m-\alpha}L)^{(r)}<cx^mW$. Indeed,
writing $\varepsilon=m-\alpha\in(0,1]$, the term with no derivatives on
$x^\varepsilon$ is nonpositive. For $k\ge1$ derivatives on this power,
the ratio to the corresponding coefficient in $x^mW$ is at most
$c\prod_{\ell=1}^{k-1}(\ell-\varepsilon)/(\alpha+\ell)\le c$.
The bound for $k=r$ is strict, and that term contains $L>0$.
For $a_j=(-1)^jL^{(j)}$, the same differentiation gives
\[
 \eta'=(1-\alpha)\sum_{j=0}^{r-1}\binom rj
 (\alpha)^{\overline{r-j-1}}x^ja_{j+1}-x^ra_{r+1}\le0,
\]
where the bar denotes a rising factorial. Let $X$ be the convex
combination of the nodes with uniform simplex weights. Formula
\eqref{hgdivided} and opposite monotonicity of $X^m$ and $\eta(X)$
under the weight $X^{-p}$ give
\[
 -\frac{J_m}{J_0}<c\frac{\mathbb E[X^{m-p}\eta(X)]}{\mathbb E[X^{-p}\eta(X)]}
 \le c\frac{\mathbb E X^{m-p}}{\mathbb E X^{-p}}\le cq.
\]
The last inequality compares power means of $X^{-1}$ at $p-m<n\le p$,
using $\mathbb E X^{-n}=(\prod_i\lambda_i)^{-1}$ by the formula for $1/x$.
The reciprocal-node identity gives the second bound of \eqref{extremecoef}.

By \eqref{shared} and TP approximation the $e_j$ are nonnegative and
strictly positive for TP inputs. Spectral interpolation and Jacobi's
identity give $F_{1n}=\sum_j e_jD_j$ and
$\iota(F)_{1n}=\sum_j e_jJ_j$, with repetitions covered by continuity.
In \eqref{cornersep}, the coefficients
of $e_0,e_m$ exceed $q^{-1}(1-c),1-c$, respectively, and the others
are nonnegative. Since $c\le1$, this proves \eqref{cornersep}.
Transposition treats the lower pair.
\end{proof}

\begin{theorem}\label{fourpatterns}
Let $n\ge4$ and $\Phi(x)=x^\alpha e^{\psi(x)}$ with $\alpha\ge n-3$
and $\psi$ real operator monotone. If $\alpha\ge n-2$, all four
corners are positive on $\TP_n$. Otherwise, fix a simple positive
spectrum $\lambda$ and let $a,b$ indicate whether $[\Phi]_\lambda<0$
and $[\Phi^\#]_{\lambda^{-1}}<0$, respectively. For $A\in\TP_n$ with
spectrum $\lambda$, group the minors in \eqref{fourcorners} into
the upper-right and lower-left pairs, with the entry first. Each
pair independently realizes $(+,+)$, also $(-,+),(0,+)$ if $a=1$,
and $(+,-),(+,0)$ if $b=1$. These are exactly $(1+2a+2b)^2$ patterns.
For symmetric TP inputs with the same spectrum, the two pairs agree
and all $1+2a+2b$ common pairs occur. In each class, the listed patterns
occur in every neighborhood of any diagonal matrix with spectrum $\lambda$.
\end{theorem}

\begin{proof}
The range $\alpha\ge n-2$ follows from Theorem~\ref{multipliers} and
Cauchy--Binet for the constant multiplier at equality. Otherwise,
Lemma~\ref{orientedcorners} excludes two nonpositive corners of one
orientation. Its expansions, with $e_0,e_m>0$, also show that
$D_0\ge0$ forces $F_{1n}>0$ and $J_m\ge0$ forces $\iota(F)_{1n}>0$.

Let $\Lambda=\diag(\lambda_1,\ldots,\lambda_n)$ and $S$ be the upper shift.
The triangular matrix $U(z)=(I+(1-z)S)(I-zS)^{-1}\Lambda$ is TN
for $0\le z\le1$. Inverting the upper bidiagonal pencil
$(I-zS)(sI-U(z))=(sI-\Lambda)-S(zsI+(1-z)\Lambda)$ gives
\[
 E_{U(z)}(s)=\lambda_n\prod_{j=2}^{n-1}(zs+(1-z)\lambda_j).
\]
If $D_0<0$, then $\Phi(U(z))_{1n}/(1-z)^m$ is a polynomial in
$z/(1-z)$ with negative constant, nonnegative other coefficients and
positive leading coefficient. Thus the corner entry $z\mapsto\Phi(U(z))_{1n}$ has a unique simple zero in $(0,1)$.
By \eqref{cornersep} the second corner is positive at the zero and
for smaller positive $z$, and both are positive just to its right.
Apply the same argument to $\Phi^\#$ on the reciprocal spectrum for
the second-corner patterns. The model allows any ordering of the
positive diagonal. Writing that diagonal explicitly, we have
$\iota(U(z;\Lambda^{-1}))=\Lambda U(1-z;\Lambda)\Lambda^{-1}$
since $ESE=-S$. This positive diagonal similarity preserves the signs
and simple zeros. If $a=b=0$, both corners are positive for any
interior $z$ by the same expansions.

For independently chosen upper models $U(u),V(v)$ with diagonal
$\Lambda$, put $U_\delta=U(u)e^{\delta S}$,
$L_\delta=(V(v)e^{\delta S})^T$, $D_\delta=\diag(1,\delta,\ldots,\delta^{n-1})$,
and $A_\delta=D_\delta L_\delta\Lambda^{-1}D_\delta^{-1}U_\delta$.
For $\delta>0$, the limit $e^{\delta S}=\lim_{N\to\infty}(I+\delta S/N)^N$ is TN. Given $J=\{j_1<\cdots<j_k\}$, factor the powers of $\delta$ from the rows and columns and multiply column $r$ by $(j_r-1)!$. The resulting rows are the falling factorials of degrees $0,\ldots,k-1$ evaluated at $j_r-1$. Thus the Vandermonde formula gives
\[
 \det(e^{\delta S})[1:k,J]=
 \frac{\delta^{\sum_{r=1}^k(j_r-r)}\prod_{a<b}(j_b-j_a)}
      {\prod_{r=1}^k(j_r-1)!}>0 .
\]
Cauchy--Binet on $1:k$
passes this to the first $k$ rows of $U_\delta$ and the first $k$
columns of $L_\delta$, and gives $A_\delta\in\TP_n$.

Read the upper signs from $A_\delta$ and the lower signs from
$D_\delta^{-1}A_\delta D_\delta$. The two families extend analytically
at $\delta=0$ with limits $U(u),V(v)^T$. Vary the diagonal entries of
$\Lambda$ and the parameters of any requested zero corners. At
$\delta=0$ the simple eigenvalues depend on the diagonal alone,
with identity Jacobian, whereas the derivative of the zero constraints
in their own parameters is diagonal and invertible by the simple
roots. The implicit function theorem fixes the prescribed spectrum
and zeros, and strict signs persist for small $\delta>0$.

Set $C_\delta=D_\delta^{-1/2}A_\delta D_\delta^{1/2}$. Its diagonal
tends to the prescribed one. Boundedness of $A_\delta$ and
$D_\delta^{-1}A_\delta D_\delta$ gives $(C_\delta)_{ij}
=O(\delta^{|i-j|/2})$ for $i\ne j$. Thus $C_\delta\to\Lambda$,
and positive diagonal similarity preserves the input and output minor signs.

For a symmetric realization use the same model in both directions
and impose only the upper zero, if required, when fixing the spectrum.
Then $L_\delta=U_\delta^T$ and
\[
 C_\delta
 =D_\delta^{1/2}U_\delta^T(\Lambda D_\delta)^{-1}
                    U_\delta D_\delta^{1/2}
\]
is symmetric and has identical upper and lower pairs.
\end{proof}

As the spectrum varies, each indicator can equal one exactly when the
corresponding $(n-1)$st derivative is negative somewhere, by
\eqref{hgdivided} and coalescence. The union has size $1,9,17$, or $25$.
The case $17$ is the union of the two nine-pattern sets when both
indicators occur but never together. With $m=n-2$, examples with $1$ and $9$ patterns
are $\Phi=x^m$ and $\Phi=x^{m-1}(1+x)$, respectively, since the latter
has vanishing $(m+1)$st derivative and its dual $x^m/(1+x)$ has strictly
negative $(m+1)$st derivative by polynomial division.

\begin{example}
In order four, $\Phi(x)=x^{1+d}e^{t(x)}$, where $t(x)=(x-x^{-1})/4$ and $0\le d\le1/8$, has exactly $17$ patterns over all spectra.
Here $t$ is operator monotone and $\Phi^\#=\Phi$. Put $H=d+(x+x^{-1})/4$ and $K=H^2-d/H>0$, since $H^3\ge(d+1/2)^3>d$.
Using $x^3\partial_x^3=D(D-1)(D-2)$ with $D=x\partial_x$ gives
\[
 \Phi'''=x^{-3}\Phi H(K+3t),\qquad -(x^2/\Phi)'''=x^{-1}\Phi^{-1}H(K-3t).
\]
At $x=1/2$, $H\le3/4$ and $3t=-9/8$, so $\Phi'''(1/2)<0$. Thus $D_0<0$ on nearby simple spectra, and reciprocating gives $J_2<0$ on other spectra.
For a fixed spectrum let $X$ be its convex combination with uniform simplex weights. With $\mathrm d\nu=x^{-1}\Phi^{-1}H\,\mathrm d\mathbb P_X$ and $W=x^{2d}e^{2t}$, \eqref{hgdivided} gives
\[
 6D_0=\int W(K+3t)\,\mathrm d\nu,\qquad 6J_2=\int(K-3t)\,\mathrm d\nu.
\]
If $D_0\le0$, the $W$-weighted mean of $t$ is negative. Since $W$ and $t$ increase, their covariance under normalized $\nu$ is nonnegative, so the unweighted mean of $t$ is also negative and $J_2>0$.
By Theorem~\ref{fourpatterns}, the two nine-pattern sets therefore meet only at the all-positive pattern.
\end{example}

\begin{corollary}\label{fournecessary}
Let $n\ge4$, $n-3<t<n-2$, and let $\Lambda$ be diagonal with distinct
positive entries. Among TP matrices $A$ with this spectrum in any
neighborhood of $\Lambda$, $A^t$ realizes all $25$ patterns in
Theorem~\ref{fourpatterns}. Any one of the four corners can be the sole
negative minor of $A^t$. Symmetric TP inputs realize all five common pairs.
\end{corollary}

\begin{proof}
For $\Phi=x^t$ we have $\Phi^\#=\Phi$ and $\Phi^{(n-1)}<0$
on $(0,\infty)$. Formula \eqref{hgdivided} gives negative divided differences
on both the prescribed and reciprocal spectra. Now apply
Theorem~\ref{fourpatterns}; for a sole negative minor choose the
other three corners positive and use Theorem~\ref{iminternal}.
\end{proof}

\begin{corollary}\label{fixedpower}
Let $\Lambda$ be diagonal of order $n\ge2$ with distinct positive entries.
Among TP matrices $A$ of this spectrum in any neighborhood of $\Lambda$,
the exponents for which every $A^t$ is TN form $T=\mathbb Z_{\ge0}\cup[n-2,\infty)$,
and those giving TP form $T\setminus\{0\}$. The same holds for symmetric
inputs. For each nonempty compact $K\subset\mathbb R\setminus T$, one
such symmetric $A$ makes at least one entry of $A^t$ negative for each $t\in K$.
Another has $A^t\in\TP_n$ for every $t>0$.
\end{corollary}

\begin{proof}
Corollary~\ref{sets} gives sufficiency, and $A^0=I$. Write
$\Lambda=\diag(\lambda_1,\ldots,\lambda_n)$ and $B=(I+S)\Lambda$.
Its bidiagonal resolvent gives, for $1\le d<n$,
\[
 (B^t)_{1,d+1}=(\lambda_2\cdots\lambda_{d+1})
                 [x^t]_{\lambda_1,\ldots,\lambda_{d+1}}.
\]
Take $d=1$ for $t<0$ and $d=\lfloor t\rfloor+2$ for noninteger
$0<t<n-2$. Since $(x^t)^{(d)}<0$, the minimum of these entries is
uniformly negative on $K$ by compactness. Using $B=U(0)$ for both
models in the proof of Theorem~\ref{fourpatterns} fixes the spectrum
along $A_\delta\to B$ in every order. Joint continuity preserves the
negative minimum on $K$, and $(C_\delta)^t$ has the same entry signs
as $(A_\delta)^t$.
For the other choice, perturb $\log\Lambda$ to a Jacobi matrix
$J_\varepsilon$ with adjacent entries $\varepsilon>0$. The implicit
function theorem fixes its spectrum by varying the diagonal, since
at $\varepsilon=0$ the labeled eigenvalues are the diagonal entries.
Then $e^{J_\varepsilon}\to\Lambda$ and Schwarz's theorem
(Section~\ref{stsection}) gives $e^{tJ_\varepsilon}\in\TP_n$ for every $t>0$.
\end{proof}

\begin{remark}
Let $A$ be nonsingular TN with $d=d(A)\ge4$ and $\alpha\ge d-3$.
For an order-$d$ block $C$, write $B_C=\Phi(A[C,C])$ with local indices.
Then $\Phi(A)$ is TN if and only if the upper pair of $B_C$ is
nonnegative for each order-$d$ block in $\mathcal U(A)$ and the lower
pair is nonnegative for each order-$d$ block in $\mathcal L(A)$.
Smaller blocks are TN by Theorem~\ref{multipliers}. For an upper block
that is not also a lower block, a lower zero cut makes both lower
corners vanish by Lemma~\ref{zerocuts}, with the transposed statement
for lower blocks. Thus Theorem~\ref{iminternal} by approximation and
Proposition~\ref{localreduction} give the criterion.
\end{remark}

\section{Endpoint flag minors}
\label{endpointsec}

For $1\le p<n$, an endpoint flag minor of order $p$ means a minor
whose row set or column set is $1:p$ or $n-p+1:n$.
For integer powers, Alseidi and Garloff considered the first power
at which all minors of a given order are positive
\cite[Definition~3]{AG21}.

\begin{theorem}\label{endpointflags}
Let $n\ge2$, $1\le p<n$, $\gamma=\max(p,n-p)-1$, and $A\in\TP_n$.
For $\Phi(x)=x^\alpha e^{\psi(x)}$ with $\alpha\ge\gamma$ and
$\psi$ real operator monotone on $(0,\infty)$, every $p$th minor
satisfies
\begin{equation}\label{endpointpositive}
 \det\Phi(A)[I,J]>0\quad\text{if $I$ or $J$ is $1:p$ or $n-p+1:n$.}
\end{equation}
The sole exception is $\alpha=0$ and constant $\psi$: then $n=2$,
$p=1$, and $\Phi(A)$ is a positive scalar matrix.
For $\Phi(x)=x^t$, the positive exponents for which
\eqref{endpointpositive} holds for every $A\in\TP_n$ form the set
\eqref{edgeexponents} for this $\gamma$, including all $t>0$
when $n=2$.
\end{theorem}

\begin{proof}
Assume $\alpha>0$ or $\psi$ is nonconstant, put $q=n-p$, and consider
the first $p$ rows. Among their minors on
cyclic intervals, the leading principal one is positive by
Lemma~\ref{principalbounds}. A column interval $r+1:r+p$, $r\ge1$,
is covered by Theorem~\ref{edge} since $r\le q$. A wrapping interval
has the form $J=(1:a)\cup(n-b+1:n)$ with $a+b=p$, $a,b>0$.
For $E=\diag(1,-1,\ldots)$, coordinate reversal $R$, and
$C=R(EA^{-1}E)^TR\in\TP_n$, complementary minors give
\begin{equation}\label{endpointchecker}
 \frac{\det\Phi(A)[1:p,J]}{\det\Phi(A)}
 =\det\Phi^\#(C)[1:q,b+1:b+q]>0.
\end{equation}
The checker factors cancel the Jacobi sign, and both index sets
are reversed. The last inequality uses Theorem~\ref{edge} with
$\Phi^\#(y)=1/\Phi(1/y)=y^\alpha e^{-\psi(1/y)}$ and
$\max(q,b)-1\le\gamma$.

Join $x^d$ to $\Phi$ by $h_v(x)=x^{(1-v)d+v\alpha}e^{v\psi(x)}$,
where $d=\max(1,\gamma)$ and $0\le v\le1$. All cyclic minors
above stay positive under the strictness hypothesis, and $A^d$ is
TP. At a first zero, the first $p$ rows would have nonnegative
maximal minors and positive cyclic ones, contrary to
Lemma~\ref{cyclicflag}. Transposition and reversal give the other
three families. In the constant case $\Phi(A)$ is scalar.
Positive integers and the tail give the stated exponents;
for each excluded $t$, Corollary~\ref{csedgepowers} with $r=q$
supplies a negative far corner.
\end{proof}

Fix $3\le p\le n-3$, set $\gamma=\max(p,n-p)-1$, and let $A\in\TP_n$.
Suppose that for every $t\in(\gamma,\gamma+1)$ the $p$th-order minors
of $A^t$ are positive whenever both index sets are cyclic intervals
other than $1:p$ and $n-p+1:n$. Then every $p$th-order minor of $A^t$
is positive for every $t>\gamma$.

Indeed, $A^\gamma$ is TP. At a first zero of a $p$th-order minor
in $(\gamma,\gamma+1)$, all minors of order $p$ would be nonnegative.
The hypothesis and Theorem~\ref{endpointflags} give strict positivity
for every pair of cyclic intervals. For each cyclic row set,
Lemma~\ref{cyclicflag} then gives positivity for every column set.
Fixing any column set and applying the same lemma after transposition
gives positivity for every row set, a contradiction.
For each positive integer $k$, Cauchy--Binet applied to
$A^{t+k}=A^tA^k$ extends positivity to $(\gamma+k,\gamma+k+1)$.
The integer powers at the endpoints are TP.
There are $n-2$ cyclic intervals other than the two endpoint intervals.
When the test is made over all $A\in\TP_n$, transposition reduces it
to at most $\binom{n-1}{2}$ pairs.

\section{Rigidity and strict positivity for preservers}\label{stsection}

We prove Theorem~\ref{strict} by realizing the positive spectrum in
the mixed determinants for the far corners of $f(A)$.

\subsection{Growth of preserving functions}

Put $S=\diag(1,-1,1,\ldots)$ and $\iota(A)=SA^{-1}S$.
As in \eqref{imchecker}, Jacobi's identity makes $\iota$ an
involution of the nonsingular TN and TP classes. Since
$\iota(A)^k=\iota(A^k)$, it also preserves oscillatory matrices.
For every positive analytic $f$ on $(0,\infty)$, define
\begin{equation}\label{checkerdual}
 f^\#(x)=\frac1{f(1/x)},\qquad f^\#(\iota(A))=Sf(A)^{-1}S.
\end{equation}
In particular, $f$ and $f^\#$ satisfy the same preservation
conditions on these classes.

\begin{proposition}\label{rigidity}
Let $n\ge2$, and let $f$ be positive and real analytic on
$(0,\infty)$. Suppose that $f(A)$ is TN for every $A\in\TP_n$.
If $n\ge3$, then either $f(x)=cx^k$ for some $c>0$ and
$k\in\{0,\ldots,n-2\}$, or
\[
 xf'(x)>(n-2)f(x)\qquad(x>0).
\]
Moreover, for every $n\ge2$, if at either zero or infinity
\[
 f(x)\sim cx^\alpha\qquad(c>0,\ \alpha\in\mathbb R),
\]
then either $\alpha\ge n-2$, or
$\alpha=k\in\{0,\ldots,n-3\}$ and $f(x)=cx^k$ for all $x>0$.
\end{proposition}

\begin{proof}
By TP approximation and \eqref{checkerdual}, both $f$ and $g=f^\#$
preserve nonsingular TN matrices. For the upper shift $N$ of order
$n$, the $(1,j+1)$ entry of $u(xI+N)$ is $u^{(j)}(x)/j!$,
where $u=f$ or $g$. Thus their derivatives through order $n-1$
are nonnegative. For $n=2$, monotonicity of $f$ excludes a negative
power exponent at either endpoint, which suffices.

Let $m=n-2\ge1$. For $1\le k\le m$, a nonzero $u^{(k)}$
is positive at some point and nondecreasing thereafter.
Integrating gives $u(x)\ge c_0x^k$ for all sufficiently large
$x$, for some $c_0>0$. If $g^{(m)}\equiv0$, then $g$ is a
polynomial of degree at most $m-1$. Let $\ell$ be its lowest
nonzero degree. The identity $f(x)g(1/x)=1$ gives $f(x)=O(x^\ell)$
as $x\to\infty$, so the growth implication with $k=\ell+1$
gives $f^{(\ell+1)}\equiv0$. Clearing denominators in
$f(x)g(1/x)=1$ now shows that both polynomials are monomials,
and $f(x)=cx^\ell$.

Otherwise $g(x)\ge c_0x^m$ near infinity, so $f(x)=O(x^m)$
near zero. The nonnegative nondecreasing derivatives $f^{(j)}$,
$0\le j\le m$, have finite limits at zero. Taylor's formula
therefore gives $f^{(j)}(0+)=0$ for $j<m$.
Write $Q(x)=xf'(x)-mf(x)$. Then $Q^{(m)}(x)=xf^{(m+1)}(x)$,
and $Q^{(j)}(0+)=0$ for $j<m$ since $f^{(j+1)}$ is bounded near
zero. Integrating gives
\[
 Q(x)=\frac1{(m-1)!}\int_0^x t(x-t)^{m-1}f^{(m+1)}(t)\dd t.
\]
The integral is finite since $f^{(m)}$ has a finite limit at zero.
If $f^{(m+1)}$ vanishes identically, then $f(x)=cx^m$ by the
vanishing initial derivatives. Otherwise its analyticity and
nonnegativity give $Q(x)>0$ for every $x>0$.

When $g^{(m)}\not\equiv0$, we have thus proved that $f(x)/x^m$
is positive and nondecreasing. It cannot tend to infinity at
zero or to zero at infinity, ruling out $\alpha<m$.
\end{proof}

In particular, if $n\ge3$ and a nonmonomial entire $f$ is positive
on $(0,\infty)$ and preserves $\TN_n$, it vanishes at zero to
order at least $n-2$.

\subsection{A TN realization with prescribed spectrum}

We use the triangular flag criterion: an upper triangular matrix
with positive diagonal is TN if every minor on its first $k$ rows
of order $k$ is positive, for each $k$. The lower triangular version
uses the first $k$ columns. Both assertions follow from the
triangular criterion of Fallat and Johnson \cite[Lemma~3.3.4]{FJ}.

\begin{lemma}\label{sttriangular}
Let $C\in\TP_r$ have eigenvalues $0<\eta_1<\cdots<\eta_r$.
There is a unit lower triangular matrix $K$ such that $K^{-1}CK$
is upper triangular and TN, with diagonal $\eta_1,\ldots,\eta_r$
in that order.
\end{lemma}

\begin{proof}
With $S=\diag(1,-1,\ldots)$ of order $r$, put $G=SC^{-1}S\in\TP_r$ and write $G=X\Lambda X^{-1}$, where
$\Lambda=\diag(\eta_1^{-1},\ldots,\eta_r^{-1})$. For each $j$,
the minors $\det X[I,1:j]$ are the coordinates of a Perron
eigenvector of the positive $j$th compound of $G$. They therefore
all have the same strict sign. In particular, $X$ has an LU
factorization $X=LU$ with $L$ unit lower triangular, and
\[
 \det L[I,1:j]=\frac{\det X[I,1:j]}{\det X[1:j,1:j]}>0
 \qquad (|I|=j).
\]
The triangular flag criterion makes $L$ TN. The matrix
$T=L^{-1}GL=U\Lambda U^{-1}$ is upper triangular with diagonal
$\Lambda$. Now $GL=LT$ is the LU factorization of the nonsingular
TN matrix $GL$, so the TN LU theorem of Section~2 gives $T\in\TN_r$. Its checkerboard inverse is $K^{-1}CK$
with $K=SLS$ and has diagonal $\eta_1,\ldots,\eta_r$.
\end{proof}

By realizing arbitrary distinct positive nodes as the spectrum of a
triangular TN matrix $U$ with the same initial spans of rational functions,
we can read the sign of the mixed determinant at those nodes from a
corner minor of $f(U)$.

\begin{lemma}\label{strealize}
Let $r\ge1$, $0<\eta_1<\cdots<\eta_r$, and let $F$ be upper
triangular with positive diagonal and
\[
 \det F[1:k,J]>0\qquad(|J|=k,\ 1\le k\le r).
\]
Suppose that $H=F\diag(\eta_1,\ldots,\eta_r)F^{-1}$ is TN, and put
\[
 \phi_j(z)=\sum_{s=1}^r\frac{F_{js}}{z+\eta_s}
 \qquad(1\le j\le r).
\]
For arbitrary positive numbers $a,d_1,\ldots,d_r,c$, there is an
upper triangular nonsingular TN matrix $U$ of order $r+2$ with that
diagonal, such that
\[
 E_U(z):=\adj(zI-U)_{1,r+2}=\prod_{s=1}^r(z+\eta_s).
\]
Moreover, with
$\psi_j^U(z)=\adj(zI-U)_{j,r+2}/E_U(z)$, there are a vector $v$ and
a unit lower triangular matrix $L_0$ such that
\begin{equation}\label{strealizedflag}
 \psi_1^U=1,\qquad
 (\psi_2^U,\ldots,\psi_{r+1}^U)^T
       =v-L_0^{-1}(\phi_1,\ldots,\phi_r)^T.
\end{equation}
\end{lemma}

\begin{proof}
Let $V(z)=\prod_s(z+\eta_s)$ and $\Lambda=\diag(\eta_1,\ldots,\eta_r)$.
Define a column $w$, a row $g$, and the columns $Q_j$ of $Q$ by
\[
 w_s=(a+\eta_s)^{-1},\qquad
 g_s=\frac{(a+\eta_s)\prod_{h=1}^r(-\eta_s-d_h)}{V'(-\eta_s)},
 \qquad (Q_j)_s=\frac{(-1)^{j-1}w_s}{\prod_{h=1}^j(d_h+\eta_s)}.
\]
The $k$-row determinants of $[Q_1,\ldots,Q_k]$ and
$[Q_1,\ldots,Q_{k-1},w]$ are positive. Clearing row denominators
gives polynomial columns of degrees $k-1,\ldots,0$ in the first
case and $k-2,\ldots,0,k-1$ in the second. In each case the column
order and leading signs $(-1)^{j-1}$ give positive Vandermonde sign.
The triangular flag criterion makes $F$ TN, so Cauchy--Binet, using
the term with $F[I,I]$ for row set $I$, gives positive initial-column
minors of $FQ$ and positive maximal minors of
$F[Q_1,\ldots,Q_{k-1},w]$. In the factorization $FQ=L_0U_0$,
the former minors give $L_0$ TN by the flag criterion, and the
latter give $v=L_0^{-1}Fw>0$. For its $i$th entry, take the
determinant on the first $i$ rows of $[(FQ)[1:r,1:i-1],Fw]$ and
divide by the leading $(i-1)$-minor of $FQ$.

Put $R=-L_0^{-1}F$, $C_0=R\Lambda R^{-1}$, $u=gR^{-1}$ and $B=vu-C_0$.
Lagrange interpolation gives $gQ_1=-1$ and $gQ_j=0$ for $j>1$:
$gQ_j$ is $(-1)^j$ times the divided difference of
$\prod_{h=j+1}^r(x+d_h)$ at the $\eta_s$. Since
$(-\Lambda-d_1I)Q_1=-w$ and $(-\Lambda-d_jI)Q_j=Q_{j-1}$ for $j>1$,
we get $(-\Lambda-wg)Q=QJ_d$, where $J_d$ has diagonal $d_j$ and unit
superdiagonal. Hence $B=U_0J_dU_0^{-1}$ is upper triangular with
diagonal $d_j$.

With $D_v=\diag(v)$, use tildes for conjugation by $D_v^{-1}$ and put
$\widetilde u=uD_v$. Then $\widetilde C_0=\mathbf1\widetilde u-\widetilde B$.
Its consecutive row differences give an upper Hessenberg matrix
$C_*=(I-N_-)\widetilde C_0$ with subdiagonal $d_1,\ldots,d_{r-1}$,
where $N_-$ is the lower shift. The upper LU factor $U_c$ of the
nonsingular TN matrix $HL_0$ is TN with positive pivots $\gamma_j$.
Since $C_0=L_0^{-1}(HL_0)$, the upper LU factor of $C_*$ is
$D_v^{-1}U_cD_v$. Its unit lower factor is bidiagonal, with
subdiagonal $d_j/\gamma_j$, so $C_*$ is TN.

For the upper shift $N_+$ of order $r+1$, row differences also give
\[
 \begin{pmatrix}\widetilde u&1\\\widetilde B&\mathbf1\end{pmatrix}
 =(I-N_+)^{-1}\diag(C_*,1)(I+d_rE_{r+1,r}).
\]
All factors are TN, with $(I-N_+)^{-1}$ the product of the
$I+E_{i,i+1}$ in increasing order. Thus
$U=\left(\begin{smallmatrix}a&u&1\\0&B&v\\0&0&c\end{smallmatrix}\right)$
is TN by undoing the positive diagonal scalings and expanding
along its first column or last row when selected.

Back substitution in the last resolvent column gives
$E_U(z)=\det(zI-B+vu)=V(z)$ by the determinant lemma. Dividing
the internal entries of that column by its first entry gives
$(z-a)(zI+C_0)^{-1}v$. As $v=-Rw$, this equals
$v+R((z+\eta_s)^{-1})_{s=1}^r$, which is \eqref{strealizedflag}.
\end{proof}

\subsection{Mixed evaluation determinants}

Let $n\ge3$, $r=n-2$, and $p+q=n$ with $p,q\ge1$. For a matrix
$B$ put $E_B(z)=\adj(zI-B)_{1n}$ and, where $E_B(z)\ne0$,
write $\psi_j^B(z)=\adj(zI-B)_{jn}/E_B(z)$. In this notation the
determinant \eqref{sourcemixed} is
\begin{equation}\label{stmixed}
 D_B^p[f](x_1,\ldots,x_n)
   =\det[\psi_1^B,\ldots,\psi_q^B\mid
          f\psi_1^B,\ldots,f\psi_p^B](x_i).
\end{equation}
If $B$ has increasing positive eigenvalues $\lambda_1<\cdots<\lambda_n$ and
$E_B(\lambda_i)\ne0$, the adjugate identity preceding
\eqref{sourcecorner} gives $\psi_j^B(\lambda_i)=X_{ji}$ for the right
eigenvector matrix normalized by $X_{1i}=1$. The same column
expansion as in \eqref{sourcecorner} therefore gives
\begin{equation}\label{stspectral}
 D_B^p[f](\lambda_1,\ldots,\lambda_n)
   =\det(X)\det f(B)[1:p,q+1:n].
\end{equation}
The required nonzero cofactors at eigenvalues hold in particular
for oscillatory matrices, by the endpoint eigenvector property.

For $A\in\TP_n$, the pole and residue formulas
\eqref{shared}--\eqref{actualflag} give
\begin{equation}\label{stactualflag}
 \begin{aligned}
 E_A(z)&=b\prod_{a=1}^r(z+\mu_a),\qquad b>0,\\
 (\psi_2^A,\ldots,\psi_{n-1}^A)^T
   &=u-W^T\left((z+\mu_a)^{-1}\right)_{a=1}^r .
 \end{aligned}
\end{equation}
Put $\Lambda=\diag(\mu_1,\ldots,\mu_r)$. By the residue formula
following \eqref{actualflag}, $W^T\Lambda(W^T)^{-1}=C^T\in\TP_r$.
Apply Lemma~\ref{sttriangular} to $C^T$ and put $F=K^{-1}W^T$.
Then $F\Lambda F^{-1}$ is upper triangular TN with diagonal
$(\mu_1,\ldots,\mu_r)$. Back substitution makes its eigenvector
matrix $F$ upper triangular. Since $K$ is unit lower
triangular, $F$ retains the positive initial-row flag minors of
$W^T$. In particular, the successive leading principal determinants
give its positive diagonal entries.

For fixed increasing positive nodes $x_1,\ldots,x_n$, let $U$ be
the matrix given by Lemma~\ref{strealize} with $F$ and these nodes
as its diagonal entries in this order. If
$\phi=F((z+\mu_a)^{-1})_{a=1}^r$, the internal ratio vectors for
$A$ and $U$ are $u-K\phi$ and $v-L_0^{-1}\phi$. In both blocks
of \eqref{stmixed}, subtract the constants using the first column.
The corresponding initial blocks of $K^{-1}$ and $L_0$ then reduce
both systems to the same functions by unit triangular operations.
Hence
\begin{equation}\label{stsameD}
 D_A^p[f](x_1,\ldots,x_n)=D_U^p[f](x_1,\ldots,x_n).
\end{equation}
For $U$, all endpoint cofactors $E_U(x_i)$ are positive by
Lemma~\ref{strealize}. An eigenvector for $x_j$ with $j$th
coordinate $1$ has zero coordinates after $j$ and nonnegative
coordinates before $j$ by back substitution, since $U\ge0$ and
$x_j-x_i>0$ for $i<j$. Its first coordinate is nonzero by
$E_U(x_j)>0$. Thus the eigenvector matrix in \eqref{stspectral}
is upper triangular with positive diagonal. If $f$ preserves
nonsingular TN matrices, \eqref{stsameD} and \eqref{stspectral}
now give
\begin{equation}\label{stweaksign}
 D_A^p[f](x_1,\ldots,x_n)\ge0
 \qquad(0<x_1<\cdots<x_n).
\end{equation}

\begin{lemma}\label{stanalytic}
Let $h_1,\ldots,h_d$ be linearly independent real analytic functions
on an open interval. Suppose that their span contains a positive
function and that every evaluation determinant at $d$ increasing
nodes is nonnegative. Then every such determinant is positive.
\end{lemma}

\begin{proof}
A function in the span cannot have $d$ strict sign changes.
Otherwise choose $d+1$ ordered points with alternating nonzero
values. A small perturbation of the first $d$ points can make
their evaluation determinant nonzero while preserving the signs.
Indeed, independence makes this determinant a nonzero analytic
function of the points. The alternating cofactors of the
$(d+1)\times d$ evaluation matrix annihilate its column space.
Their pairing with the alternating values has terms of one sign,
with at least one nonzero term, a contradiction.

If a nonzero function $h$ in the span has $d$ distinct zeros,
isolate them in disjoint intervals with nonzero endpoint values.
For positive $g$ in the span, both perturbations $h\pm\varepsilon g$
retain a sign change at each odd-order zero when $\varepsilon$
is small. At an even-order zero one of the perturbations creates
two sign changes. Thus one perturbation has at least $d$ strict
sign changes, a contradiction. A singular evaluation matrix
would supply such an $h$.
\end{proof}

\subsection{Uniform corner comparisons}

\begin{proposition}\label{stuniform}
Let $n\ge2$, and let $f:(0,\infty)\to(0,\infty)$ be real analytic
and preserve nonsingular $\TN_n$. Suppose that $f$ is not $cx^k$ for
any $c>0$ and $k\in\{0,\ldots,n-2\}$. Let $K$ be a compact subset
of $\{0<\lambda_1<\cdots<\lambda_n\}$. There are constants
$0<c\le C<\infty$, depending only on $n,f,K$, such that every
nonsingular TN matrix $A$ with ordered spectrum in $K$ satisfies
\[
 cG_p(A;n-1)\le G_p(f(A))\le CG_p(A;n-1)
 \qquad(1\le p<n).
\]
The same bounds hold with $F_p$ in place of $G_p$.
\end{proposition}

\begin{proof}
For $n=2$, the upper bidiagonal test $xI+sE_{12}$ gives
$f'(x)\ge0$, so nonconstancy and analyticity make $f$ strictly
increasing. Two-point interpolation multiplies both off-diagonal
entries of $A$ by $f[\lambda_1,\lambda_2]$, which is continuous
and strictly positive on $K$.

Let $n\ge3$, $r=n-2$, and $g(x)=x^{n-1}$. For $A\in\TP_n$, the
polynomials $P_j(x)=\adj(xI-A)_{jn}$, $1\le j<n$, have degree
at most $r$ and are independent: at the eigenvalues,
$P_j(\lambda_i)=E_A(\lambda_i)X_{ji}$ by \eqref{stspectral} and
its eigenvector normalization. Orthonormalize their coefficient
vectors in order, taking positive pivots, and keep the resulting
column polynomials $u_1,\ldots,u_{n-1}$. Their coefficient matrices
range in the orthogonal group; denote their compact closure by
$\mathcal F$. For $U\in\mathcal F$ and $p+q=n$, write
\[
 D_h(U;x)=\det[u_1,\ldots,u_q,hu_1,\ldots,hu_p](x_i).
\]
Before taking the closure, clearing the common row denominators
$E_A(x_i)$ and changing basis within both prefixes multiplies
\eqref{stmixed} by a positive factor independent of $h$.
Hence \eqref{stweaksign} gives $D_h(U;x)\ge0$ at all increasing
positive nodes for $h=f,g$, also on the closure.

For every limiting basis, $u_1$ has norm one and nonnegative
coefficients, inherited from $P_1=E_A$ in \eqref{stactualflag}.
In particular $u_1(x)>0$ for $x>0$. A dependence in a mixed
family would give $P+hQ=0$ with nonzero polynomials of degree
at most $r$, which is impossible by degree for $h=g$.
For $h=f$, Proposition~\ref{rigidity} gives $f(x)=O(x^r)$ at
zero and $f(x)\ge f(1)x^r$ at infinity. In a reduced quotient
with both degrees at most $r$, the first bound forces the numerator
to be a multiple of $x^r$, and the second forces the denominator
to be constant. This would give the excluded $f(x)=cx^r$.
The positive span function $u_1$ now allows Lemma~\ref{stanalytic}
to give $D_f,D_g>0$ at every strictly ordered positive node tuple.

Choose a compact neighborhood $K'$ of $K$ inside that region.
The continuous ratio $D_f/D_g$ on $\mathcal F\times K'$ has a
positive minimum and a finite maximum. At the eigenvalues of a
TP input it equals $G_p(f(A))/G_p(g(A))$ by \eqref{stspectral},
since the basis factors cancel. Take common bounds for the finitely
many $p$. Any TN input with spectrum in $K$ has TP approximants
with spectra eventually in $K'$. Continuity gives the inequalities
even when the reference corner vanishes. Transposition gives them
for $F_p$ with the same constants.
\end{proof}

\begin{proof}[Proof of Theorem~\ref{strict}]
Suppose $f$ preserves nonsingular $\TN_n$ and is outside the
stated monomial exceptions. For $A\in\Osc_n$, the simple spectrum
allows Proposition~\ref{stuniform} to compare the far corners of
$f(A)$ with those of $A^{n-1}$. The latter are positive by the
oscillation theorem. Since $\det f(A)>0$, Lemma~\ref{recognition}
makes $f(A)$ TP.

Conversely, the matrix $J_n(a)=\operatorname{tridiag}(1,a,1)$
used in \eqref{sharp} is oscillatory for $a>2$, also in order two.
The bandwidth of $J_n(a)^k$ is at most $k$, so
$(cJ_n(a)^k)_{1n}=0$ whenever $0\le k\le n-2$.
\end{proof}

There is a direct test for the two extreme corner orders when
eigenvalues coincide. Let $B$ be a nonsingular TN matrix of order
$d\le n$, $d\ge2$, with positive superdiagonal, and let
$E(z)=\adj(zI-B)_{1d}$. Its coefficients are nonnegative by
\eqref{shared} and TP approximation (for $d=2$, $E=B_{12}$).
It is nonzero because the resolvent expansion at infinity contains
$(B^{d-1})_{1d}>0$. At the eigenvalues $\lambda_1,\ldots,\lambda_d$
of $B$, residues of the resolvent, followed by confluence, give
\[
 f(B)_{1d}=[fE]_{\lambda_1,\ldots,\lambda_d}.
\]
Here and below $f$ satisfies the hypotheses of
Proposition~\ref{stuniform}, as does $f^\#$ by \eqref{checkerdual}.
The shift test in Proposition~\ref{rigidity} gives
$f^{(j)}\ge0$ for $j\le n-1$. The last derivative is not
identically zero: its vanishing would make $f$ a polynomial of degree
at most $n-2$, hence an excluded monomial by the growth bounds
(or a constant for $n=2$). Each earlier derivative is therefore
positive by monotonicity and analyticity.

By Leibniz' rule, $(fE)^{(d-1)}\ge E f^{(d-1)}$ on the positive
axis. By \eqref{hgdivided}, the divided difference is positive if $d<n$, or if the nodes are not all
equal: in the latter case their convex combinations meet an open
interval of positivity. If $d=n$ and all nodes equal $\lambda$, the derivative is positive
exactly when $f^{(d-1)}(\lambda)>0$ or $E$ is nonconstant: in the
latter case $E'(\lambda)f^{(d-2)}(\lambda)>0$. The same argument for $f^\#$ and
$\iota(B)$ gives the order $d-1$ upper corner, since
$(f^\#(\iota(B)))_{1d}=G_{d-1}(f(B))/\det f(B)$ by Jacobi's identity.
Indeed, each upper zero cut is invariant under inversion by block
triangularity and hence under $\iota$ by Lemma~\ref{zerocuts}.

\begin{corollary}\label{stzeropattern}
Let $n\ge2$, and let $f:(0,\infty)\to(0,\infty)$ be real analytic
and preserve nonsingular $\TN_n$. Suppose that
$f(x)\not\equiv cx^k$ for $c>0$ and $k=0,\ldots,n-2$.
Let $A\in\TN_n$ be nonsingular. Assume that each directional block
of order $d\ge4$ has simple spectrum. For a block of order $d=2,3$
with a single eigenvalue $\lambda$, assume
$f^{(d-1)}(\lambda)>0$ and $(f^\#)^{(d-1)}(1/\lambda)>0$;
these conditions are automatic if $d<n$.
Then for every pair of equal-size index sets $I,J$,
\[
 \det f(A)[I,J]>0
 \quad\Longleftrightarrow\quad
 \det A^{n-1}[I,J]>0.
\]
These determinants are positive exactly when \eqref{zeropath}
holds, and otherwise both are zero.
\end{corollary}

\begin{proof}
Embedding a nonsingular $B\in\TN_q$ as $B\oplus I_{n-q}$ gives
preservation in every smaller order. In an upper directional block
$B=A[C,C]$ of order $q\ge2$, its superdiagonal is positive, so
Lemma~\ref{directioncorners} gives all $G_p(B^{q-1})>0$.
For $q\ge4$, Proposition~\ref{stuniform} transfers strict positivity
to $G_p(f(B))$. For $q=2,3$, all proper corners are extreme, so the
preceding divided-difference test applies. Lower blocks follow by
transposition, and scalar blocks are positive. Thus
Proposition~\ref{localreduction} gives \eqref{zeropath} for
$f(A)$, and Theorem~\ref{zeropattern} gives it for $A^{n-1}$.
\end{proof}

Conversely, if this zero-pattern assertion holds for every input
with simple directional block spectra, it holds for all
oscillatory inputs. Their $(n-1)$st powers are TP, so
Theorem~\ref{strict} excludes precisely the stated monomials.

\begin{remark}
In order $d=2$ or $3$, under the same hypotheses on $f$,
the zero pattern agrees with that of $A^{d-1}$ for every nonsingular
TN input if and only if both $f^{(d-1)}$ and $(f^\#)^{(d-1)}$ are
positive on $(0,\infty)$. For necessity, take $J=\lambda I+N$ with
$N$ the upper shift: $f(J)_{1d}=f^{(d-1)}(\lambda)/(d-1)!$ and
$(J^{d-1})_{1d}=1$. Testing $\iota(J)$ and taking complementary
minors gives the same test for $f^\#$. For $d=2$ the two derivative
conditions are equivalent by \eqref{checkerdual}. For example,
$f(x)=\exp((x-1)^3)$ preserves nonsingular $\TN_2$: its divided
differences are nonnegative, its diagonal entries after evaluation
are positive by Lemma~\ref{principalbounds}, and the determinant is
positive. Yet $f(I+sE_{12})=I$ for $s>0$, since $f'(1)=0$.
\end{remark}

\begin{remark}\label{repeatedfour}
In order four the middle corner can vanish even if $f'''$ and $(f^\#)'''$
are positive everywhere. For $0\le s\le1$ set $f_s(x)=e^{F_s(\log x)}$, where
$F_s(y)=(4+s/60)y-(s/1800)\tanh(30y)$. Since $F_s$ is odd, $f_s^\#=f_s$.
For $(p,q,d)=(F_s',F_s'',F_s''')$ evaluated at $\log x$, differentiation
gives $4\le p\le4+1/60$, $|q|\le1$, and $-10\le d\le30$, with
$(p,q,d)=(4,0,30)$ at $s=x=1$. We have $xf_s'/f_s=p$,
$x^2 f_s''/f_s=p(p-1)+q>0$, and
$x^3 f_s'''/f_s=p(p-1)(p-2)+3(p-1)q+d\ge24-3(3+1/60)-10>0$.
The Schwarzian $\mathcal S(f)=f'''/f'-\tfrac32(f''/f')^2$ satisfies
\[
 x^2\mathcal S(f_s)=\frac{1-p^2}{2}+\frac{d}{p}-\frac{3q^2}{2p^2}\le0,
 \qquad \mathcal S(f_1)(1)=0.
\]
For $A\in\TP_4$, $G_1(f_s(A))=[f_sE_A]_\lambda>0$ by the derivative bounds
and the coefficients of $E_A$, and \eqref{checkerdual} gives $G_3>0$. By
\eqref{actualflag} and \eqref{stspectral} the middle corner uses $(1,r,f_s,rf_s)$ with
$r=u-\sum_{i=1}^2w_i/(x+\mu_i)$, $w_i>0$ and distinct $\mu_i$. Its nontrivial
initial Wronskians are $r'$, $r'f_s''-r''f_s'$ and
$2(r'f_s')^2(\mathcal S(r)-\mathcal S(f_s))$, all positive since $r'>0>r''$
and $\mathcal S(r)>0$ by strict Cauchy--Schwarz. Lemma~\ref{ECT} and
transposition give both middle corners. Starting at $f_0(A)=A^4$, these
strict corners and the positive determinant exclude a first exit from
$\TP_4$ by Lemma~\ref{recognition}. Approximation gives preservation of
nonsingular $\TN_4$. For the upper shift $N$, put $J=I_4+N$. The confluent
extreme tests give $G_1(f_1(J)),G_3(f_1(J))>0$, whereas Taylor expansion
gives $G_2(f_1(J))=-f_1'(1)^2\mathcal S(f_1)(1)/6=0<G_2(J^3)$.
\end{remark}

\subsection{Preservation in all orders}

For symmetric matrices with positive spectrum, orthogonal
diagonalization defines $f(A)$ without a regularity assumption on
$f$. If $f$ is positive and analytic, Proposition~\ref{rigidity}
already forces every all-order TN preserver to be a monomial,
by choosing $n-2>f'(1)/f(1)$. The symmetric case also allows
a smaller family of test matrices. Here a Jacobi matrix is a real
symmetric tridiagonal matrix with positive adjacent entries.
Schwarz's theorem gives $e^{tJ}\in\TP$ for every such matrix $J$
and every $t>0$ \cite[Theorem~3]{SchwarzTPDS}.

For a real sequence $(a_j)_{j\ge0}$, write $Da_j=a_j-a_{j+1}$.
Hausdorff's moment theorem \cite{HausdorffMom} states that if
$D^k a_j\ge0$ for all $j,k\ge0$, then
$a_j=\int_0^1 u^j\,\mathrm d\mu(u)$ for a finite nonnegative
measure $\mu$. Applying Cauchy--Schwarz to $1$ and $u$ gives
$a_1^2\le a_0a_2$.

\begin{theorem}\label{heat-rigidity}
Let $f:(0,\infty)\to(0,\infty)$. The matrices $f(e^J)$ are TN for
every Jacobi matrix $J$ of every order if and only if $f(x)=cx^a$
with $c>0$ and $a\ge0$. The matrices $f(A)$ are TN for every real
symmetric TP matrix $A$ of every order if and only if $f(x)=cx^k$
with $c>0$ and $k$ a nonnegative integer.
\end{theorem}

\begin{proof}
Put $F(y)=f(e^y)$ and $G(y)=1/F(-y)$. Every set of $n$ distinct
real numbers is the spectrum of a Jacobi matrix. Use multiplication
by the coordinate in the orthonormal polynomial basis for a measure
with positive weights on these $n$ points. The tridiagonal form
follows from orthogonality, and positive leading coefficients give
positive adjacent entries.
For a Jacobi matrix with eigenvalues $\lambda_1,\ldots,\lambda_n$,
polynomial interpolation gives
\[
 (F(J))_{1n}=\left(\prod_{i=1}^{n-1}J_{i,i+1}\right)
                  F[\lambda_1,\ldots,\lambda_n],
\]
because reaching $n$ from $1$ requires $n-1$ adjacent steps,
whose product is the displayed coefficient of the divided difference.
Hence every divided difference of $F$ at distinct nodes is nonnegative.
Writing $E=\diag(1,-1,1,\ldots)$, the matrix $E(-J)E$ is Jacobi and
$G(J)=E F(E(-J)E)^{-1}E$. Jacobi's complementary-minor identity gives
the same divided-difference inequality for $G$. In fact, for
$H=f(e^J)$, the two tests $H_{1n}\ge0$ and
$\det H[1:n-1,2:n]\ge0$ in each order $n\ge2$ suffice here. The
second supplies the far entry of the checkerboard inverse, since
$\det H>0$.

For $y\in\mathbb R$ and $h>0$ put $a_j=F(y-jh)$. At $x=y-jh$,
$D^k a_j=k!h^k[F]_{x,x-h,\ldots,x-kh}\ge0$ for $k\ge1$, and
$a_j\ge0$. The moment inequality gives $F(y-h)^2\le F(y)F(y-2h)$.
The same inequality for $G$ at the reflected points gives the reverse
inequality for $F$.
Thus $\log F$ is midpoint affine. It is nondecreasing by the
first divided differences, so it is affine with nonnegative slope.
This gives $f(x)=cx^a$ with $a\ge0$. Conversely, $f(e^J)=ce^{aJ}$
is TN, including the case $a=0$.

Since $e^J$ is symmetric and TP, the second assertion also forces
$f(x)=cx^a$ with $a\ge0$. If $a$ is nonintegral, take
$n=\lfloor a\rfloor+3$ and the tridiagonal matrix $B=J_n(b)$,
$b>2$, of \eqref{sharp}. Then $(B^a)_{1n}<0$.
For $\varepsilon>0$, the matrix $B_\varepsilon=e^{\varepsilon B}B$
is symmetric by commutation and TP by Cauchy--Binet. In each
minor the TP left factor and the nonsingular TN right factor
give a positive term.
Its $a$th power retains the negative entry for small $\varepsilon$.
Nonnegative integer powers preserve TN by Cauchy--Binet.
\end{proof}

If $f$ is real valued, the zero function is the only additional
case in either assertion. A real symmetric TN matrix is positive
semidefinite, so the freely prescribed Jacobi spectra give $F\ge0$.
The midpoint inequality for $F$ still holds and a zero at $y_0$
forces $F(z)^2\le F(y_0)F(2z-y_0)=0$ for every $z$.

\section{Fixed multipliers and finite spectral interpolation}

\subsection{Exponent conditions}

Theorem~\ref{multipliers} supplies sufficient exponent bounds;
Proposition~\ref{rigidity} makes them necessary for the following
fixed multipliers.

\begin{corollary}\label{fixedbuffer}
Fix $n\ge2$ and a nonconstant real operator-monotone function $\psi$ on $(0,\infty)$. Suppose that $\psi(0+)\in\mathbb R$ or $\psi(+\infty)\in\mathbb R$. For $\alpha\in\mathbb R$, the following conditions are equivalent:
\begin{enumerate}
\item $A^\alpha e^{\psi(A)}$ is TN for every nonsingular $A\in\TN_n$;
\item $A^\alpha e^{\psi(A)}$ is TN for every $A\in\TP_n$;
\item $A^\alpha e^{\psi(A)}$ is TP for every $A\in\TP_n$;
\item $A^\alpha e^{\psi(A)}$ is TP for every $A\in\Osc_n$;
\item $\alpha\ge n-2$.
\end{enumerate}
\end{corollary}

\begin{proof}
Theorem~\ref{multipliers} proves all preservation assertions when $\alpha\ge n-2$. Conversely, at an endpoint where $\psi$ has a finite real limit, $f(x)=x^\alpha e^{\psi(x)}\sim cx^\alpha$ for some $c>0$. If condition~(2) holds and $\alpha<n-2$, Proposition~\ref{rigidity} forces $\alpha=k$ and $f(x)=cx^k$. Hence $e^{\psi(x)}=c$ on $(0,\infty)$, contrary to nonconstancy. Since each preservation condition implies condition~(2), all five are equivalent.
\end{proof}

\begin{remark}
The finite-endpoint hypothesis cannot be omitted: for $\psi(x)=\log x$, the function is the pure power $x^{\alpha+1}$ and has the exceptional exponent set of Corollary~\ref{sets}, shifted by one. Nonconstancy is also essential, since a constant $\psi$ retains every nonnegative integer exponent.
\end{remark}

For each fixed $s>0$, Corollary~\ref{fixedbuffer} gives, in particular,
\[
 A^\alpha e^{sA}\in\TN_n\text{ for all nonsingular }A\in\TN_n
 \quad\Longleftrightarrow\quad \alpha\ge n-2.
\]

\begin{corollary}\label{shifted}
Fix $n\ge2$ and $b>0$. The function $f(x)=x^\alpha(x+b)^{-\beta}$ preserves all nonsingular TN matrices of order $n$ if and only if
\begin{equation}\label{shiftset}
 (\alpha,\beta)\in
 \{\alpha\ge n-2,\ \alpha-\beta\ge n-2\}
 \ \cup\ \{(k,0):k=0,\ldots,n-3\}.
\end{equation}
The second set is empty when $n=2$. If $\beta\ne0$, the same condition is necessary and sufficient for $f(A)$ to be TP for every $A\in\TP_n$, and also for every $A\in\Osc_n$.
\end{corollary}

For $\beta\ne0$, the condition is $\alpha\ge n-2+\max\{\beta,0\}$. For doubly nonnegative matrices, Guillot, Khare and Rajaratnam~\cite[Theorem~1.6]{GKR} obtained this tail threshold together with a finite set of possible exceptions. Under total nonnegativity, every exponent below this threshold is excluded when $\beta\ne0$.

\begin{proof}
When $\beta\ne0$, the function $f$ is not a monomial, and its endpoint asymptotics are
\[
 f(x)\sim b^{-\beta}x^\alpha\quad(x\downarrow0),\qquad
 f(x)\sim x^{\alpha-\beta}\quad(x\to\infty).
\]
Proposition~\ref{rigidity} forces both $\alpha\ge n-2$ and $\alpha-\beta\ge n-2$, even if we require only that $f(A)$ be TN for every $A\in\TP_n$.

For sufficiency, write $f(x)=x^{n-2}e^{\psi(x)}$. If $\beta=-q<0$, take
\[
 \psi(x)=(\alpha-n+2)\log x+q\log(x+b).
\]
If $\beta>0$, take
\[
 \psi(x)=(\alpha-n+2-\beta)\log x
       +\beta\log\frac{x}{x+b}.
\]
In each case the coefficient of $\log x$ is nonnegative. The functions $\log x$, $\log(x+b)$, and $\log(x/(x+b))$ are operator monotone. In the noncompact form following \eqref{Pick}, their representing measures are Lebesgue measure on $[0,\infty)$, $[b,\infty)$, and $[0,b]$, respectively, with the appropriate real constants at $x=1$. Thus each displayed $\psi$ is nonconstant and operator monotone. Theorem~\ref{multipliers} gives TN preservation and TP preservation on oscillatory matrices throughout the closed parameter range. When $\beta=0$, Corollary~\ref{sets} gives the exact pure-power assertion.
\end{proof}

For example, Theorem~\ref{multipliers} gives
$A^{n-2}(A+aI)\in\TP_n$ for $a>0$ and
$A^{n-2}\log(I+A)\in\TP_n$ whenever $A\in\Osc_n$.
The first assertion uses $\psi(x)=\log(x+a)$; the second uses
$\psi(x)=\log\log(1+x)$, which is operator monotone and nonconstant
because $\log(1+x)$ is positive and operator monotone.

For nonsingular $A\in\TN_n$ with $d(A)\le d$, $2\le d\le n$,
Proposition~\ref{localreduction} gives the exact TN-preserving ranges
in Corollaries~\ref{fixedbuffer} and~\ref{shifted} with $d$ in place of $n$.

\subsection{Finite spectral applications}

Corollary~\ref{csfinite} uses a separate interpolant for each
spectral subset in the Gram expansion, so it places no upper bound
on the ambient matrix order. For $n\le N$, a single interpolant
at all the eigenvalues of the given matrix gives the following
consequence of Theorems~\ref{multipliers} and~\ref{iminternal}.

\begin{corollary}\label{finitePick}
Let $N\ge2$, and let $\psi$ be real analytic on $(0,\infty)$.
Suppose that, on every set of $N$ distinct positive numbers, $\psi$
agrees with some real operator-monotone function on $(0,\infty)$.
For $2\le n\le N$ and $\tau\ge0$, set $\Phi(x)=x^\alpha e^{\tau\psi(x)}$.
\begin{enumerate}
\item If $\alpha\ge n-2$, then $\Phi(A)$ is TN for every
nonsingular $A\in\TN_n$.  For $A\in\Osc_n$, the matrix $\Phi(A)$ is TP
if $\alpha>n-2$, or if $\tau>0$ and $\psi$ is nonconstant.
If $\Phi$ extends analytically across zero, the TN assertion
also includes singular TN matrices.
\item If $n\ge4$ and $\alpha\ge n-3$, then $\Phi(A)$ is a
$P$-matrix for every nonsingular $A\in\TN_n$. All its minors
other than its $(1,n)$ and $(n,1)$ entries and its two far corners
of order $n-1$ are nonnegative, and are strictly positive if $A$ is TP.
\end{enumerate}
\end{corollary}

\begin{proof}
For a multiset of $2\le n\le N$ nodes in $(0,\infty)$, there is an
operator-monotone interpolant $\chi$ matching $\psi$ to the
multiplicities of those nodes. To see this, approximate the multiset
by distinct nodes in a compact interval $[a,b]\subset(0,\infty)$
and take the given interpolants, padding to $N$ nodes if needed.
For two of the approximating nodes $u,v$, the secant slope of
each interpolant equals $[\psi]_{u,v}$, bounded on $[a,b]$ by
$\max_{[a,b]}|\psi'|$. In \eqref{Pick} this slope is the integral of
$(1+s)^2/((u+s)(v+s))$, uniformly bounded below by a positive
constant for $u,v\in[a,b]$. Thus the representing measures have
bounded mass, and the values at $1$ are bounded by interpolation
and \eqref{Pick}. Take a weakly convergent subsequence of the measures
on $[0,\infty]$ and a convergent subsequence of the constants.
The bounded masses and joint continuity of the kernels and their
$x$ derivatives, including at $s=\infty$, give convergence locally
uniform in $x>0$, with all derivatives. The resulting operator-monotone
limit matches the divided differences on each cluster by
\eqref{hgdivided}, hence the required derivatives.

Take the multiset of eigenvalues of a nonsingular TN matrix $A$.
Hermite interpolation now gives $\chi(A)=\psi(A)$, including on
Jordan blocks. Theorems~\ref{multipliers} and~\ref{iminternal}
apply with multiplier $\tau\chi$. For oscillatory input, the nodes
are distinct, and nonconstant $\psi$ is strictly increasing by
the two-point consequence of the hypothesis, as in
Corollary~\ref{csfinite}. Thus $\tau\chi$ is nonconstant if $\tau>0$.
When $\Phi$ is analytic across zero, TN approximation includes
singular inputs.
\end{proof}

For nonsingular TN inputs of arbitrary order with $2\le d(A)\le N$
and $\alpha\ge d(A)-2$, part~(1) gives total nonnegativity on each
nonscalar directional block. Scalar blocks are positive, so
Proposition~\ref{localreduction} gives $\Phi(A)\in\TN$.

A function $h$ is matrix monotone of order $q$ on an interval if
$X\le Y$ implies $h(X)\le h(Y)$ for Hermitian $q\times q$ matrices
with spectra in that interval.

Theorem~3.1 of Ameur, Kaijser and Silvestrov with $n=q-1$ sends the
positive order-$q$ matrix-monotone class $P'_q$ into $C_{2q-1}$.
Corollary~2.4 identifies membership in $C_{2q-1}$ with matching by
a positive Pick function on each such set of nodes. These Pick
functions are positive operator-monotone functions by L\"owner's
theorem \cite[Corollary~2.4, Theorem~3.1 and (2.1)--(2.2)]{AKS}.
A smaller set can be enlarged to $2q-1$ distinct positive nodes.

\begin{lemma}[Finite matching]\label{finitepositiveinterpolation}
Let $q\ge2$. If $h:(0,\infty)\to(0,\infty)$ is matrix monotone of order
$q$, then on any set of at most $2q-1$ distinct positive nodes it agrees
with a positive operator-monotone function on $(0,\infty)$.
\end{lemma}

\begin{example}
For each $N\ge2$, choose $q\ge2$ with $2q-1\ge N$. In the polynomial
construction of Hansen, Ji and Tomiyama \cite[Sec.~2]{HJT},
$g_q(t)=\sum_{j=0}^{q-1}t^{2j+1}/(2j+1)$ is matrix monotone of order
$q$ on $[0,\delta)$ for some $\delta>0$ and fails matrix monotonicity of
order $q+1$ on every nondegenerate subinterval of $[0,\delta)$.
For $0<\beta<\delta$, set $\varphi_\beta(x)=\beta x/(1+x)$ and
$\psi=g_q\circ\varphi_\beta$. The bijection $\varphi_\beta$ from
$(0,\infty)$ to $(0,\beta)$ and its inverse
$\varphi_\beta^{-1}(t)=\beta/(\beta-t)-1$ are operator monotone on
their respective intervals, by order reversal under inversion.
Thus $\psi$ is positive and analytic on $(0,\infty)$ and matrix monotone
of order $q$ there. If $\psi$ were operator monotone, then
$g_q=\psi\circ\varphi_\beta^{-1}$ would be operator monotone on
$(0,\beta)$, contrary to its failure of order $q+1$. By
Lemma~\ref{finitepositiveinterpolation}, $\psi$ satisfies the hypothesis
of Corollary~\ref{finitePick}.
\end{example}

\section*{Concluding remarks}\addcontentsline{toc}{section}{Concluding remarks}
For ordinary powers, the uniform threshold $n-2$ for positivity
of minors of all orders remains sharp on every prescribed simple
positive spectrum, even for symmetric TP inputs
(Corollary~\ref{fixedpower}). Fixing the order of the minors
can give a smaller threshold.
For endpoint flag minors of order $p$, Theorem~\ref{endpointflags}
gives the sharp threshold $\gamma_{n,p}=\max(p,n-p)-1$.
For all minors of that order, this is a lower bound for $g_{n,p}$.
Corollary~\ref{imexact} supplies the upper bound $n-3$ for
$2\le p\le n-2$. The bounds agree for $p=2,n-2$, while
Theorem~\ref{main} gives equality with $\gamma_{n,p}$ for $p=1,n-1$.
For $3\le p\le n-3$, Section~\ref{endpointsec} reduces equality to
positivity for pairs of nonendpoint cyclic intervals on
$(\gamma_{n,p},\gamma_{n,p}+1)$.

\begin{conjecture}[Fixed-order critical exponent]\label{fixedorderconj}
For $n\ge2$ and $1\le p<n$, the sharp uniform tail threshold for
positivity of all $p$th-order minors of $A^t$ over $A\in\TP_n$ is
$g_{n,p}=\max(p,n-p)-1$.
\end{conjecture}

\end{document}